\documentclass[11pt,leqno]{amsart} 
\usepackage[margin=1.2in]{geometry} 
\usepackage{kkn} 
  
\title[Syntomic cohomology of Morava K-theory]{Syntomic cohomology of Morava K-theory}

\author[G. Angelini-Knoll, J. Hahn, and D. Wilson]{Gabriel Angelini-Knoll, Jeremy Hahn, and Dylan Wilson}
\address{Universit{\'e} Paris 13, LAGA, CNRS, UMR 7539, F-93430, Villetaneuse, France}
\email{angelini-knoll@math.univ-paris13.fr}

\address{Department of Mathematics, Massachusetts Institute of Technology, Cambridge, MA, USA}
\email{jhahn01@mit.edu}

\address{K1X, Inc. Technology, Information and Internet, Morristown, NJ}
\email{dylwil3@gmail.com}

\usepackage{graphicx}   
\graphicspath{{./EU-symbol.jpg}}

\usepackage{hyperref}
\hypersetup{
    hypertexnames=false,
    colorlinks=true,
    linkcolor=blue,
    filecolor=magenta,      
    urlcolor=cyan,
    citecolor=seagreen, 
    pdftitle={Syntomic cohomology of Morava K-theory},
    pdfpagemode=FullScreen,
}

\begin{document} 

\begin{abstract}     
We compute the $\MU$-based syntomic cohomologies, mod $(p,v_1,\cdots,v_{n+1})$, of all $\bE_{1}$-$\MU$-algebra forms of connective Morava K-theory $k(n)$. As qualitative consequences, we deduce the Lichtenbaum--Quillen conjecture, telescope conjecture, and redshift conjecture for the algebraic K-theories of all $\bE_{1}$-$\bS$-algebra forms of $(2p^n-2)$-periodic Morava $K$-theory. Notably, the motivic spectral sequence computing $\pi_*\TC(k(n))_p^{\wedge}$ is concentrated on at most three lines, independently of $n$. 
\end{abstract}         

\maketitle  

\setcounter{tocdepth}{1} 
\tableofcontents     
        
\section{Introduction}

We use syntomic cohomology to investigate the algebraic K-theories of associative ring spectra~$K(n)$, the $(2p^n-2)$-periodic Morava K-theories with homotopy groups~$\bF_p[v_n^{\pm 1}]$. These Morava K-theories are basic objects in higher algebra -- together with $\bQ$ and~$\bF_p$, they constitute the smallest $p$-local division rings.
To the best of our knowledge, this is the first application of syntomic cohomology to the algebraic K-theory of noncommutative rings (some Morava K-theories do not admit even homotopy commutative ring structure).

In addition to proving qualitative theorems such as redshift (which seem difficult to access via simpler arguments), we give detailed computational results.  The computations suggest in a precise way that the algebraic K-theories of Morava K-theories~$\K(K(n))$ are relatively simple as spectra, or at least simpler than the algebraic K-theories of Lubin--Tate theories~$E_n$ or truncated Brown--Peterson spectra~$\BP\langle n \rangle$.

\subsection{Three qualitative questions}

In the early 2000s, Ausoni--Rognes proposed higher chromatic height analogs of the Lichtenbaum--Quillen conjectures in number theory, known as the redshift conjectures~\cite{Rog00,AR08}. These conjectures form the core of a program to study the arithmetic of ring spectra through chromatically localized algebraic K-theory, and were motivated by explicit computation of the algebraic K-theories of the Adams summand~$\ell$ and the first connective Morava K-theory~$k(1)=\ell/p$~\cite{AR02,AR12}. At least from our point of view, three of the most essential questions arising from the Ausoni--Rognes program are as follows:
 
\begin{quest}[Chromatic redshift for $\K(R)$]
\label{quest:redshift}
Given an $\bE_1$-ring~$R$ of chromatic height~$n$, does its algebraic K-theory~$\K(R)$ have chromatic height~$n+1$? 
\end{quest}

Here, we say a spectrum~$X$ has chromatic height~$n$ if its $n$-th telescopic localization~$L_{T(n)}X$ is non-zero, but its $m$-th telescopic localization $L_{T(m)}X$ vanishes for all $m>n$. 

For $\bE_\infty$-rings, Question~\ref{quest:redshift} has been spectacularly and completely answered in the affirmative by~\cite{BSY22}, building on the recent breakthroughs~\cite{Yua21,CMNN20b,LMMT24}. On the other hand, Question~\ref{quest:redshift} is known to have a negative answer for some $\bE_1$-rings~$R$, including degenerate cases where $\K(R)=0$. Nonetheless, Question~\ref{quest:redshift} has been answered in the affirmative for several of the most important $\bE_1$-rings, including $\ell/p$ and $\bE_3$-$\MU$-algebra forms of the truncated Brown--Peterson spectrum~$\BP\langle n \rangle$~\cite{AR02,AR12,AKS20,HW22,AKACHR22}. The fact that a $T(m-1)\oplus T(m)$-equivalence $A\to B$ of $\bE_1$-rings induces a $T(m)$-local equivalence $\K(A)\to \K(B)$ was first considered for the family of $\bE_1$-ring maps $y(n)\to \bF_p$ in~\cite{AKQ19} and proven in full generality in~\cite{LMMT24}. One upshot is that algebraic K-theory is known to increase chromatic height by at most $1$ for all $\bE_1$-rings. 
 
\begin{quest}[Lichtenbaum--Quillen for $\TC(R)$] 
\label{quest:fp-type}
Given an $\bE_1$-ring~$R$ of chromatic height~$n$, and a type~$n+2$ $p$-local finite complex~$V$, is $V_*\TC(R)$ bounded? 
\end{quest}

Here, we say a graded abelian group is bounded if it is only non-trivial in finitely many degrees. Notably, an affirmative answer to
Question~\ref{quest:fp-type} 
often implies the Lichtenbaum--Quillen property for $\K(R)$, which is the statement that the natural map
\[
    \K(R)_{(p)} \longrightarrow L_{n+1}^{f} \K(R)_{(p)}
\]
has bounded above fiber~\cite[Lemma 7.22]{BHLS23}.
Here, $L_{n+1}^{f}$ refers to the localization of $p$-local spectra that kills all finite complexes of type larger than $n+1$. Question~\ref{quest:fp-type} is known to hold for many fundamental ring spectra, including $\bE_3$-$\MU$-algebra forms of $\BP \langle n \rangle$~\cite{AR02,HW22,AKACHR22}, the mod $p$ Adams summand~$\ell/p$~\cite{AR12}, and some bounded below models for height $1$ local fields \cite[\S~7]{BHLS23}.

When Question~\ref{quest:fp-type} has an affirmative answer, one can furthermore ask if $V_*\TC(R)$ is \emph{finite}, in which case $\TC(R)_p^{\wedge}$ is fp in the sense of~\cite{MR99}. 
It is known that $\TC(R)_{p}^{\wedge}$ is fp when $R$ is an $\bE_3$-$\MU$-algebra form of $\BP \langle n \rangle$, and when $R$ is $\ell/p$, but $\TC(R)_p^{\wedge}$ notably need not be fp when $R$ is a bounded below model of a height $1$ local field.  

\begin{quest}[Telescope conjecture for $\K(R)$]
\label{quest:telescope}
Given an $\bE_1$-ring~$R$ of chromatic height~$n$, is the map
\[
L_{n+1}^{f} \K(R)_{(p)} \to L_{n+1} \K(R)_{(p)}
\]
an equivalence? Here, $L_{n+1}$ refers to the localization of $p$-local spectra at the Lubin--Tate theory $E_{n+1}$. 
\end{quest}

While Question~\ref{quest:telescope} has an affirmative answer for $\bE_3$-$\MU$-algebra forms of $\BP \langle n \rangle$ and for $\ell/p$, it has a negative answer when $R$ is a bounded below model of a height $1$ local field~\cite[\S~7]{BHLS23}. This is in line with a striking prediction of Mahowald--Rezk~\cite[Conjecture~7.3]{MR99} that fp spectra always satisfy the telescope conjecture. 

\subsection{Detailed results}
In this paper we make detailed computations of the algebraic $K$-theories of Morava $K$-theories. In particular, we settle the above three questions in the case that $R$ is the connective cover of $K(n)$. 
As explained in Remark \ref{remark:intro-3-lines}, our ultimate goal is to use such computations to identify the spectrum $\mathrm{K}(K(n))$. 

To begin, our formal definition of $K(n)$ is as follows:

\begin{defin} \label{intro-defin-form}
We fix for the remainder of this paper a chromatic height $n$ $\bE_1$-$\MU$-algebra $K(n)$ with homotopy algebra $\pi_*K(n) \cong \bF_p[v_n^{\pm 1}]$, where $|v_n|=2p^n-2$.  We let $k(n)$ denote its connective cover.
\end{defin}

For more details on the possible choices of Morava $K$-theory $K(n)$, see Appendix~\ref{forms}.  In particular, using work of Angeltveit we observe there that every $\bE_1$-$\bS$-algebra form of $K(n)$ admits an $\bE_1$-$\MU$-algebra structure.

In addition to studying $(2p^n-2)$-periodic $K(n)$ as we do here, another natural project would be to study the algebraic $K$-theories of residue fields of $2$-periodic Lubin--Tate theories. For some preliminary results in that direction, see Section~\ref{sec:perfect-fields}.  

Our first three main results answer Questions~\ref{quest:redshift}-\ref{quest:telescope} for $K(n)$:
\begin{thmx}[Redshift]\label{thm:redshift}
The spectrum $\K(K(n))$ has height exactly $n+1$.
\end{thmx} 

\begin{thmx}[Lichtenbaum--Quillen]\label{thm:LQ}
The spectrum $\TC(k(n))_p^{\wedge}$ is fp and the localization map 
\[ 
\K(K(n))_{(p)}\to L_{n+1}^f\K(K(n))_{(p)}
\]
has bounded above fiber. 
\end{thmx}

\begin{thmx}[Telescope]\label{thm:telescope-conj}
The localization map 
\[ 
L_{n+1}^f\K(K(n))_{(p)}\to L_{n+1}\K(K(n))_{(p)}
\]
is an equivalence. 
\end{thmx}


\begin{remark}
Suppose that $A$ is an $\bE_2$-ring and $B$ is an $\bE_1$-$A$-algebra, so that $\TC(B)$ is a unital $\TC(A)$-module. Then the telescope conjecture for $\TC(A)$ implies the telescope conjecture for $\TC(B)$, since $L_{n+1}^{f}$ is a smashing localization. Also, knowing Lichtenbaum--Quillen for $\TC(A)$ can greatly simplify proofs of Lichtenbaum--Quillen for $\TC(B)$ \cite[Proposition 3.3.7]{HW22}. For this reason, Theorems \ref{thm:LQ} and~\ref{thm:telescope-conj} are most difficult when $k(n)$ is not an algebra over an $\bE_2$-algebra form of $\BP \langle n \rangle$, and we do not assume such structure.

On the other hand, knowing redshift for $\TC(B)$ is enough to deduce it for $\TC(A)$, because any module over a zero ring must be zero. Of the above three theorems, we work hardest to prove redshift. Currently (and despite substantial effort), there is no proof of redshift for $\TC(k(n))$ that does not rely on detailed information about syntomic cohomology.
\end{remark}

Our fourth main result is an extension of a result of~\cite{AR12} to arbitrary $\bE_1$-$\MU$-algebra forms of $k(1)$ and to the prime $p=3$. This is in line with a prediction of Rognes~\cite[Chromatic redshift problem]{Rog00} from an influential talk at Oberwolfach in 2000. 

\begin{thmx}[Pure~fp-type]\label{thm:low-height}
The mod $(p,v_1)$-topological cyclic homology of $k(1)$ is an explicit finitely generated free $\bF_p[v_2]$-module at all primes~$p\ge 5$. At the prime $p=3$, it is an explicit finitely generated free $\bF_p[v_2^9]$-module.
\end{thmx}

Using \cite{DGM13}, we additionally explicitly compute the mod $(p,v_1)$ algebraic K-theory of $k(1)$ at all primes $p\ge 3$, in Theorem~\ref{thm:kk1}.

Our proofs of each of Theorem~\ref{thm:redshift}, Theorem~\ref{thm:LQ}, Theorem~\ref{thm:telescope-conj}, and Theorem~\ref{thm:low-height} rely on recent advances in prismatic cohomology and syntomic cohomology~\cite{BMS19,BL22,HRW22}. By syntomic cohomology, we will mean $\MU$-based syntomic cohomology in the sense of~\cite[Variant~A.1.9]{HRW22}. In Theorem~\ref{thm:syntomic-cohomology-kn}, we give a complete and explicit description of the mod $(p,v_1,\cdots,v_{n+1})$-syntomic cohomology of $k(n)$.  We view Theorem~\ref{thm:syntomic-cohomology-kn} as the main technical result of the paper, and  expect it to be a key step in determining the spectrum $\mathrm{K}(K(n))$ (see Remark~\ref{remark:intro-3-lines}). It has in particular the following consequence:

\begin{thmx}[Finite syntomic cohomology]\label{thm:finiteness}
The mod $(p,v_1,\cdots ,v_{n+1})$-syntomic cohomology of $k(n)$ is finite for all primes~$p$ and all heights~$n\ge 1$. 
\end{thmx}

We say that a bigraded abelian group is finite if it is finite in each bidegree and only non-trivial in finitely many bidegrees.

\begin{figure}[ht!]
\resizebox{\textwidth}{!}{ 
\begin{tikzpicture}[radius=1,yscale=2]
\foreach \n in {-2,-1,...,26} \node [below] at (\n,-.8-3) {$\n$};
\foreach \s in {-3,-2,...,3} \node [left] at (-.3-2,\s) {$\s$};
\draw [thin,color=lightgray] (-2,-3) grid (26,3);
\node [below] at (-1,1) {$\partial$};
\node [below] at (0,0) {$1$};
\node [below] at (2,0) {$\partial \bar{\varepsilon}_1$};
\node [below] at (6,0) {$\partial \bar{\varepsilon}_2$};
\node [below] at (9,-1) {$\partial \bar{\epsilon}_1\bar{\varepsilon}_2$};
\node [below] at (3,-1) {$\bar{\varepsilon}_1$};
\node [below] at (7,-1) {$\bar{\varepsilon}_2$};
\node [below] at (10,-2) {$\bar{\varepsilon}_1\bar{\varepsilon}_2$};
\node [above] at (15,1) {$\lambda_3$};
\node [above] at (13,1) {$t\lambda_3$};
\node [above] at (11,1) {$t^2\lambda_3$};
\node [above] at (9,1) {$t^3\lambda_3$};
\node [above] at (7,1) {$t^4\lambda_3$};
\node [above] at (14,2) {$\partial\lambda_3$};
\node [below] at (14,0) {$t^2\bar{\varepsilon}_1\lambda_3$};
\node [below] at (12,0) {$t^3\bar{\varepsilon}_1\lambda_3$};
\node [below] at (10,0) {$t^4\bar{\varepsilon}_1\lambda_3$};
\node [below] at (8,0) {$t^5\bar{\varepsilon}_1\lambda_3$};
\node [above] at (16,0) {$t^3\bar{\varepsilon}_2\lambda_3$};
\node [above] at (14,0) {$t^4\bar{\varepsilon}_2\lambda_3$};
\node [above] at (12,0) {$t^5\bar{\varepsilon}_2\lambda_3$};
\node [above] at (10,0) {$t^6\bar{\varepsilon}_2\lambda_3$};
\node [below] at (17,-1) {$t^4\bar{\varepsilon}_1\bar{\varepsilon}_2\lambda_3$};
\node [below] at (15,-1) {$t^5\bar{\varepsilon}_1\bar{\varepsilon}_2\lambda_3$};
\node [below] at (13,-1) {$t^6\bar{\varepsilon}_1\bar{\varepsilon}_2\lambda_3$};
\node [below] at (11,-1) {$t^7\bar{\varepsilon}_1\bar{\varepsilon}_2\lambda_3$};
\node [above] at (18,0) {$\bar{\varepsilon}_1\lambda_3$};
\node [above] at (17,1) {$\partial\bar{\varepsilon}_1\lambda_3$};
\node [above] at (22,0) {$\bar{\varepsilon}_2\lambda_3$};
\node [above] at (21,1) {$\partial\bar{\varepsilon}_2\lambda_3$};
\node [above] at (25,-1) {$\bar{\varepsilon}_1\bar{\varepsilon}_2\lambda_3$};
\node [above] at (24,0) {$\partial\bar{\varepsilon}_1\bar{\varepsilon}_2\lambda_3$};
\end{tikzpicture}
}
\caption{The mod $(2,v_1,v_2,v_3)$-syntomic cohomology of $k(2)$.}
\label{fig:k2}
\end{figure}

\begin{rem2}
Figure~\ref{fig:k2} exhibits rotational symmetry akin to Poitou--Tate duality, as conjectured by Rognes~\cite[Example~5.2]{Rog14} and exhibited previously in the cases of the Adams summand~\cite{HRW22} and real topological K-theory~\cite{AKAR23}.  The duality is closely related to a Lagrangian refinement of Tate duality in the setting of prismatic F-gauges~\cite[Theorem~1.3.1]{Bha22}. 
\end{rem2}
\begin{rem2}
Work of Pstr\k{a}gowski~\cite{Pst23} beautifully extends the motivic filtration on $\THH(R)$, defined for 
$\bE_\infty$-rings $R$ in~\cite{HRW22}, to a filtration defined on $\bE_2$-rings. Pstr\k{a}gowski's work does not apply to topological Hochschild homology of $\bE_1$-algebras like $k(n)$. Nonetheless, one of the main morals of this paper is that syntomic cohomology can be a very practical tool for understanding the algebraic K-theory of such $\bE_1$-ring spectra. 

In particular, if $A$ is an $\mathbb{E}_{\infty}$-ring spectrum and $B$ is an $\bE_1$-$A$-algebra, then $\THH(B)$ is an $\bE_0$-$\THH(A)$-algebra.  This leads to a spectral sequence computing $\pi_*\TC(B)$ beginning with coherent cohomology over the syntomification $A^{\Syn}$. 

In this paper, we execute (a Nygaard completed variant of) this strategy in the case $A=\MU$ and $B=k(n)$.  That the resulting spectral sequence is powerful in practice is closely related to the fact that the unit map $\pi_*\MU \to \pi_*k(n)$ is surjective. 
This should be contrasted with the case of discrete rings, where a discrete noncommutative ring may never receive a surjective map from a discrete commutative ring.  It is a phenomenon unique to higher algebra that a noncommutative ring may receive a $\pi_*$-surjection from a commutative ring.
\end{rem2}

\begin{rem2} \label{remark:intro-3-lines}
We will prove in Proposition~\ref{prop:integral-collapse} that the motivic spectral sequence computing the integral homotopy groups~$\pi_*\TC(k(n))_p^{\wedge}$ is concentrated on at most three horizontal lines.  It follows that there are no differentials in this motivic spectral sequence, so the homotopy groups of the spectrum $\TC(k(n))_p^{\wedge}$ are essentially determined by the $\EE_2$-page.  Analogously, understanding the spectrum $\TC(k(n))_p^{\wedge}$ should be very closely related to understanding the sheaf over the moduli stack of formal groups $\cM_{\fg}$ whose cohomology produces the $\EE_2$-page. The main purpose of our calculations here, beyond proving the qualitative Theorems~\ref{thm:redshift}, \ref{thm:LQ} and~\ref{thm:telescope-conj} and the quantitative Theorems~\ref{thm:low-height}, \ref{thm:finiteness}, and~\ref{thm:syntomic-cohomology-kn}, is to assist in the eventual identification of this sheaf.

  All this may be contrasted with the motivic spectral sequence computing $\pi_*\TC(\BP \langle n \rangle)_p^{\wedge}$ or the motivic spectral sequence computing the telescopic homotopy of $\TC(\BP \langle n \rangle)_p^{\wedge}$ (which is closely related to the telescopic homotopy of $\K(E_n)$).  The number of lines in those spectral sequences grows as a linear function of $n$.
\end{rem2}

\begin{remark}
The situation is even better when studying Morava K-theories base changed over an algebraically closed field. As we point out at the end of this paper, the motivic spectral sequence computing 
\[
\pi_*\TC(k(n) \otimes \bS_{\W(\overline{\bF}_p)})_p^{\wedge}
\] 
is concentrated on just two horizontal lines.\footnote{Here $\bS_{\W(\overline{\bF}_p)}$ denotes the spherical Witt vectors from~\cite[Example~5.2.7]{Lur21}.} This strongly suggests that, for any height $n>0$, the spectrum~$\TC(k(n) \otimes \bS_{\W(\overline{\bF}_p)})_p^{\wedge}$ is the fiber of a map between two even $\MU$ modules.
\end{remark}

\subsection{Acknowledgements}
The authors benefited from conversations with Christian Ausoni, Robert Burklund, Ishan Levy, Piotr Pstr\k{a}gowski, Arpon Raksit, John Rognes, Andrew Senger, and Tristan Yang. 
Hahn was
supported by the Sloan Foundation and by a grant from the Institute for Advanced Study
School of Mathematics. Angelini-Knoll is grateful to Max Planck Institute for Mathematics in
Bonn for its hospitality and financial support.
This project received funding from the European Union's Horizon 2020 research and innovation programme under the Marie Sk\l{}odowska-Curie grant agreement No 1010342555.
\thinspace \includegraphics[scale=0.1]{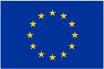} 
\section{Hochschild homology and the motivic filtration}\label{sec:hochschild}

Recall from Definition~\ref{intro-defin-form} that we fix an (arbitrary) $\bE_1$-$\MU$-algebra form of $K(n)$, which induces an $\bE_1$-$\MU$-algebra structure on the connective cover $k(n)$. More details about the possible choices of $K(n)$ are contained in Appendix~\ref{forms}. Following that appendix, we name generators 
\[
p,v_1,\cdots \in \pi_*\BP
\]
such that, for all $i \ne n$, $v_i$ maps to zero under the composite map
\[
\pi_*\BP \to \pi_*\MU_{(p)} \to \pi_* k(n).
\]

In Section~\ref{Hochschildstructures}, we review various structures on Hochschild homology. We compute the motivically filtered Hochschild homologies of $k(n)$ and $\BP$, with $\bF_p$ coefficients, in Section~\ref{sec:HochschildBPnknFp-coeff}. We discuss the mod $(p,v_1,\cdots,v_n)$ Hochschild homology of $\bF_p$ in Section~\ref{HochschildFp} in order to fix notation for certain classes. Finally, we compute the mod $(p,v_1,\cdots,v_n)$ motivic filtration on the Hochschild homology of $k(n)$, identifying it explicitly in terms of its image in the mod $(p,v_1,\cdots,v_n)$ Hochschild homology of $\bF_p$ in Section~\ref{sec:modHH}. Specifically, we will identify the map
\[
\pi_*\gr^*_\mot \THH(k(n)) / (p,v_1,\cdots,v_n) \to \pi_*\gr^*_{\mot} \THH(\bF_p) / (p,v_1,\cdots,v_n)
\]
induced by the map of $\bE_1$-$\MU$-algebras $k(n)\to \tau_{\le 0}k(n)=\bF_p$ (hereafter, we will often call maps associated to $k(n) \to \bF_p$ reduction maps).

\subsection{Structures on Hochschild homology}\label{Hochschildstructures}
In this section, we review various structures on Hochschild homology such as the motivic filtration in Section~\ref{sec:mot}, the $\sigma$ operator in Section~\ref{sec:sigma}, the cap product in Section~\ref{sec:cap}, and the double suspension map in Section~\ref{sec:suspension}.

\subsubsection{The motivic filtration}\label{sec:mot}
We review the definition of the motivic filtration based on $\MU$, defined on the $\THH$ of any $\bE_1$-$\MU$-algebra in~\cite[Appendix~A]{HRW22}. 

\begin{rec}
If $A$ is an $\bE_{\infty}$-ring and $M$ is an $A$-module, we  define
\[
\fil_{\ev/A}^*M
\] 
following~\cite[Construction~A.1.3.]{HRW22}. 
In the case that $A=\THH(\MU)$, we have an explicit description 
\[
\fil^*_{\ev/\THH(\MU)} M = \mathrm{Tot}\left( \tau_{\ge 2*} \left ( {M}\otimes_{\THH(\MU)}{\MU}^{\otimes_{\THH(\MU)}\bullet+1} \right) \right) 
\]
by~\cite[Example~4.2.3]{HRW22} and~\cite[Corollary~A.2.6]{HRW22}. 
Moreover, if $M$ is an $\bE_2$-$\MU$-algebra and $N$ is an $\bE_1$-$M$-algebra then 
\[
\fil_{\ev/\THH(\MU)}^*\THH(N)
\] 
is an $\bE_0$-$\fil_{\ev/\THH(\MU)}^*\THH(M)$-algebra.
\end{rec}

\begin{remark}
Note that the construction $\fil_{\ev/A}^*M$ is functorial in maps of pairs~$(A,M)\to (B,N)$ where $A$ and $B$ are $\bE_{\infty}$-rings, $M$ is an $A$-module, and $N$ is a $B$-module.  
\end{remark}

In this particular paper, we will always study motivic filtrations relative to $\MU$.  Thus, we make the following further simplifying convention:

\begin{convention}\label{conv:motfilt}
Given a $\THH(\MU)$-module $M$, we write
\[
\filmot^*M := \fil^*_{\ev/\THH(\MU)}M\,. 
\]
\end{convention}

\begin{convention}
If $M$ is a graded spectrum and $y\in  \pi_{s}M^t$, then we write 
\[ 
\| y \|: =(s,2t-s)
\]
and refer to $s$ as the \emph{degree} and $2t-s$ as the \emph{Adams weight}. When we refer to the bidegrees of $\pi_*M^*$, we will always mean the bidegree $(\textup{degree},\textup{Adams weight})$.
\end{convention}

\begin{defin}\label{ss-not}
In a spectral sequence 
\[ 
\EE_1^{*,*}\implies G_*
\]
we write $\{x\}$ for the coset of elements $\xi\in G_*$ that are detected by $x$. Sometimes we will write $[[x]]$ for a specific choice of such an element so that $[[x]]\in \{x\}$ (cf.~\cite[Notation~7.2]{AKACHR22}). 
\end{defin}
\begin{remark}
Note that the homotopy groups of $\gr_{\ev/\bS}^*\MU\simeq \gr_{\ev}^*\MU$ can be identified with the $\EE_2$-page of the Adams--Novikov spectral sequence converging to $\pi_*\MU$. Consequently, every element of $\pi_*\grev^*\MU\cong \MU_*$ has Adams weight $0$, and we say that the bigraded homotopy groups are concentrated on the $0$-line. 
\end{remark}

\begin{defin}\label{def:vi}
For each $i>0$, let $v_i\in \pi_*\grev^*\MU$ be the class in bidegree 
\[ 
\| v_i\| =(2p^i-2,0)
\]
fixed in Convention~\ref{conv:vi}. Let $\grev^*\MU/v_{i}$ denote the cofiber of $v_{i}$ considered as a self-map of $\grev^*\MU$. For $M$ a $\grev^*\MU$-module, we then define
\[ 
M/v_{i}:=M\otimes_{\grev^*\MU} (\grev^*\MU/v_{i}) \,.
\]
We sometimes write $M/(p,v_1,\cdots ,v_{i})$ for the iterated tensor product
\[
M\otimes_{\grev^*\MU} (\grev^*\MU/p) \otimes_{\grev^*\MU} \cdots \otimes_{\grev^*\MU} (\grev^*\MU/v_i) \,.
\]
\end{defin}

\begin{rmk}\label{rem:vi-thh}
Via the unit map $\MU \to \THH(\MU)$, any $\grmot^* \THH(\MU)$-module is also a $\grev^*\MU$-module.  Thus, if $R$ is an $\bE_1$-$\MU$-algebra, we may refer to
\[
\gr^*_{\mot}\THH(R) / (p,v_1,\cdots,v_i) \,.
\]
\end{rmk}

\begin{rmk}\label{rem:vi-elements}
Our convention that $v_i \in \pi_{2p^i-2} \MU$ is in the image of the unit map 
\[
\pi_{2p^i-2}\BP \to \pi_{2p^i-2}\MU
\]
for each $i \ge 1$ ensures that the natural map
\[
\grev^* \bS / (p,v_1,\cdots,v_{i-1})  \to \grev^* \MU / (p,v_1,\cdots,v_{i-1})
\]
sends $v_i$ to $v_i$. In particular, if $M$ is a $\grev^*\MU$-module, then
\[
M \otimes_{\grev^*\bS} \grev^*\bS/(p,v_1,\cdots,v_i) \cong M \otimes_{\grev^* \MU} \grev^*\MU / (p,v_1,\cdots,v_i)\,.
\]
\end{rmk}

\subsubsection{The sigma operator}\label{sec:sigma}
Here we construct an operator
\[
\sigma: \pi_*\grmot^* \THH(R) \longrightarrow \pi_{*+1} \grmot^{*+1}\THH(R)
\]
that is natural in maps of $\bE_1$-$\MU$-algebras~$R$. This is a motivic refinement of the classical sigma operator, which records the circle action on $\THH(R)$.

To begin, we observe that $\THH(R)$ is an $\bS[S^1]$-module, which formally leads to an action of $\fil_{\ev/\bS}^*\bS[S^1]$ on $\filmot^* \THH(R)$. We then define the motivic $\sigma$ operator to be multiplication by the fundamental class of $S^1$, which is a specific class in $\pi_1 \gr_{\ev/\bS}^1 \bS[S^1]$. 

This action may be understood explicitly in terms of choices of eff covers, as we now explain: The map of $\bE_{\infty}$ rings $\bS[S^1]\to \MU$ is eff and $\MU$ is even, so the fundamental class of $S^1$ is a specific element in 
\[
\pi_1\gr^1_{\ev/\bS} \bS[S^1] \cong \pi_1 \left(\lim_{\Delta} \pi_2(\MU\otimes \MU[(BS^1)^{\otimes \bullet}]) \right) \,,
\]
which is detected by the bottom cell of $BS^1$ in cosimplicial level~$1$. It then suffices to observe that the cosimplicial $\bE_{\infty}$-ring
\[
\bS[S^1] \otimes_{\bS[S^1]} \MU^{\otimes_{\bS[S^1]} \bullet +1} =\MU\otimes \MU[(BS^1)^{\otimes \bullet}]
\]
acts naturally on the cosimplicial spectrum
\[
\THH(R) \otimes_{\THH(\MU)} \MU^{\otimes_{\THH(\MU)} \bullet+1} \,,
\]
inducing an action of the filtered $\bE_\infty$-ring
\[
\fil_{\ev/\bS}^* \bS[S^1] = \lim_{\Delta} \left(\tau_{\ge 2*} \left(\bS[S^1] \otimes_{\bS[S^1]} \MU^{\otimes_{\bS[S^1]} \bullet +1} \right) \right) 
\]
on 
\[ 
\filmot^*\THH(R)= \lim_{\Delta} \left (\tau_{\ge 2*} \left (\THH(R)\otimes_{\THH(\MU)}\MU^{\otimes_{\THH(\MU)} \bullet+1} \right )\right) \,.
\] 

\subsubsection{The cap product}\label{sec:cap}
Let $R$ denote an $\bE_1$-$\MU$-algebra equipped with an augmentation $R \to \bF_p$ of $\bE_1$-$\MU$-algebras. We regard $\bF_p$ as a $R\otimes R^{\op}$-module by restriction along the map 
\[
R\otimes R^{\op}\to  \bF_p\otimes  \bF_p^{\op}\to \bF_p 
\] 
of $\bE_1$-$\MU\otimes \MU^{\op}$-algebras. 
There is then a cap product pairing  
\[ 
\THH(R;\bF_p)\otimes_{\THH(\MU;\bF_p)} \THC(\bF_p)\longrightarrow \THH(R;\bF_p)
\]
using the identifications $\THH(R;\bF_p):=\bF_p\otimes_{R\otimes R^{\op}}R$ and $\THC(\mathbb{F}_p):=\map_{\bF_p\otimes \bF_p^{\op}}(\bF_p,\bF_p)$ as $\THH(\MU;\bF_p)$-modules (cf.~\cite[\S~5]{AHL10}). 

Recall that $\THH_*(\bF_p) = \THH_*(\bF_p;\bF_p) \cong \bF_p[\mu]$ by~\cite{Bre78,Bok87a}, where $|\mu|=2$. Furthermore, 
\[ 
\map_{\bF_p\otimes \bF_p^{\op}}(\bF_p,\bF_p)\simeq  \map_{\bF_p} (\THH(\bF_p),\bF_p),
\]
so $\pi_*\THC(\bF_p)$ is the $\bF_p$-linear dual of $\THH(\bF_p)$. In particular, writing $c_j \in \pi_{-2j} \THC(\bF_p)$ for the dual of $\mu^j \in \pi_{2j} \THH(\bF_p)$, then the upshot of this discussion is a natural map 
\[
\-- \cap c_j: \THH(R;\bF_p) \to \Sigma^{2j} \THH(R;\bF_p)
\]
of $\THH(\MU)$-modules. Applying $\grmot^*$, one obtains a motivic cap product 

\[
\-- \cap c_j: \pi_* \grmot^* \THH(R;\bF_p) \to \pi_* \grmot^* \THH(R;\bF_p) \,,
\]
which changes bidegree by $(-2j,0)$.

\subsubsection{The double suspension map}\label{sec:suspension}
Following~\cite[Appendix~A]{HW22}, given an $\bE_\infty$-ring $A$ and an $\bE_1$-$A$-algebra $B$, there is a canonical map 
\[ 
\sigma^2 : \Sigma \mathrm{cof}(A \to B)  \longrightarrow \THH(B/A) 
\]
where $\mathrm{cof}(A \to B)$ denotes the cofiber of the unit map~$A\to B$. We will later use the following basic fact:

\begin{lem}\label{lem:sigma-squared-commutes-with-power-operations}
Let $A$ be an $\bE_\infty$-ring and let $B$ be an $\bE_\infty$-$A$-algebra. Then the (twice desuspended) $\sigma^2$ map
\[ 
\mathrm{fib}(A \to B) \to \Sigma^{-2} \THH(B/A)
\]
is a map of non-unital $\bE_{\infty}$-ring spectra.
In particular, $\sigma^2$ commutes with stable $\bE_{\infty}$-power operations.
\end{lem}
\begin{proof}
Consider the diagram of $S^1$-equivariant $\bE_{\infty}$-rings
\begin{equation}\label{diag:ofE2-rings}
\begin{tikzcd}
A\ar[r]
 \ar[d]& \THH(B/A)\ar[d] \\ 
B\ar[r] & \map(S^1_+,\THH(B/A))
\end{tikzcd}
\end{equation}
The induced map on vertical fibers is the map 
\[
\Sigma^{-2}\sigma^2 : \mathrm{fiber}(A \to B) \to \Sigma^{-2} \THH(B/A)
\]
from~\cite[Construction~A.1.2,~Example~A.2.4]{HW22} (up to suspension). The result then follows because the fiber of a map of $\bE_{\infty}$-rings is naturally a non-unital $\bE_{\infty}$-ring. 
\end{proof}

\begin{warning}
The suspension map $\sigma^2$ is distinct from the composite $\sigma \circ \sigma$ of the $\sigma$-operator from Section~\ref{sec:sigma}. 
\end{warning}

\subsection{Hochschild homology with $\bF_p$ coefficients}\label{sec:HochschildBPnknFp-coeff}


We compute Hochschild homology of $k(n)$ with $\bF_{p}$-coefficients in Section~\ref{HH} along with its motivic filtration in Section~\ref{motHH}. 

\subsubsection{Hochschild homology with $\bF_p$-coefficients}\label{HH}
We begin with results that do not involve the motivic filtration. 

\begin{rec}\label{rec:THHMU}
There is an isomorphism 
\[ 
\THH_*(\MU;\bF_p)=\Lambda(\lambda_i^{\prime}: i\ge 1)
\]
of graded $\bF_p$-algebras, where $|\lambda_i^{\prime}|=2i+1$ by \cite[Remark~4.3]{MS93}. The $\bE_3$-algebra retract $\THH(\BP;\bF_p)$ of $\THH(\MU;\bF_p)$ has homotopy groups
\[
\THH_*(\BP;\bF_p) \cong  \Lambda(\lambda_i: i \ge 1) 
\]
where we write 
\[ 
\lambda_j\coloneqq\lambda_{p^j-1}^{\prime} \,.
\]
\end{rec}
\begin{proposition}\label{prop:Hocschild-MaySSBP}
The $\bF_p$-Adams filtration on $\BP$ gives rise to a spectral sequence 
\[
\THH_*(\bF_p[v_0,v_1,\cdots];\bF_p) \implies \THH_*(\BP;\bF_p) \,.
\]
This has $\EE_1$-page $\bF_p[\mu] \otimes \Lambda(d v_0, d v_1,\cdots)$. It has differentials
\[
d_{2p^k-2}(\mu^{p^k}) = d v_k \,,
\]
leaving an $\EE_{\infty}$-page of
\[
\Lambda(\mu^{p-1} d v_0, \mu^{p^2-p} d v_1,\cdots) \,,
\]
where $\lambda_{k+1}$ is detected by $\mu^{p^k-p^{k-1}}dv_k$ without indeterminacy for $k\ge 0$.
\end{proposition}

\begin{proof}
We apply the Hochschild--May spectral sequence~\cite{AKS18,Kee20,JLL23} to the $\bE_{4}$ algebra $R_{*}=\lim_{\Delta} \tau_{\ge 2*} (\BP\otimes \bF_p^{\otimes \bullet+1})$ in filtered spectra. Here we know $R_*$ is a $\bE_{4}$ algebra by \cite[Theorem~1.1]{BM13} and~\cite[Theorem~5.64]{PP21}. By~\cite[Proposition~4.2.1]{HW18}, we can identify the associated graded of $R_{*}$ with $\bF_p\otimes \bS[v_0,v_1,\cdots ]$ as an $\bE_2$ algebra and the identification of the $\EE_1$-page follows as in~\cite[Lemma 4.1.3]{HW18}. Since we know that the abutment is $\THH_*(\BP;\bF_p)=\Lambda(\lambda_i : i\ge 0)$ by Recollection~\ref{rec:THHMU}, the classes in stems $2p^k$ for $k\ge 0$ must not survive by elementary degree arguments. Note that the Hochschild--May spectral sequence is a first quadrant spectral sequence with differentials satisfying the Adams convention. Since the classes $\mu^{p^k}$ are in stem $2p^{k}$ and they are on the $0$-line, the differentials $d_{2p^k-2}(\mu^{p^k})=d v_k$ are forced for $k\ge 0$. Any further differentials would cause the abutment to have a smaller dimension as an $\bF_p$-vector space than the known abutment. There is no room for hidden multiplicative extensions. The fact that $\lambda_{k+1}$ is detected by $\mu^{p^{k}-p^{k-1}} d v_k$ without indeterminacy follows because there are no classes in higher Hochschild--May filtration than $\mu^{p^{k}-p^{k-1}} d v_k$  in degree $2p^{k}-1$ by inspection. 
\end{proof}

\begin{proposition}\label{prop:Hochschild-MaySSkn}
As an $\bE_0$-$\THH_*(\BP)$-algebra, there is a preferred isomorphism 
\[
\THH_*(k(n);\bF_p) \cong \bF_p [\mu^{p^{n+1}}]\otimes \bF_{p}[\mu]/(\mu^{p^{n}})\otimes \Lambda ( \lambda_{n+1} ).
\]
Classes are named such that the reduction map of $\mathbb{F}_p$ vector spaces
\[\THH_*(k(n);\bF_p)\to \THH_*(\bF_p)\]
is given by the quotient by $(\lambda_{n+1})$ followed by the injection of named classes into the polynomial ring $\mathbb{F}_p[\mu]$. 
\end{proposition}

\begin{proof}
Consider the Hochschild--May spectral sequence, with $\EE_1$-page 
\[
\THH_*(\bF_p[v_n];\bF_p) \,,
\] 
as an $\bE_0$-algebra over the Hochschild--May spectral sequence for $\BP$ from Proposition~\ref{prop:Hocschild-MaySSBP}. Explicitly, we know $k(n)$ is an $\bE_1$-$\BP$-algebra and consequently, we know that 
\[
K_{*}:=\lim_{\Delta}\tau_{\ge 2*}(k(n) \otimes \bF_p^{\otimes \bullet+1})
\] 
is an $\bE_1$-$R_{*}$-algebra by~\cite[Theorem~5.64]{PP21}, where $R_{*}=\lim_{\Delta}\tau_{\ge 2*}(\BP \otimes \bF_p^{\otimes \bullet+1})$. 
The $\EE_1$-page is
\[
\bF_p[\mu]\otimes \Lambda(d v_n)
\] 
and the $\bE_{0}$-algebra structure over the Hochschild--May spectral sequence for $\BP$ implies a differential~$d_{2p^n-2}(\mu^{p^n})=d v_n$ that generates the differentials using the Leibniz rule for the module structure over the Hochschild--May spectral sequence for $\BP$. At the $\EE_{2p^n-1}$-page, the map of Hochschild--May spectral sequences induced by the reduction map $k(n)\to \bF_p$ is injective mod $(\mu^{p^{n+1}-p^n}d v_n)$ and the target spectral sequence collapses. We know that $\mu^{p^{n+1}-p^n}d v_n$ is a permanent cycle in the Hochschild--May spectral sequence for $\BP$ by Proposition~\ref{prop:Hocschild-MaySSBP}. 
Therefore, the Hochschild--May spectral sequence for $k(n)$ collapses at this page.  
We then observe that the map 
$\THH_*(k(n);\bF_p)\longrightarrow \THH_*(\bF_p)$
induced by the reduction map $k(n)\to \bF_p$ is the canonical quotient by $(\lambda_{n+1})$ followed by the canonical inclusion $\bF_p[\mu^{p^{n+1}}]\otimes \bF_p[\mu]/(\mu^{p^n})\subset \bF_p[\mu]$ of graded $\bF_p$-vector spaces and we therefore name classes accordingly. 
There is no room for hidden $\bE_0$-$\THH(\BP)$-algebra extensions. 
For bidegree reasons, there is no indeterminacy in the choice of classes $\mu^k$ and  $\mu^{k}\mu^{p^{n+1}-p^n}d v_{n}$ detecting $\mu^k$ and $\mu^k\lambda_{n+1}$ for $k\equiv 0,1,\cdots,p^n\mod p^{n+1}$. The fact that the map
\[ \THH_*(k(n);\bF_p) \to \THH_*(\bF_p) \]
is given as described follows from the fact that there is a map of filtered objects 
\[
K_{*} \longrightarrow \lim_{\Delta }\tau_{\ge 2*}\bF_p \otimes \bF_p^{\otimes  \bullet+1},
\] 
producing a map of spectral sequences as $\bE_0$-algebras over the spectral sequence for $\BP$.  The target collapses at the $\bE_2$-term, and the map on $\bE_2$-terms is an injection modulo $dv_n$. Consequently, we observe that it is an injection modulo $(\lambda_{n+1})$ on the abutment. 
\end{proof}

\subsubsection{The motivic filtration with $\bF_p$ coefficients}\label{motHH}
The purpose of this section is to compute 
\[
\pi_*\grmot^* \THH(k(n);\bF_p) \,,
\]
defined using the $\THH(\MU)$-module structure on $\THH(k(n);\bF_p)$ as in Convention~\ref{conv:motfilt}.
These bigraded homotopy groups form the $\EE_2$-page of a motivic spectral sequence
\[
\pi_*\grmot^*\THH(k(n);\bF_p)\implies \THH_*(k(n);\bF_p)
\]
abutting to the known result from Proposition~\ref{prop:Hochschild-MaySSkn}.  We will deduce that this motivic spectral sequence degenerates at the $\EE_2$-page without extensions. 
We first introduce some notation using the following remark. 
\begin{rem2}\label{rem:Fpevenflat}
By Recollection~\ref{rec:THHMU}, we can identify 
\[
\pi_*(\bF_p\otimes_{\THH_*(\MU;\bF_p)}^{\mathbb{L}}\bF_p)  =\Gamma (d \lambda'_i : i\ge 1) 
\]
with $|d \lambda'_i|=2i$. 
\end{rem2}
\begin{defin} 
We write $(\bF_p,\Gamma)$ for the Hopf algebra 
\[ 
(\pi_*\bF_p ,  \pi_*(\bF_p\otimes_{\THH_*(\MU;\bF_p)}^{\mathbb{L}}\bF_p) ) =(\bF_p, \Gamma (d\lambda'_i : i\ge 1) )
\]
where $|d \lambda'_i|=2i+2$ and we define 
\[ H^*(\Gamma ,N):= \Cotor_{(\bF_p,\Gamma )}^{*}(\bF_p,N)
\]
where $N$ is a $\Gamma$-comodule. 
\end{defin}

We now consider the case of $\MU$. 
\begin{proposition}\label{prop:thh-mu} 
We can identify 
\[
	\pi_*\grmot^*\THH(\MU;\bF_p )
	\cong
\Lambda(\lambda_i^{\prime} : i\ge 1)
\] 
as a bigraded $\bF_p$-algebra, where the bidegrees of indecomposable algebra generators are as follows:
\begin{align*}
	\|\lambda_i^{\prime}\| &= (2i+1, 1) \,.
\end{align*}
Moreover, the motivic spectral sequence collapses at the $\EE_2$-page without extensions.
\end{proposition}
\begin{proof}
We filter the Hopf algebroid 
$(\bF_p ,\pi_{*}(\bF_p \otimes_{\THH(\MU;\bF_p )}\bF_p ))$ 
by  
\[ 
(\bF_p ,\pi_{*}(\tau_{\ge \bullet }\bF_p \otimes_{\tau_{\ge \bullet}\THH(\MU;\bF_p )}\tau_{\ge \bullet}\bF_p )),
\] 
with associated graded 
$(\bF_p,\Gamma)$ 
by collapse of the K\"unneth spectral sequence. 
The cohomology of $(\bF_p,\Gamma)$ is then the input of a May--Ravenel spectral sequence\footnote{Here we refer to the spectral sequence from~\cite[Theorem A1.3.9]{Rav86} as the May--Ravenel spectral sequence. Note that Salch~\cite{Sal18} refers to a specific instance of this spectral sequence as the Ravenel--May spectral sequence and this is a different instance of this spectral sequence than the one we use.}~\cite[Theorem A1.3.9]{Rav86} 
\[
H^*(\Gamma,\bF_p) \implies \pi_{*}\grmot^{*}\THH(\MU;\bF_p ) \,.
\]
Since we can compute that 
\begin{equation}\label{eq:cotor=abutment-MU}
H^*(\Gamma,\bF_p)=\Lambda(\lambda_{i}^{\prime} : i\ge 1), 
\end{equation}
both the May--Ravenel spectral sequence and the motivic spectral sequence 
\[ 
\pi_{*}\grmot^{*}\THH(\MU;\bF_p )\implies \THH_{*}(\MU;\bF_p )
\]
degenerate at the $\EE_2$-page without extensions. For $i\ge 1$, we observe that the Adams weight of $\lambda_{i}^{\prime}$ is $1$ and there is no indeterminacy in the choice of classes $\lambda'_i\in \{\lambda'_i\}$ for bidegree reasons. 
\end{proof}

\begin{remark}\label{rem: grading}
A priori, the left-hand side of Equation~\eqref{eq:cotor=abutment-MU} is trigraded. However, we can also regard it as bigraded using the total grading coming from the the May--Ravenel spectral sequence. Moreover, we can further regard it as singly graded using the total grading coming from the motivic spectral sequence. We abuse notation and we do not distinguish between these three different grading conventions and the meaning will be clear from the context. 
\end{remark}

We will the give proofs of the following three results at the end of this subsection. 
\begin{proposition}\label{prop:thh-fp} 
We can identify 
\[
	\pi_*\grmot^*\THH(\bF_p )
	\cong
	\bF_p [\mu]
\] 
as a bigraded $\bF_p$-algebra where the bidegree of $\mu$ is 
\[ 
\|\mu\| = (2, 0) \,.
\]
Moreover, the motivic spectral sequence degenerates at the $\EE_2$-page without extensions.
\end{proposition}

\begin{proposition}\label{prop:thh-bp} 
We can identify 
\[
	\pi_*\grmot^*\THH(\BP;\bF_p )
	\cong
	\Lambda(\lambda_i : i\ge 1)
 \] 
	as a bigraded $\bF_p$-algebra where the bidegrees of indecomposable algebra generators follow: 
\begin{align*}
	\|\lambda_i\| &= (2p^i+1, 1) \,.
\end{align*}
Moreover, the motivic spectral sequence degenerates at the $\EE_2$-page without extensions.
\end{proposition}

\begin{proposition}\label{prop:thh-kn}
There is an isomorphism of $\bE_{0}$-$\pi_*\grmot^*\THH(\MU;\bF_p )$-algebras
\[
	\pi_*\grmot^*\THH(k(n);\bF_p )
	\cong 
	\bF_p[\mu]/(\mu^{p^n})\otimes \bF_p[\mu^{p^{n+1}}] \otimes \Lambda(\lambda_{n+1})
\] 
where the bidegrees are as follows:
\begin{align*}
	\|\mu^j\| &= (2j,0)\,, \\
 \|\lambda_{n+1}\| &=(2p^{n+1}-1, 1)\,.
\end{align*}
The $\bE_{0}$-$\pi_*\grmot^*(\THH(\MU;\bF_p )$-algebra structure map 
\[
\pi_*\grmot^*(\THH(\MU;\bF_p ))\to \pi_*\grmot^*(\THH(k(n);\bF_p ))
\] 
is given by 
\[
		\lambda_{i}^{\prime}\mapsto \begin{cases} \lambda_{n+1} & \text{ if } i=p^{n+1}-1 \\ 0 & \text{ otherwise.} \end{cases} 
\]
Moreover, the motivic spectral sequence degenerates at the $\EE_2$-page without extensions.
\end{proposition}

\begin{corollary}\label{cor:vanishing-line}
The motivic spectral sequence 
\[ 
\pi_*\grmot^*\THH(k(n))\implies \THH_*(k(n))
\]
is supported on the $0$-line and the $1$-line. 
\end{corollary}
\begin{proof}
Note that $\|v_n\|=(2p^n-2,0)$ so the abutment of the $v_n$-Bockstein spectral sequence 
\[ 
\pi_*\grmot^*\THH(k(n);\mathbb{F}_p)[v_n]\implies \pi_*\grmot^*\THH(k(n))
\]
is concentrated in Adams weights $[0,1]$ by Proposition~\ref{prop:thh-kn}. 
\end{proof}

\begin{defin}
Given a $\THH_*(\MU;\bF_p )$-module $M$, we will say that $M$ is even flat if 
\[  
M_{*}\otimes_{\THH_{*}(\MU;\bF_p )}^{\mathbb{L}} \bF_p 
\]
is concentrated in even total degrees. 
\end{defin}
\begin{remark}
By \cite[Theorem~4.19]{Pst23}, this corresponds to the definition of even flat in \cite[Definition~4.2]{Pst23} where we regard $\THH_*(\MU;\bF_p)$ as an $\bE_{\infty}$-$\bF_p$-algebra. 
\end{remark}
\begin{exm}
By Remark~\ref{rem:Fpevenflat}, $\bF_p$ is even flat as a $\THH_*(\MU;\bF_p)$-module. 
\end{exm}
\begin{defin}
Given an augmentation $R\to \bF_p$ of $\bE_1$-$\MU$-algebras, we will write 
\[ \THH(R;\bF_p)=\bF_p\otimes_{R\otimes R^{\op}}R \,, \]
where $\bF_p$ is regarded as a $R\otimes R^{\op}$-module by restriction along the map 
\[ R\otimes R^{\op}\to \bF_p\otimes \bF_p^{\op}\longrightarrow \bF_p
\]
of $\bE_1$-$\MU\otimes \MU^{\op}$-algebras. This exhibits $\THH(R;\bF_p)$ as a $\bE_0$-$\THH(\MU;\bF_p)$-algebra. 
\end{defin}
\begin{prop}\label{prop: identification of associated graded of motivic filtration}
Suppose $R\to \bF_p$ is an augmentation of $\bE_1$-$\MU$-algebras and $\THH_*(R;\bF_p )$ is an even flat $\THH_*(\MU;\bF_p)$-module. Then, there is a conditionally convergent spectral sequence
\[ 
	H^*(\Gamma,\THH_*(R;\bF_p )\otimes_{\THH_*(\MU;\bF_p)}^{\mathbb{L}}\bF_p  )\implies \pi_{*}\grmot^{*}\THH_*(R;\bF_p) \,.
\]
\end{prop}
\begin{proof}
Since $\THH_*(R;\bF_p )$ is even flat, we know that $\THH_*(R/\MU^{\otimes q+1} ;\bF_p)$ is concentrated in even degrees for each $q\ge 0$. 
Since the map $\THH(\mathrm{MU})\to \mathrm{MU}$ is eff by \cite[Example~4.2.3]{HRW22}, we know $\THH(\mathrm{MU};\mathbb{F}_p)\to \mathbb{F}_p$ is eff by base-change. We can therefore identify $\pi_*\grmot^*\THH(R;\bF_p)$ with the cohomology of the cobar complex
\[\lim_{\Delta}\pi_{2*}\left ( \THH(R;\bF_p)\otimes_{\THH(\MU;\bF_p)} \bF_p^{\otimes_{\THH(\MU;\bF_p)} \bullet+1}  \right )\,.\]
We filter this cobar complex using the Whitehead filtration $\tau_{\ge \bullet}$ as follows
\[
	\lim_{\Delta}\pi_{2*}\left ( \tau_{\ge\bullet}\THH(R;\bF_p)\otimes_{\tau_{\ge\bullet}\THH(\MU;\bF_p)} \bF_p^{\otimes_{\tau_{\ge \bullet}\THH(\MU;\bF_p)} \bullet +1} \right ) \,,
\]
which has associated graded 
\begin{equation}\label{eq: derived cobar}
	\THH_{*}(R;\bF_p)\otimes_{\THH_{*}(\MU;\bF_p)}^{\mathbb{L}}\bF_p^{\otimes_{\THH_*(\MU;\bF_p)}^{\mathbb{L}} \bullet+1} \,.
\end{equation}
Since $\THH_*(R;\bF_p)$ and $\bF_p$ are even flat $\THH_*(\MU;\bF_p)$-modules, each of the iterated K\"unneth spectral sequences collapse at the $\EE_2$-term and we can regard the associated graded in \eqref{eq: derived cobar} as the associated graded of a filtration on the cochain complex 
\[ 
\THH_*(R/\MU^{\otimes \bullet+1};\bF_p)
\]
whose cohomology is $\pi_*\grmot^*\THH(R;\bF_p)$. 
We can identify the cohomology of the derived cobar complex \eqref{eq: derived cobar}
with $H^*(\Gamma, \THH_{*}(R)\otimes_{\THH_{*}(\MU;\bF_p)}^{\mathbb{L}}\bF_p )$. 
This produces the desired spectral sequence 
as the spectral sequence associated to a filtered cochain complex, or alternatively,
as an instance of the May--Ravenel spectral sequence \cite[Theorem A1.3.9]{Rav86}. 
\end{proof}

\begin{lem}\label{collapsing}
Suppose that $R\to \bF_p$ is an augmentation of $\bE_{1}$-$\MU$-algebras, such that $\THH_*(R;\mathbb{F}_p)$ is finite in each degree $*$. If 
\begin{enumerate}
\item \label{it:collapsing1} $\THH_*(R;\bF_p)$ is an even flat $\THH_*(\MU;\bF_p)$-module, and 
\item  \label{it:collapsing2} there is an isomorphism of $\bE_{0}$-$H^*(\Gamma,\bF_p)$-algebras
\begin{equation*}
H^*(\Gamma,\THH_{*}(R;\bF_p )\otimes_{\THH_{*}(\MU;\bF_p )}^{\mathbb{L}}\bF_p )\cong \THH_{*}(R;\bF_p ) \,,
\end{equation*}
\end{enumerate}
then the motivic spectral sequence 
\[ 
\pi_{*}\grmot^{*}\THH(R;\bF_p )\implies \pi_{*}\THH(R;\bF_p )
\]
collapses.
\end{lem}
\begin{proof}
By~\eqref{it:collapsing1} and Proposition~\ref{prop: identification of associated graded of motivic filtration}, there is a May--Ravenel spectral sequence 
\[
	H^*(\Gamma,\THH_{*}(R;\bF_p )\otimes_{\THH_{*}(\MU;\bF_p )}^{\mathbb{L}}\bF_p )\implies \pi_{*}\grmot^{*}\THH(R;\bF_p ) \,.
\]
By~\eqref{it:collapsing2}, the input of the May--Ravenel spectral sequence can be identified with the abutment of the motivic spectral sequence as a graded~$\bF_p$-module so both the May--Ravenel spectral sequence and the motivic spectral sequence collapse. 
\end{proof}

\begin{proposition}\label{algebraic eff map}
For each $R \in  \{\BP,k(n),\bF_p\}$, $\THH_*(R;\bF_p)$ is an even flat~$\THH_*(\MU;\bF_p)$-module. 
\end{proposition}
\begin{proof}
Let $R\to \bF_p$ be the unique augmentation of $\bE_1$-$\MU$-algebras. It suffices to check that 
\[
\Tor_{*}^{\THH_{*}(\MU;\bF_p )}(\THH_{*}(R;\bF_p ),\bF_p )
\] 
is concentrated in even total degrees. This follows from Proposition \ref{prop:Hocschild-MaySSBP}, Proposition~\ref{prop:Hochschild-MaySSkn}, and the fact that $\THH_*(\bF_p)=\bF_p[\mu]$ by~\cite{Bre78,Bok87a}. 
\end{proof}

\begin{cor}\label{cor:even-THH}
The spectrum $\THH(k(n)/\MU)$ is even.
\end{cor}

\begin{proof}
By Proposition~\ref{algebraic eff map}, we know 
\[ \THH_*(k(n);\bF_p)\otimes_{\THH_*(\MU;\bF_p)}^{\mathbb{L}}\bF_p
\]
is concentrated in even total degrees. By the K\"unneth spectral sequence, this implies that 
\[ \THH(k(n)/\MU;\bF_p)=\THH(k(n);\mathbb{F}_p)\otimes_{\THH(\MU;\mathbb{F}_p)}\mathbb{F}_p
\]
is even. Since $|v_n|=2p^n-2$ the $v_n$-Bockstein spectral sequence
\[
\THH_*(k(n)/\MU;\mathbb{F}_p)[v_n]\implies \THH_*(k(n)/\MU)
\]
collapses and the abutment is concentrated in even degrees.
\end{proof}

\begin{proof}[Proof of Proposition \ref{prop:thh-fp}]
This is a special case of~\cite[Proposition 6.1.6]{HW22}, but we nevertheless provide a proof using the setup from this section. By Lemma~\ref{collapsing} and Proposition~\ref{algebraic eff map}, it suffices to observe that there is a preferred isomorphism 
\[ 
H^*(\Gamma,\THH_{*}(\bF_p ) \otimes_{\THH_*(\MU;\bF_p)}^{\mathbb{L}}\bF_p )\cong H^*(\Gamma,\THH_{*}(\bF_p ) \otimes_{\bF_p}\Gamma )\cong \bF_p[\mu]
\]
of $\bE_{\infty}$-$H^*(\Gamma, \bF_p)$-algebras using the fact that $\THH_*(\bF_p)\cong \bF_p[\mu]$~\cite{Bre78,Bok87a}. 
\end{proof}

\begin{proof}[Proof of Proposition \ref{prop:thh-bp}]
By Lemma~\ref{collapsing} and Proposition~\ref{algebraic eff map}, it suffices to observe that there is a preferred isomorphism
\[ 
H^*(\Gamma,\THH_{*}(\BP ;\bF_p ) \otimes_{\THH_*(\MU;\bF_p)}^{\mathbb{L}}\bF_p )\cong \Lambda(\lambda_i: i\ge 1)
\]
of $\bE_3$-$H^*(\Gamma,\bF_p)$-algebras.
\end{proof}

\begin{proof}[Proof of Proposition~\ref{prop:thh-kn}]
By Lemma~\ref{collapsing} and Proposition~\ref{algebraic eff map}, it suffices to observe that 
\begin{align*}
H^*(\Gamma,\THH_*(k(n);\bF_p)\otimes_{\THH(\mathrm{MU};\bF_p)}^{\mathbb{L}}\bF_p)\cong \bF_p[\mu^{p^{n+1}-1}]\otimes \bF_p[\mu]/\mu^{p^{n-1}}\otimes \Lambda(\lambda_{n+1})
\end{align*}
as a $\bE_0$-$H^*(\Gamma,\bF_p)$-algebra and this identification is compatible with the map 
\[
H^*(\Gamma,\THH_*(k(n);\bF_p)\otimes_{\THH(\MU;\bF_p)}^{\mathbb{L}}\bF_p)\to H^*(\Gamma, \THH_*(\bF_p)\otimes_{\THH(\MU;\bF_p)}^{\mathbb{L}}\bF_p)
\]
of $\bE_0$-$H^*(\Gamma,\bF_p)$-algebras. 
\end{proof}

\subsection{The mod $p$ motivic filtrations on $\THH(\MU)$ and $\THH(\BP)$}

In this subsection, we give a brief discussion of $\pi_*\gr^*_{\mot} \THH(\MU) / p$, as well as its retract $\pi_*\gr^*_{\mot} \THH(\BP)/p$.

First, an immediate consequence of Recollection \ref{rec:THHMU} is that there is an isomorphism
\[\THH(\MU) / p \cong \pi_*(\MU/p) \otimes_{\mathbb{F}_p} \Lambda(\lambda_1',\lambda_2',\cdots),\]
with retract
\[\THH(\BP) / p \cong \mathbb{F}_p[v_1,v_2,\cdots] \otimes_{\mathbb{F}_p} \Lambda(\lambda_1,\lambda_2,\cdots).\]
We observe first that similar descriptions apply to mod $p$ motivic associated gradeds.

\begin{proposition}\label{prop:thhbp}
There is a preferred isomorphism
\[
 \pi_{*}\grmot^{\ast}\THH(\BP)/p\cong (\BP/p)_*\otimes \Lambda (\lambda_{i} : i\ge 1) \,,
\]
with bidegrees of indecomposable algebra generators as follows:
\begin{align*}
	\|\lambda_i\| &= (2p^i-1, 1) \,,\\
     \|v_i\| &= (2p^i-2, 0) \,.
\end{align*}
\end{proposition}
\begin{proof}
First, note that since $\MU$ is an $\bE_\infty$-ring there is a splitting 
\[ \pi_*\gr_{\mot}^*\MU/p\to \pi_*\grmot^*\THH(\MU)/p\to \pi_*\grmot^*\MU/p\]
and since $\BP$-splits off of $\MU_{(p)}$ as an $\bE_2$-ring there is a splitting 
\[ \pi_*\gr_{\mot}^*\BP/p\to \pi_*\grmot^*\THH(\BP)/p\to \pi_*\grmot^*\BP/p\]
where the second map in the composite above is the composite
\[\pi_*\grmot^*\THH(\BP)/p\to \pi_*\grmot^*\THH(\MU)/p \to \pi_*\grmot^*\MU/p\to \pi_*\grmot^*\BP/p \,.
\]
We then observe that the spectral sequences 
\[
 \pi_{*}\grmot^{\ast}\BP/(p,v_j,v_{j+1},\cdots )[v_j]\implies  \pi_{*}\grmot^{\ast}\BP/(p,v_{j+1},v_{j+2},\cdots ) \,.
\]
each collapse at the $\EE_1$-page. This implies that the classes $v_1,\cdots ,v_{j-1}$ are permanent cycles in the spectral sequence 
\[
 \pi_{*}\grmot^{\ast}\THH(\BP)/(p,v_j,v_{j+1},\cdots )[v_j]\implies  \pi_{*}\grmot^{\ast}\THH(\BP)/(p,v_{j+1},v_{j+2},\cdots ) \,.
\]
with $\EE_2$-page 
\[
\bF_p[v_1,\cdots ,v_{j-1}]\otimes \Lambda(\lambda_1,\lambda_2,\cdots )[v_j]\,.
\]

Consequently, we observe that the Bockstein spectral sequence
\[
 \pi_{*}\grmot^{\ast}\THH(\BP;\bF_p)[v_1]\implies \pi_{*}\grmot^{\ast}\THH(\BP)/(p,v_2,v_3,\cdots ) \,.
\]
collapses at the $\EE_1$-page. 
Similarly, the spectral sequences 
\[
 \pi_{*}\grmot^{\ast}\THH(\BP)/(p,v_j,v_{j+1},\cdots )[v_j]\implies  \pi_{*}\grmot^{\ast}\THH(\BP)/(p,v_{j+1},v_{j+2},\cdots ) \,.
\]
each collapse at the $\EE_1$-page. 

Moreover, the maps 
\[  \grmot^{\ast}\THH(\BP)/(p,v_{j+1}\cdots )\to  \grmot^{\ast}\THH(\BP)/(p,v_j,v_{j+1},\cdots )\]
are surjective on homotopy groups. Therefore, the sequence is Mittag-Leffler and 
\begin{align*}
\pi_*\lim_{j}(\grmot^{\ast}\THH(\BP)/(p,v_{j+1},v_{j+2},\cdots ))&\cong \lim_j (\BP/(p,v_{j+1},v_{j+2},\cdots ))_*\otimes \Lambda (\lambda_i : i\ge 1)  \\
&\cong (\BP/p)_*\otimes \Lambda (\lambda_{i} : i\ge 1) \,.
\end{align*}
The maps 
\[
\grmot^{\ast}\THH(\BP)/p\to \grmot^{\ast}\THH(\BP)/(p,v_{j+1},v_{j+2},\cdots )
\]
are also increasingly connective as $j$ grows to infinity so the map
\[
\grmot^{\ast}\THH(\BP)/p\to \lim_j\grmot^{\ast}\THH(\BP)/(p,v_{j+1},v_{j+2},\cdots )
\]
is an equivalence. 
\end{proof}

\begin{defin}\label{lambda-thh-cobar}
Using Proposition~\ref{prop:thhbp} and the fact that $\BP$ splits off of $\MU$ as an $\bE_2$-ring by~\cite[Corollary~1.3]{CM15}, we fix explicit cobar complex cocycles 
\[
\sigma^2t_{k} \in \pi_{2p^k}\THH(\MU/\MU^{\otimes 2})/p \,.
\]
for $k\ge 1$ using notation from~\cite[Example~A.2.4]{HW22} providing a choice of generators  
\[ 
\lambda_k \in \grmot^{*}\THH(\MU)/p\,.
\]
We also fix a preferred choice of generator $\mu\in \pi_2\grmot^*\THH(\bF_p)$ detected by the $0$-cocycle   
\[ \sigma^2p \in \THH_2(\bF_p/\MU)\]
in the cobar complex without indeterminacy. 
\end{defin}

\subsection{Mod $(p,v_1,\cdots,v_n)$ Hochschild homology of $\mathbb{F}_p$}
\label{HochschildFp}

A key subtlety in this paper is that $\THH(\bF_p)$ may be viewed as a $\bZ_p$-algebra in two incompatible ways. On the one hand, it is an $\bF_p$-algebra (and therefore a $\bZ_p$-algebra) via the natural unit map $R \to \THH(R)$ that exists for any $\bE_\infty$-ring $R$. On the other hand, it obtains a $\bZ_p$-algebra structure via the trace map 
\[
\bZ_p=\tau_{\ge 0}\TC(\bF_p)_p^\wedge \to \TC(\bF_p)_p^\wedge \to \THH(\bF_p) \,.
\]
In this section, we use these two $\bZ_p$-algebra structures to produce and investigate two different sequences of classes in $\pi_* \grmot^* \THH(\bF_p) / (p,v_1,\cdots,v_n)$.  

We will make use of Dyer--Lashof operations here and in Section~\ref{sec:lambdas}. 
\begin{defin}\label{def:Dyer-Lashof}
Given an $\bE_{2}$-$\bF_p$-algebra there exists a Dyer--Lashoff operation $Q^{(i+1)/2}$  ($Q^{i+1}$ at $p=2$) that raises degree by $(i+1)(p-1)$ \cite[III~\S~1.~Theorem~1.1]{CLM76}. In the setting of graded $\bE_{2}$-$\bF_p$-algebras these operations send classes $x$ in odd degree $i$ and grading $j$ to classes $Q^{i+1}(x)$ (resp. $Q^{(i+1)/2}(x)$) in degree $pi+p-1$ and grading $pj$~\cite[\S~16.1.2]{GKRW19}. 
\end{defin}

First, we record a non-motivic way of understanding the difference between these algebra structures, in terms of homology with $\bF_p$-coefficients. 

\begin{lemma}
Suppose that we identify 
\[H_*\THH(\bF_p) \cong \mathcal{A}_*[\mu]\]
so that the unit map
\[\mathcal{A}_*=H_*\bF_p \to H_*\THH(\bF_p)\]
is a map of $\mathcal{A}_*$-modules.  The element $\mu$ refers to the Hurewicz image of the B\"okstedt class. 

Then, for each $i \ge 1$, the trace map
\[H_*\bZ_p \to H_*\THH(\bF_p)\]
sends:
\begin{enumerate}
\item $\overline{\tau}_i$ in the domain to $\overline{\tau}_i -\mu^{p^i-p^{i-1}} \overline{\tau}_{i-1}$
in the codomain, when $p>2$.
\item $\overline{\xi}_{i+1}$ in the domain to $\overline{\xi}_{i+1} - \mu^{2^i-2^{i-1}} \overline{\xi}_{i}$ in the codomain, when $p=2$.
\end{enumerate}
\end{lemma}

\begin{proof}

We begin by checking the case that $i=1$.  
At odd primes, we must check that $\overline{\tau}_1$ maps to $\overline{\tau}_1 - \mu^{p-1} \overline{\tau}_0$.  To do so, we first observe that the composite map
\[\bZ_p \to \THH(\bF_p) \to \bF_p\]
must be the unique ring map from $\bZ_p \to \bF_p$.  It follows that the trace map
\[H_* \bZ_p \to H_* \THH(\bF_p)\]
must map $\overline{\tau}_1$ to a class $\overline{\tau}_1 + a \mu^{p-1} \overline{\tau}_0$ for some $a \in \mathbb{F}_p$.  Now, any class in $H_*\THH(\bF_p)$ that is in the image of $H_*(\bZ_p)$ must be invariant under the circle action, and we may conclude from \cite[Lemma 3.3]{Rog98} that $a=-1$; cf.~\cite[Proposition~3.2, Proposition 6.1]{BR05}.

The cases when $i>1$ can be iteratively concluded from the case when $i=1$, by applying Dyer--Lashof operations. Specifically, we observe that 
\begin{itemize}
\item $Q^{p^i} \overline{\tau}_i =  \overline{\tau}_{i+1}$ (resp. $Q^{2^{i+1}}\overline{\xi}_{i+1}=\overline{\xi}_{i+2}$ at $p=2$)  for $i\ge 1$, and 
\item $Q^{p^i} \mu^{p^i} = \mu^{p^{i+1}}$ (resp. $Q^{2^{i+1}}(\mu^{2^i})=\mu^{2^{i+1}}$ at $p=2$) for $i\ge 0$
\end{itemize}
by \cite[III.2.2,~III.2.3]{BMMS86} and~\cite[Proposition 5.9]{AR05}. The lemma follows from these formulas and the Cartan formula~\cite[III.1.1]{BMMS86} since the trace map $\bZ_p\to \THH(\bF_p)$ is a map of $\bE_\infty$-algebras and therefore preserves Dyer--Lashof operations.
\end{proof}

\subsubsection{Classes $\varepsilon_i$ and $\overline{\varepsilon}_i$}

We first fix certain classes: 

\begin{convention} \label{cnv:epsilon-def}
For each $n\ge 0$, we fix isomorphisms
\begin{align*}
\pi_*\grev^*\bZ_p/ (p,v_1,\cdots,v_n) & \cong \Lambda(\varepsilon_1,\varepsilon_2,\cdots,\varepsilon_n) \text{ and } \\
\pi_*\grev^* \bF_p / (p,v_1,\cdots,v_n) & \cong \Lambda(\varepsilon_0,\varepsilon_1,\varepsilon_2,\cdots,\varepsilon_n),
\end{align*}
compatible with the mod $p$ reduction map $\bZ_p \to \bF_p$, where 
\[
\| \varepsilon_i \| = (2p^i-1,-1) 
\] 
for each $0\le i\le n$.  We do so such that the natural map
\begin{eqnarray*}
    \pi_* \gr^*_{\ev} \bZ_p / (p,v_1,\cdots,v_n) &=& \pi_*\left(\gr^*_{\ev} \bZ_p \otimes \gr^*_{\ev} \mathbb{S} /(p,v_1,\cdots,v_n)\right) \\
&\to& \pi_*(\gr^*_{\ev} \bZ_p \otimes \gr^*_{\ev} \mathbb{F}_p) \\
&\cong& (\gr^*_{\ev} \bF_p)_*(\gr^*_{\ev} \bZ_p)
\end{eqnarray*}
sends $\epsilon_i$ to $\overline{\tau}_i$.
\end{convention}

\begin{remark}
The map
\[\pi_*\left(\gr^*_{\ev} \mathbb{F}_p \otimes \gr^*_{\ev} \mathbb{S} / (p,v_1,\cdots)\right) \to \pi_*\left(\gr^*_{\ev} \mathbb{F}_p \otimes \gr^*_{\ev} \BP / (p,v_1,\cdots) \right)\cong (\gr^*_{\ev} \mathbb{F}_p)_* (\gr^*_{\ev} \mathbb{F}_p)\]
sends exterior generators in the domain to the elements  $\tau_i$ in the codomain, by \cite[\S~4]{GIKR22} and~\cite[\S~6.2]{Pst23a}.  Postcomposing with the conjugation map in the codomain, we obtain a map sending exterior generators to $\overline{\tau}_i$.  This postcomposition is the mod $p$ reduction of the map mentioned at the end of the preceding convention, so the requirements of the convention can be attained.
\end{remark}

\begin{defin}\label{def:epsilonclasses}
We denote by 
\[
\varepsilon_0,\cdots,\varepsilon_n \in \pi_*\grmot^*\THH(\bF_p) /(p,v_1,\cdots,v_n)
\]
the image of the classes from Convention \ref{cnv:epsilon-def} under the mod $(p,v_1,\cdots,v_n)$ reduction of the unit map
\[
\grev^* \bF_p \to \grmot^*\THH(\bF_p),
\]
where this unit map exists because $\grmot^*\THH(\bF_p) \cong \grev^*\THH(\bF_p)$ by e.g. \cite[Example~6.25]{Pst23}.
We denote by
\[
\overline{\varepsilon}_1,\cdots,\overline{\varepsilon}_n \in \pi_*\grmot^*\THH(\bF_p) /(p,v_1,\cdots,v_n)
\]
the image of the classes from Convention \ref{cnv:epsilon-def} under the mod $(p,v_1,\cdots,v_n)$ reduction of the trace map
\[
\grev^* \bZ_{p} \to \grmot^*\THH(\bF_p)
\]
obtained from applying $\gr_{\ev}^*$ to the trace map $\bZ_p=\tau_{\ge 0} \TC(\bF_p) \to \THH(\bF_p)$. 
\end{defin}

One can write every class in $\pi_* \grmot^*\THH(\bF_p) / (p,v_1,\cdots,v_n)$ in terms of products of $\varepsilon_i$ and $\mu$, according to the isomorphism
\[
\pi_* \grmot^*\THH(\bF_p) / (p,v_1,\cdots,v_n) \cong \bF_p[\mu] \otimes \Lambda(\varepsilon_0,\cdots,\varepsilon_n)\,.
\]
In particular, the classes $\overline{\varepsilon}_i$ are themselves expressible in terms of $\varepsilon_i$.

\begin{proposition} \label{prop:epsilonvsoverlineepsilon}
For each $i \ge 1$, in 
\[\pi_*\grmot^*\THH(\bF_p) / (p,v_1,\cdots,v_n) \cong \bF_p[\mu] \otimes \Lambda(\varepsilon_0,\varepsilon_1,\cdots,\varepsilon_n),\]
there is an equality
\[
\overline{\varepsilon}_i=\varepsilon_i-\mu^{p^i-p^{i-1}} \varepsilon_{i-1} \,.
\]
\end{proposition}

\begin{proof}
It suffices to check these relations in the image of the (injective) map
\[
\pi_* \gr^*_{\mot} \THH(\bF_p)  / (p,v_1,\cdots,v_n) \to \pi_* ((\gr^*_{\ev}\bF_p) \otimes_{\gr^*_{\ev} \mathbb{S}} (\gr^*_{\ev}\THH(\bF_p))),
\]
where they are a consequence of the non-motivic lemma.  Note that the motivic homology of $\fil^*_{\ev} \THH(\bF_p)$ is $\tau$-torsion free, since $\fil^*_{\ev} \THH(\bF_p)$ is a direct sum of bigraded shifts of $\fil^*_{\ev} \bF_p$, and thus motivic computations can be read off from non-motivic computations.

\end{proof}

\begin{warning}
In \cite[Proposition 4.2]{AR12}, Ausoni--Rognes define a class in $V(1)_*\mathrm{THH}(k(1))$ that they denote by $\overline{\varepsilon}_1$. This class is \textbf{not} detected by the class we denote $\overline{\varepsilon}_1$, but rather by the negative of the class we denote $\overline{\varepsilon}_1$. 
In other words, our sign convention differs from~\cite{AR12}.
\end{warning}

We end this subsection by recording the action of the motivic $\sigma$ operator, for future reference.

\begin{proposition} \label{prop:sigmaFp}
In 
\[
\pi_*\grmot^*\THH(\bF_p) / (p,v_1,\cdots,v_n) \cong \bF_p[\mu] \otimes \Lambda(\varepsilon_0,\cdots,\varepsilon_n) \,,
\]
the $\sigma$ operator satisfies the following formulas:
\begin{itemize}
    \item $\sigma(\mu)=0$.
    \item For each $0 \le i \le n$, $\sigma(\varepsilon_i)=\mu^{p^i}$. 
    \item For each $1 \le i \le n$, 
 $\sigma(\overline{\varepsilon}_i)=0$.
\end{itemize}
\end{proposition}

\begin{proof}
The motivic spectral sequence for $\pi_*\grmot^*\THH(\bF_p)/p$ collapses at the $\mathrm{E}_2$-term, so for $\mu$ and $\varepsilon_0$ it suffices to prove the non-motivic statement. 
The fact that $\sigma(\mu)=0$ follows because $\mu$ is in the image of $\sigma$ and $\sigma$ acts as a derivation on $\bF_p$-algebras. The fact that $\sigma(\varepsilon_0)=\mu$ follows from the fact that $p$ is detected by $t\mu$ modulo higher filtration in $\pi_0\TC^{-}(\bF_p)$. The fact that $\sigma(\overline{\varepsilon}_i)=0$ follows from the fact that $\overline{\epsilon}_i$ lifts to a class in $\pi_*\grmot^*\TC(\bF_p)/(p,\cdots ,v_n)$. The fact that $\sigma(\varepsilon_i)=\mu^{p^i}$ then follows from Proposition~\ref{prop:epsilonvsoverlineepsilon}. 
\end{proof}

\begin{rem2}
While the $\epsilon_i$ classes are convenient in the study of Hochschild homology, the $\overline{\epsilon}_i$ classes have the major advantage that they are defined in topological cyclic homology. 

In the remaining sections of this paper, we will thus generally prefer as basis for the bigraded vector space $\pi_*\grmot^*\THH(\bF_p)/(p,\cdots ,v_n)$ the following identification
\[
\pi_*\grmot^*\THH(\bF_p)/(p,\cdots ,v_n) = \bF_p[\mu] \otimes \Lambda (\varepsilon_0,\overline{\varepsilon}_1,\cdots ,\overline{\varepsilon}_n) \,. 
\]
\end{rem2}

\subsubsection{The cap product}

Below, we discuss formulas for the cap product action on 
\[
\pi_*\grmot^*\THH(\bF_p) / (p,v_1,\cdots,v_n) \cong \bF_p[\mu] \otimes \Lambda(\varepsilon_0,\cdots,\varepsilon_n) \,. 
\]

\begin{defin}
Given $S\subset \{1,\cdots ,n\}$, we write 
$\varepsilon_S:=\prod_{s\in S}\varepsilon_s$ where $\varepsilon_{\emptyset}=1$. 
\end{defin}


\begin{lem}\label{lem: cap product formula}
The formula 
\begin{align}\label{eq: cap product formula}
	(\varepsilon_{S} \cdot \mu^j)\cap c_k=\binom{j}{k}\cdot \varepsilon_{S} \cdot \mu^{j-k} 
\end{align}
holds in $\pi_* (\gr^*\THH(\bF_p ) /(p,v_1,\cdots ,v_n))$.  
\end{lem}

\begin{proof}
By restricting the cap product $-\cap c_{k}$ to  
\[
\pi_{*}\gr^*_{\mot} \THH(\bF_p)/(v_i)=\bF_p[\mu]\otimes \Lambda (\varepsilon_i)\,,
\] 
we observe that $\varepsilon_i \cap c_{k}=0$ for each $0\le i\le n$ and $k>0$.  Indeed, the cap product lands in a zero group.

The formula~\eqref{eq: cap product formula} then follows from the fact that $\--\cap c_k$ satisfies the Leibniz rule, so long as we determine $\mu^j \cap c_k$.  To determine $\mu^j \cap c_k$, it suffices to work non-motivically (since the motivic cap product preserves Adams weight).  Since 
\[\THC_*(\bF_p)=\Hom_{\bF_p}(\THH_*(\bF_p),\bF_p)\]
there is a pairing given by evaluation
\[ \langle \-- ,\-- \rangle \colon\thinspace \THH_*(\bF_p)\otimes \THC_*(\bF_p)\longrightarrow \bF_p\]
which satisfies 
\[  
\langle a\cap b, c\rangle = \langle a, b\cdot c \rangle 
\]
for $a\in \THH_*(\bF_p)$ and $b,c\in \THC_*(\bF_p)$. The proof that 
\[  \mu^j \cap c_k=\binom{j}{k}\mu^{j-k}\]
then follows by the computation 
\begin{align*}
\langle \mu^j \cap c_k, c_{j-k} \rangle = & \langle \mu^j , c_k\cdot c_{j-k} \rangle \\
=& \langle \mu^j , \binom{j}{k} c_j \rangle \\
=& \binom{j}{k}\langle \mu^j , c_j \rangle \\
= & \binom{j}{k} . 
\end{align*}

\end{proof}
It will be useful to also record a consequence of this computation in terms of the classes~$\overline{\varepsilon}_i$ from Definition~\ref{def:epsilonclasses}.

\begin{defin}\label{e(S)-f(S)}
Given $S\subset \{1,\cdots ,n\}$, we write $\bar{\varepsilon}_S:=\prod_{s\in S}\bar{\varepsilon}_s$ where $\bar{\varepsilon}_{\emptyset}=1$.  We 
let $f(S)$ denote the integer $\sum_{s\in S}p^{s-1}-p^s$ where $f(\emptyset)=0$.
\end{defin}

\begin{cor}\label{formula for trace classes}
Let $S\subset \{1,\cdots, n\}$. In $\pi_*\grmot^*\THH(\bF_p)/(p,v_1,\cdots ,v_n)$, the formula 
\[
	\bar{\varepsilon}_S\cdot \mu^\ell\cap c_{p^n}=0
\]
holds when $0 \le \ell < p^n+f(S)$.  If it holds for an integer $\ell$, it also holds for all integers in the same congruence  class as $\ell$ modulo $p^{n+1}$.
\end{cor}

\begin{proof}
When $S=\emptyset$, $f(S)=0$ by convention and the formula holds by Lemma~\ref{lem: cap product formula}.

Suppose $S\ne \emptyset$. Writing $\bar{\varepsilon}_S$ in terms of the $\epsilon_i$ and $\mu$ variables 
\[
\bar{\varepsilon}_S=(-1)^{|S|}\prod_{s\in S}\varepsilon_{s-1}\mu^{-f(S)}+(-1)^{|S|-1} \sum_{s' \in S} \varepsilon_{s'}\prod_{s\in S-\{s'\}}\varepsilon_{s-1}\mu^{-f(S-\{s'\})}+\ldots+\varepsilon_{S}
\]
we observe that, since $f(S) \le f(T)\le 0$ for all $T\subset S$, 
\[ 
\bar{\varepsilon}_S\cdot \mu^{\ell}\cap c_{p^n}= 0 \mod p
\]
 when $0\le \ell<p^{n}+f(S)$, by Lemma~\ref{lem: cap product formula}.  
\end{proof}

\subsection{The mod~$(p,v_1,\dots,v_n)$ Hochschild homology $k(n)$}\label{sec:modHH}
In this section, we compute the image of the map
\[ 
f_n : \pi_*\grmot^*\THH(k(n))/(p,v_1,\cdots,v_{n})\to \pi_*\grmot^*\THH(\mathbb{F}_p)/(p,v_1,\cdots,v_{n})
\]
induced by the reduction map $k(n)\to \bF_p$. To accomplish this, we first compute the dimension of the image in Section~\ref{sec:cardinality}, and then the image itself in Section~\ref{sec:image}. 

The basic strategy will be to consider the factorization
\[
\begin{tikzcd}
\pi_*\gr^*_{\mot} \THH(k(n)) / (p,v_1,\cdots,v_n) \arrow{d} \arrow{r}{f_n} & \pi_*\gr^*_{\mot} \THH(\bF_p) / (p,v_1,\cdots,v_n) \\
\pi_*\gr^*_{\mot} \THH(k(n);\mathbb{F}_p) / (p,v_1,\cdots,v_n) \arrow{ur}
\end{tikzcd}
\]

As we will explain in more detail, all maps in this diagram are injections modulo $\lambda_{n+1}$. The map pointing diagonally up and to the right is a $(\pi_*\gr^*_{\ev} \mathbb{F}_p / (p,v_1,\cdots,v_n))$-module map, or in other words a $\Lambda(\epsilon_0,\epsilon_1,\cdots,\epsilon_n)$-module map, making it easy to identify its image.  It is straightforward to show that the vertical downward map is an injection, and to calculate the exact dimensions of its image as a bigraded $\mathbb{F}_p$ vector space.  However, this dimension count alone is not enough to determine the image of $f_n$.

To finish, we note that $f_n$ respects the $\sigma$ operation.  Thus, for any class $x$ in the image of $f_n$, both $x$ and $\sigma(x)$ must be in the image of the map pointing up and to the right.  It will turn out that this last sentence alone is enough to precisely determine the image of $f_n$.

\subsubsection{The cardinality of the image}\label{sec:cardinality}

In this subsection, we compute the dimension of the image of $f_n$ as a bigraded $\bF_p$-module.

Note first that each map in the diagram

\[
\begin{tikzcd}
\pi_*\gr^*_{\mot} \THH(k(n)) / (p,v_1,\cdots,v_n) \arrow{d} \arrow{r}{f_n} & \pi_*\gr^*_{\mot} \THH(\bF_p) / (p,v_1,\cdots,v_n) \\
\pi_*\gr^*_{\mot} \THH(k(n);\mathbb{F}_p) / (p,v_1,\cdots,v_n) \arrow{ur}
\end{tikzcd}
\]
is a $\gr^*_{\mot} \THH(\MU)$-module map, and in particular a $\gr^*_{\ev} \MU$-module map.  We give a preliminary lemma about this module structure:

\begin{lem}\label{lem:acting-trivially}
Each of $p,v_1,\cdots,v_n \in \pi_* \gr^*_{\ev} \MU$ acts trivially on $\pi_* \gr^*_{\mot} \THH(k(n);\mathbb{F}_p)$.  Furthermore, the natural map of $\gr^*_{\ev} \MU$-modules
\[\pi_* \gr^*_{\mot} \THH(k(n)) \to \pi_* \gr^*_{\mot} \THH(k(n);\bF_p)\]
is given by killing $v_n$.
\end{lem}

\begin{proof}
To check the first sentence, note that $\gr^*_{\mot}\THH(k(n);\bF_p)$ is a module over the algebra $\gr^*_{\mot} \THH(\MU;\bF_p)$, so it suffices to check that each of $p,v_1,\cdots,v_n$ act trivially on the latter object.  Since $\gr^*_{\mot}\THH(\MU;\bF_p)$ is a unital $\gr^*_{\mot}\THH(\MU)$-algebra (unlike $\gr^*_{\mot} \THH(k(n);\bF_p)$, which is a unital module but does not have an obvious ring structure), it suffices to note that each of $p,v_1,\cdots,v_n$ is sent to zero under the unit map 
\[\pi_*\gr^*_{\mot}\THH(\MU) \to \pi_*\gr^*_{\mot}\THH(\MU;\bF_p).\]

To check the second statement, we need to observe for each non-negative integer $q$ that 
\[\pi_*\THH(k(n)/\MU^{q+1})\] 
is a free $\mathbb{F}_p[v_n]$-module concentrated in even degrees.  Via a $v_n$-Bockstein spectral sequence, this is equivalent to the already proven fact that
\[\pi_*\THH(k(n)/\MU^{q+1}) / v_n \cong \pi_*\THH(k(n)/\MU^{q+1};\mathbb{F}_p)\]
is concentrated in even degrees.
\end{proof}

\begin{cor}
There is an isomorphism of $\gr^*_{\ev} \MU$-modules between
\[\gr^*_{\mot}\THH(k(n);\bF_p) / (p,v_1,\cdots,v_n)\]
and
\[\gr^*_{\mot}\THH(k(n)) / (p,v_1,\cdots,v_n) \oplus \Sigma^{2p^n-1,-1}\gr^*_{\mot}\THH(k(n)) / (p,v_1,\cdots,v_n),\]
such that the mod $(p,v_1,\cdots,v_n)$ reduction of the map
\[\gr^*_{\mot} \THH(k(n)) \to \gr^*_{\mot} \THH(k(n);\bF_p)\]
is the inclusion of the first summand.
\end{cor}

\begin{proof}
By the Lemma~\ref{lem:acting-trivially}, the mod $(p,v_1,\cdots,v_n)$ reduction of the map
\[\gr^*_{\mot} \THH(k(n)) \to \gr^*_{\mot} \THH(k(n);\bF_p)\]
can be obtained by taking the cofiber of the multiplication by $v_n$ map on
\[\gr^*_{\mot} \THH(k(n)) / (p,v_1,\cdots,v_n),\]
which is a trivial map since we have already killed $v_n$.
\end{proof}

\begin{defin}\label{W}
Let 
\[ W \subset \bZ\times \bZ\]
denote the subset consisting of pairs of integers 
\[ 
 (a_{0}+(2p-1)a_1+\dots +(2p^{n-1}-1)a_{n-1}+2k,-a_{0}-a_1-\ldots -a_{n-1})
\] 
where $k\ge 0$, $k$ is congruent to an element of $\{0,1,\cdots ,p^n-1\}$ modulo $p^{n+1}$, and each $a_i\in \{0,1\}$ for $i=0,1,\cdots,n-1$. 
\end{defin}

\begin{prop}\label{prop:lower-bound}
The image of $f_n$ has dimension $1$ in bidegrees contained in $W$. It has trivial image in all other bidegrees.
\end{prop}
\begin{proof}
The image of the map 
$\pi_{*}\grmot^*\THH(k(n);\bF_p)\to \pi_{*}\grmot^*\THH(\bF_p)$
has dimension exactly $1$ in bidegrees $(2k,0)$ for each $k=0,1,\cdots ,p^n-1\mod p^{n+1}$ by Proposition~\ref{prop:thh-kn}. 
By the preceding proposition, the dimension of the image of $f_n$ is exactly one in bidegree 
\[
(a_{0}+(2p-1)a_1 +\ldots +(2p^{n-1}-1)a_{n-1}+2k,-a_{0}-\ldots -a_{n-1})
\] 
for each $k=0,1,\cdots ,p^n-1\mod p^{n+1}-1$ and each $a_i\in \{0,1\}$ for $i=0,1,\cdots ,n-1$. 
\end{proof}

\begin{defin}\label{U}
Let 
\[ 
U\subset \bZ\times \bZ 
\]
be the subset consisting of pairs of integers 
\[
(a_0+(2p-1)a_1+\ldots +(2p^n-1)a_n+2k,-a_0-\ldots-a_n)
\] 
where either 
\begin{enumerate}
    \item \label{it:1notU} $k\equiv p^n,p^n+1,\cdots ,p^{n+1}-1\mod p^{n+1}$ and $a_i\in \{0,1\}$ for $i=0,1,\cdots , n$, or
    \item \label{it:2notU} $k \equiv 0, 1, \cdots ,p^n-1 \mod p^{n+1}$, $a_n=1$, and $a_i\in \{0,1\}$ for $i=0,1,\cdots ,n-1$.  
\end{enumerate} 
\end{defin}

\begin{cor}\label{cor:codimension}
The subspace 
\[ 
\mathrm{im} f_n \subset  \pi_{\ast}\grmot^{\ast}\THH(\bF_p) / (p,v_1,\cdots,v_n)
\]
has codimension exactly one in bidegrees contained in $U$.  It is surjective in all other bidegrees.  
\end{cor}
\begin{proof}
This follows from Proposition~\ref{prop:lower-bound} and the fact that 
\[ 
\pi_{\ast}\grmot^{\ast}\THH(\bF_p)/(p,v_1,\dots ,v_n)\cong\Lambda(\varepsilon_0,\cdots, \varepsilon_n)\otimes \bF_p[\mu] 
\]
with $\|\mu^k\|=(2k,0)$ and $\|\varepsilon_i\|=(2p^i-1,-1)$. 
\end{proof}

\subsubsection{Determining the image}\label{sec:image}
If $R$ is any $\bE_1$-$\MU$-algebra with an augmentation $R\to \bF_p$ map of $\bE_1$-$\MU$-algebras, then the induced map
\[
\pi_* \grmot^*\THH(R)/(p,v_1,\cdots,v_n) \to \pi_*\grev^*\THH(\bF_p)/(p,v_1,\cdots,v_n)
\]
is closed under the $\sigma$ operator. Here, we use this closure property to determine the image of the map
\[
f_n : \pi_*\grmot^* \THH(k(n)) / (p,v_1,\cdots,v_n) \longrightarrow \pi_*\grev^* \THH(\bF_p) / (p,v_1,\cdots,v_n), 
\]
 in terms of our preferred basis 
\[
\pi_*\grmot^*\THH(\bF_p)/(p,v_1,\cdots ,v_n) \cong \Lambda (\varepsilon_0,\overline{\varepsilon}_1,\cdots ,\overline{\varepsilon}_n)\otimes \bF_p[\mu] 
\]
for $\pi_*\grmot^*\THH(\bF_p)/(p,v_1,\cdots ,v_n)$.


\begin{lem}\label{lem:ker-cap}
The image of $f_n$ is contained in the kernel of the map 
\[\-- \cap c_{p^n} : \pi_* \gr^*_{\mot} \THH(\bF_p) / (p,v_1,\cdots,v_n) \to \THH(\bF_p) / (p,v_1,\cdots,v_n).\]
This kernel is exactly given by the inclusion of $\bF_p$ vector spaces
\[
\ker (\--\cap c_{p^n})\cong \mathbb{F}_p[\mu^{p^{n+1}}]\otimes \bF_p[\mu]/(\mu^{p^n})\otimes \Lambda (\varepsilon_0,\varepsilon_1,\cdots ,\varepsilon_n) \subset \mathbb{F}_p[\mu] \otimes \Lambda(\varepsilon_0,\cdots,\varepsilon_n)\,.
\]
\end{lem}

\begin{proof}
This is immediate from Proposition~\ref{prop:lower-bound} and Lemma~\ref{lem: cap product formula}. 
\end{proof}


\begin{lem}\label{lem: contained in intersection}
The image of $f_n$ is contained in the intersection  
\[\ker ((\sigma(-)\cap c_{p^n})\cap \ker( -\cap c_{p^n} ) \,.
\]
\end{lem}
\begin{proof}
By Lemma~\ref{lem:ker-cap}, the image of $f_n$ is contained in the kernel of $\-- \cap c_{p^n}$.  Furthermore, the image is closed under the $\sigma$-operator, so by the preceding sentence if $x$ is in the image then $\sigma x$ is in the kernel of $\-- \cap c_{p^n}$.
\end{proof}


\begin{defin}\label{def:B}
We define 
\begin{align*}
B:= \coprod_{\overset{S\subset \{1,\cdots ,n\}}{|S|<n}}\{ \bar{\varepsilon}_S\mu^{j}  &:0\le j < p^{n}-1+f(S)\}  
      \}
\end{align*}
as a subset of $\pi_*\grev^*\THH(\bF_p)/(p,\cdots ,v_n)$ using Definition~\ref{e(S)-f(S)}. 
\end{defin}

\begin{lem}\label{lem: height n}
We have 
\[ (\Lambda(\overline{\varepsilon}_1,\dots , \overline{\varepsilon}_n)\oplus\bF_p\{x\mu,x\varepsilon_0 : x\in B\})\otimes \bF_p[\mu^{p^{n+1}}]\subset \ker ( \sigma (-)\cap c_{p^n})\cap \ker ( - \cap c_{p^n}) 
\]
where $B$ is defined as in Definition~\ref{def:B}
\end{lem}

\begin{proof}
This follows from Proposition~\ref{prop:sigmaFp} and Lemma~\ref{lem: cap product formula} as observed in Corollary~\ref{formula for trace classes}. 
\end{proof}
We now use the subspace above to identify 
\[ 
\ker (\sigma(\--)\cap c_{p^n} ) \cap \ker (\--\cap c_{p^n})
\]
in the relevant bidegrees. Note that we determined the action of the $\sigma$ operator in Proposition~\ref{prop:sigmaFp} and we determined $\ker(-\cap c_{p^n})$ in Corollary~\ref{formula for trace classes}. 
\begin{prop}\label{prop:also-codimension-one}
Let $B$ be defined as in Definition~\ref{def:B}. 
\begin{enumerate}
    \item The canonical inclusion
\[
(\Lambda(\overline{\varepsilon}_1,\cdots , \overline{\varepsilon}_n)\oplus\bF_p\{x\mu,x\varepsilon_0 : x\in B\})\otimes \bF_p[\mu^{p^{n+1}}] \hookrightarrow \ker (\sigma (\--) \cap c_{p^n})\cap\ker (\-- \cap c_{p^n}) 
\]
is an isomorphism in bidegrees contained in $W$ as in Definition~\ref{W}. 
\item In bidegrees contained in $W$, the subspace 
\[
\ker (\sigma (\--) \cap c_{p^n})\cap\ker (\-- \cap c_{p^n}) \subset \pi_*\grmot^*\THH(\bF_p)/(p,v_1,\cdots,v_n)
\]
has codimension exactly one in bidegrees contained in $U$ as in Definition~\ref{U} and zero otherwise. 
\end{enumerate}
\end{prop}

\begin{proof}
First, we claim that the subspace 
\[ 
	(\Lambda(\overline{\varepsilon}_1,\dots , \overline{\varepsilon}_n)\otimes \bF_p\{x\mu,x\varepsilon_0 : x\in B\})\otimes \bF_p[\mu^{p^{n+1}}]\subset \pi_{*}\grmot^{*}\THH(\bF_{p})/(p,\cdots ,v_{n})
\]
has codimension exactly one in bidegrees $U$. 

Second, we claim that the subspace 
\[ 
\ker ( \sigma (-)\cap c_{p^n})\cap \ker ( - \cap c_{p^n})\subset  \pi_{*}\grmot^{*}\THH(\bF_{p})/(p,\cdots ,v_{n})
\]
has codimension at most one. 

To prove the first claim, observe that $\varepsilon_0^{a_0}\cdot \ldots \cdot \varepsilon_n^{a_n}\mu^k$ is not in $\ker (-\cap c_{p^n})$ when $k\equiv p^n,\dots ,p^{n+1}-1\mod p^{n+1}$ by Corollary~\ref{formula for trace classes}. In bidegrees
\[
(2k+a_0+a_1(2p-1)+\cdots +a_n(2p^n-1),-a_0-\cdots-a_n)
\] 
where $k\equiv p^n,\cdots ,p^{n+1}-1\mod p^{n+1}$ and $a_i\in \{0,1\}$ for $i=0,\cdots , n$, the $\bF_p$-vector space $\pi_{*}\grev^*\THH(\bF_p)/(p,\cdots ,v_n)$ is spanned by $\varepsilon_0^{a_0}\cdot \ldots \cdot \varepsilon_n^{a_n}\mu^k$ and classes in 
\[
(\Lambda(\overline{\varepsilon}_1,\dots , \overline{\varepsilon}_n)\oplus \bF_p\{\mu,\varepsilon_0x : x\in B\})\otimes \bF_p[\mu^{p^{n+1}}]\,.
\]
This handles Case~\eqref{it:1notU} of Definition~\ref{U}.

When $a_m=0$ for some $0\le m\le n-1$ and then the class 
\[
\varepsilon_0^{a_0}\cdot \varepsilon_1^{a_1}\cdot \ldots \cdot \varepsilon_{m-1}^{a_{m-1}}(\varepsilon_m\mu^{p^n-p^m})\cdot \varepsilon_{m+1}^{a_{m+1}} \cdot \ldots \cdot \varepsilon_{n-1}^{a_{n-1}}\mu^{k}
\]
where $k\not \equiv p^n,\dots ,p^{n+1}-1 \mod p^{n+1}$ is not in $\ker (\sigma (-)\cap c_{p^n})$. This class and classes in 
\[
(\Lambda(\overline{\varepsilon}_1,\overline{\varepsilon}_2,\cdots , \overline{\varepsilon}_n)\oplus \bF_p\{x\mu,\varepsilon_0x : x\in B\})\otimes \bF_p[\mu^{p^{n+1}}]
\] 
span $\pi_*\grev^*\THH(\bF_p)/(p,v_1,\cdots,v_n)$ in this bidegree. This handles Case~\eqref{it:2notU} of Definition~\ref{U} when $a_m=0$ for some $0\le m \le n-1$. 

When $a_{m}\ne 0$ for all $m=0,1,\cdots ,n-1$, then the subspace is zero in this bidegree and it is a subspace of a one dimensional $\bF_{p}$-vector space. This handles Case~\eqref{it:2notU} of Definition~\ref{U}. 

To prove the second claim, we simply observe that in each bidegree in $W$ where the inclusion
\[ 
(\Lambda(\overline{\varepsilon}_1,\overline{\varepsilon}_2, \cdots, \overline{\varepsilon}_n)\oplus\bF_p\{\mu x,\varepsilon_0x : x\in B\})\otimes \bF_p[\mu^{p^{n+1}}]\subset \pi_{*}\grev^*\THH(\bF_{p})
\]
has codimesion $1$, the class spanning the complement is also not in $\ker ( \sigma (-)\cap c_{p^{n}})\cap \ker ( -\cap c_{p^{n}})$. 
The result then follows from Lemma~\ref{lem: height n}. 
\end{proof}

\begin{thm}\label{thm:height-n}
We can identify
\[ 
	\pi_*\grmot^*\THH(k(n))/(p,v_1,\cdots ,v_n)=(\Lambda(\overline{\varepsilon}_1,\cdots ,\overline{\varepsilon}_n)\oplus \bF_2\{  x\mu,x\varepsilon_0 : x\in B\})\otimes \Lambda(\lambda_{n+1})\otimes \bF_2[\mu^{p^{n+1}}]
\]
where $B$ is defined as in Definition~\ref{def:B}. Here the bidegrees are determined by the bidegrees 
\begin{align*}
 \| \overline{\varepsilon}_i \| &= (2p^i-1,-1)\,, \text{ and }\\ 
 \| \mu \| & =(2,0) \,, \\
  \| \lambda_{n+1} \| & =(2p^{n+1}-1,0) 
\end{align*}
and the fact that the image of $f_n$ is 
\[
(\Lambda(\overline{\varepsilon}_1,\cdots , \overline{\varepsilon}_n)\oplus \bF_p\{x\mu,x\varepsilon_0 : x\in B\})\otimes \bF_p[\mu^{p^{n+1}}]\,.
\]
Therefore, we can identify $f_n$ as the composite of the canonical quotient by $(\lambda_{n+1})$ followed by the canonical inclusion as a map of $\mathbb{E}_0$-$\grmot^*\THH_*(\MU)/(p,v_1,\cdots ,v_n)$-algebras. 
\end{thm}
\begin{proof}
By Proposition~\ref{prop:also-codimension-one}, we can identify $\ker (\sigma (\--)\cap c_{p^n})\cap \ker (\--\cap c_{p^n})$ in bidegrees $W$, which is exactly in the bidegrees where we know $\im f_n$ is non-trivial by Proposition~\ref{prop:lower-bound}. 
Proposition~\ref{prop:also-codimension-one} and Corollary~\ref{cor:codimension} demonstrate that in these bidegrees the subspace $\im f_n$ and the subspace $\ker (\sigma (\--)\cap c_{p^n})\cap \ker (\--\cap c_{p^n})$ have the same codimension in 
\[
\pi_*\grmot^*\THH(\bF_p)/(p,v_1,\cdots ,v_n)\,.
\]
The result then follows from Lemma~\ref{lem: contained in intersection}. 
\end{proof}

We now introduce further notation. 
\begin{defin}\label{Kn-Knt}
Fix $n\in \mathbb{N}$. Then we define $K_n$ to be the kernel 
 \[  
 \ker \left ( \sigma : \pi_*\grmot^*\THH(k(n))/(p,\dots ,v_{n+1}) \to \pi_{*+1}\grmot^*\THH(k(n))/(p,\dots ,v_{n+1}) \right )
 \]
 and we define $K_{n}^t$ to be the kernel
  \[  
 \ker \left ( \sigma : \pi_*\grmot^*\THH(k(n))^{tC_p}/(p,\dots ,v_{n+1}) \to \pi_{*+1}\grmot^*\THH(k(n))^{tC_p}/(p,\dots ,v_{n+1}) \right ) \,. 
 \]
We also define $K_{\infty}$ to be the kernel
 \[  
 \ker \left ( \sigma : \pi_*\grmot^*\THH(\bF_p)/(p,\dots ,v_{n+1}) \to \pi_{*+1}\grmot^*\THH(\bF_p)/(p,\dots ,v_{n+1}) \right )
 \]
 and  $K_{\infty}^t$ to be the kernel
   \[  
\ker \left ( \sigma : \pi_*\grmot^*\THH(\bF_p)^{tC_p}/(p,\dots ,v_{n+1}) \to \pi_{*+1}\grmot^*\THH(\bF_p)^{tC_p}/(p,\dots ,v_{n+1}) \right ) \,.  
 \]
These definitions are compatible in the sense that there is a commutative diagram 
 \[
\begin{tikzcd}
  K_n \ar[r,"\varphi_p |_{K_n}"] \ar[d,"f_n/(v_{n+1})|_{K_n}",swap] &   K_n^t \ar[d,"f_n^{tC_p}/(v_{n+1})|_{K_n^t}"] \\ 
    K_\infty \ar[r,"\varphi_p |_{K_{\infty}}",swap] &   K_\infty^t  \,. \\ 
 \end{tikzcd}
 \]
 \end{defin}
The main result we will need from this section is the following. 
\begin{thm}\label{thm:Kn}
There is an identification
 \[
 K_n =\left (
   \Lambda (\overline{\varepsilon}_1,\overline{\varepsilon}_2,\cdots, \overline{\varepsilon}_n)  \oplus \bF_p\{x\mu ,xz_{n+1}: x\in B \} \right ) \otimes \Lambda( \lambda_{n+1})\otimes \bF_p[\mu^{p^{n+1}}]
\]
where $K_n$ is defined in Definition~\ref{Kn-Knt}, 
\[
z_{n+1} =\overline{\varepsilon}_{n+1}\mu+\overline{\varepsilon}_{n}\mu^{p^{n+1}-p^{n}+1}+\cdots +\overline{\varepsilon}_1\mu^{p^{n+1}-p+1}
\]
and $B$ is defined in Definition~\ref{def:B}. Moreover, the map 
\[ 
f_n/(v_{n+1})|_{K_n} : K_n\longrightarrow K_\infty 
\]
is the canonical quotient by $(\lambda_{n+1})$ composed with the canonical inclusion. 
\end{thm}
\begin{proof}
It is clear that $K_n\to K_\infty$ is injective modulo $(\lambda_{n+1})$ and $\sigma \lambda_{n+1}=0$, so it suffices to consider the image of $K_n$ in $K_{\infty}$. 

The commutative diagram of pairs given by 
\[
\begin{tikzcd}
(\MU,\MU) \ar[r] \ar[d] & (\THH(\MU), \THH(\MU)) \ar[d] \\
(\MU,k(n)) \ar[r] & (\THH(\MU),\THH(k(n)))
\end{tikzcd}
\]
induced by the unit maps and identity maps implies that the map 
\[ 
v_{n+1} : \Sigma^{2p^{n+1}-2,0}\grmot^*\THH(k(n))\to \grmot^*\THH(k(n) 
\]
induced by the map 
\[
v_{n+1} : \Sigma^{2p^{n+1}-2,0}\grev^*\MU \to  \grev^*\MU 
\]
from Definition~\ref{def:vi} is nullhomotopic. There is therefore a unique class $\varepsilon_{n+1}$ in bidegree $(2p^{n+1}-1,-1)$ that is compatible with the reduction map $k(n)\to \mathbb{F}_p$ and therefore maps to $\varepsilon_{n+1}$ as fixed in Convention~\ref{cnv:epsilon-def}. This implies that 
\[ \pi_*\grmot^*\THH(k(n))/(p,v_1,\cdots ,v_{n+1})\cong \pi_*\grmot^*\THH(k(n))/(p,v_1,\cdots ,v_n)\otimes \Lambda(\varepsilon_n)
\]
and the map $f_n/(v_{n+1})$ is given by $f_n\otimes \Lambda(\varepsilon_{n+1})$. 

From Theorem~\ref{thm:height-n}, we determine that 
\[ 
K_n\cap \grmot^*\THH(k(n))/(p,v_1,\cdots ,v_n)\{1\}
\] 
is exactly 
\[ \left ( \Lambda (\overline{\varepsilon}_1,\overline{\varepsilon}_2,\cdots,\overline{\varepsilon}_n)\oplus \bF_p\{x\mu : x\in B \} \right )\otimes \Lambda (\lambda_{n+1})\otimes \bF_p[\mu^{p^{n+1}}] \,.
\]
We also observe that $\sigma(z_{n+1})=0$, $z_{n+1}=\varepsilon_{n+1}\mu-\varepsilon_0\mu^{p^{n+1}}$, and 
\[ K_n\cap \grmot^*\THH(k(n))/(p,v_1,\cdots ,v_n)\{\varepsilon_{n+1}\} \]
is exactly 
$
\bF_p\{xz_{n+1} : x\in B \} \otimes \Lambda (\lambda_{n+1})\otimes \bF_p[\mu^{p^{n+1}}]
$. 
\end{proof}

\section{Prismatic cohomology and syntomic cohomology}\label{sec:syntomic}
The goal of this section is to compute the mod $(p,v_1,\cdots,v_{n+1})$ syntomic cohomology of $k(n)$ and prove Theorem~\ref{thm:finiteness}. First, we provide some setup in Section~\ref{sec:setup}. In Section~\ref{sec:classes-and-differentials}, we fix certain permanent cycles and determine certain differentials in (algebraic) $t$-Bockstein spectral sequences. We then compute the Nygaard-completed Hodge--Tate cohomology of $k(n)$ in Section~\ref{sec:Hodge--Tatekn}, the Nygaard-completed prismatic cohomology of $k(n)$ in Section~\ref{sec:prismatickn}, and finally the mod $(p,v_1,\cdots, v_{n+1})$ syntomic cohomology of $k(n)$ in Section~\ref{sec:syntomic-kn}.

\subsection{Setup}\label{sec:setup}
We first define the main objects of study, recalling the relevant parts of~\cite{HRW22}. Throughout this section, fix a connective $\bE_1$-$\MU$-algebra $R$.

\begin{defin}
We define
\[ 
\filmot^*\THH(R)^{tC_p} :=\lim_{\Delta} \tau_{\ge 2*}\left (\THH(R)^{tC_p}\otimes_{\THH(\MU)^{tC_p}}(\MU^{tC_p})^{\otimes_{\THH(\MU)^{tC_p}}\bullet+1} \right )_p^{\wedge} \,,
\]
\[ 
\filmot^*\TC^{-}(R)_p^{\wedge} :=\lim_{\Delta} \tau_{\ge 2*}\left (\TC^{-}(R)\otimes_{\TC^{-}(\MU)}(\MU^{hS^1})^{\otimes_{\TC^{-}(\MU)}\bullet+1} \right )_p^{\wedge}  \,,
\]
and 
\[ 
\filmot^*\TP(R)_p^{\wedge}:= \lim_{\Delta} \tau_{\ge 2*}\left ( \TP(R)\otimes_{\TP(\MU)}(\MU^{tS^1})^{\otimes_{\TP(\MU)}\bullet+1}\right)_p^{\wedge} \,.
\]
\end{defin}
\begin{remark}
There is a map of  cyclotomic spectra
\[
\THH(\MU)\longrightarrow\THH(\MU/\MW)
\]
which is an eff map of $\mathbb{E}_{\infty}$-algebras where $\THH(\MU/\MW)$ is even. As in the proof of \cite[Theorem~4.2.9]{HRW22}, this implies that there are maps of filtered spectra
\[ 
\can, \varphi : \filmot^*\TC^{-}(R)_p^{\wedge}\longrightarrow \filmot^*\TP(R)_p^{\wedge}
\]
that converge to the canonical and Frobenius map respectively from~\cite[Corollary~1.5]{NS18}.\footnote{Here we write $\varphi$ as shorthand for $\varphi_p^{hS^1}$.} 
\end{remark}

\begin{defin}
We define
\[
\filmot^*\TC(R)_p^{\wedge}:=\mathrm{eq} \left ( \begin{tikzcd} 
\filmot^*\TC^{-}(R)_p^{\wedge}
\ar[r,"\can",shift left=1ex] \ar[r,"\varphi", swap, shift right=1ex] & \filmot^*\TP(R)_p^{\wedge}
\end{tikzcd} \right ) \,. 
\] 
We write 
\[ \grmot^*F(R):= \filmot^*F(R)/\filmot^{*+1}F(R)
\]
for $F\in \{\THH^{tC_p},\TC^{-},\TP,\TC\}$. We refer to $\pi_*\grmot^*\THH(R)^{tC_p}$ as \emph{Nygaard completed Hodge--Tate cohomology}, $\pi_*\grmot^*\TP(R)$ as \emph{Nygaard completed prismatic cohomology}, and $\pi_*\grmot^*\TC(R)$ as \emph{syntomic cohomology}. 
\end{defin}
By construction, $\filmot^*F(R)$ is complete for each $F\in \{ \TC^{-},\TP,\TC\}$ and we can consider the corresponding spectral sequence
\[
\pi_*\grmot^*F(R) \implies \pi_*F(R)\,,
\] 
which we call the \emph{motivic spectral sequence}. 

By filtering the chain complexes 
\[
\pi_{2*}\left (\TC^{-}(R)\otimes_{\TC^{-}(\MU)}(\MU^{hS^1})^{\otimes_{\TC^{-}(\MU)}\bullet+1} \right )_p^{\wedge}  \,,
\]
\[
\pi_{2*}\left ( \TP(R)\otimes_{\TP(\MU)}(\MU^{tS^1})^{\otimes_{\TP(\MU)}\bullet+1}\right)_p^{\wedge}  \,,
\]
and 
\[
\pi_{2*}\left ( (\THH(R)^{tC_p})^{hS^1}\otimes_{(\THH(\MU)^{tC_p})^{hS^1}}((\MU^{tC_p})^{hS^1})^{\otimes_{(\THH(\MU)^{tC_p})^{hS^1}}\bullet+1}\right)_p^{\wedge} 
\]
using the skeletal/Tate filtrations, we produce a diagram of spectral sequences 
\[
\begin{tikzcd}
\pi_*\grmot^*\THH(R)[t^{\pm 1}] \ar[r,Rightarrow] & \pi_*\grmot^*\TP(R)  \\ 
\pi_*\grmot^*\THH(R)[t] \ar[u] \ar[d] \ar[r,Rightarrow] & \pi_*\grmot^*\TC^{-}(R)  \ar[d,"\varphi"] \ar[u,"\can",swap]\\ 
\pi_*\grmot^*\THH(R)^{tC_p}[t] \ar[r,Rightarrow] & \pi_*\grmot^*\TP(R)  \,. 
\end{tikzcd}
\]
We refer to the middle spectral sequence as the \emph{$t$-Bockstein spectral sequence} and the top spectral sequence as the \emph{periodic $t$-Bockstein spectral sequence}. We also refer to the (compatible) filtrations on the middle and top spectral sequences as the \emph{Nygaard filtration}. 

\subsection{Classes in syntomic cohomology and differentials}\label{sec:classes-and-differentials}

In this section, we record some useful information about the syntomic and Nygaard filtered prismatic cohomology of $\BP$. The goal of Section~\ref{sec:lambdas} is to produce certain permanent cycles called $\lambda_k$ in the (periodic) $t$-Bockstein spectral sequence. The goal of Section~\ref{sec:diff-on-powers-t} is to determine certain differentials on powers of $t$ in the (periodic) $t$-Bockstein spectral sequence. 

\subsubsection{The classes $\lambda_k$}\label{sec:lambdas}
First, observe that the $\bE_4$-ring map 
\[
\BP\to \tau_{\le 0}\BP=\bZ_{(p)}
\] 
is $(2p-2)$-connective, and consequently induces a $(2p-1)$-connective map 
\[
\TC(\BP)_p^{\wedge}\to\TC(\bZ_{(p)})_p^{\wedge}=\TC(\bZ_{p})_p^{\wedge}
\]
by~\cite[Proposition~10.9]{BM94} and~\cite{Dun97}. By~\cite{BM94,Rog00}, there is a class 
\[
\lambda_1\in \pi_{2p-1}\TC(\bZ_{p})^{\wedge}_p
\] 
mapping to a generator $\lambda_1\in \pi_{2p-1}(\THH(\bZ_p)/p)=\bF_p$. By connectivity, there is a corresponding class $\lambda_1\in \pi_{2p-1}\TC(\BP)^{\wedge}_p$ mapping to $\lambda_1\in \pi_{2p-1}\TC(\bZ_{(p)})_p^{\wedge}$ and to a generator $\lambda_1$ of $\pi_{2p-1}\THH(\BP)/p=\bF_p$. 

\begin{proposition}\label{prop:lambda1}
The class $\lambda_1$ is detected by a class in bidegree 
$(2p-1,1)$ in 
\[
\pi_{*}\grmot^*\TC(\BP)_p^{\wedge}\,,
\] mapping to the generator of  
\[
\pi_{*}\grmot^*\THH(\BP)/p=\bF_p
\]
in the same bidegree.
\end{proposition}

\begin{proof}
Observe that the motivic spectral sequence 
\[ 
\pi_*\grmot^*\TC(\bZ_p)_p^{\wedge}\implies \pi_*\TC(\bZ_p)_p^{\wedge}
\]
is concentrated on lines $0$, $1$, and $2$ and consequently $\lambda_1$ must be detected in Adams weight $1$ because it is in odd dimension. The motivic spectral sequence 
\[ 
\pi_*\grmot^*\THH(\bZ_p)/p\implies \pi_*\THH(\bZ_p)/p
\]
is concentrated on lines $0$ and $1$ so the class detecting $\lambda_1$ maps to a generator of 
\[
\pi_{*}\grmot^*\THH(\bZ_p)/p=\bF_p
\] in bidegree $(2p-1,1)$. Since the maps 
\[ 
\grmot^*\TC(\BP )\to \grmot^*\TC(\bZ_{(p)})
\]
and 
\[ 
\grmot^*\THH(\BP )\to \grmot^*\THH(\bZ_{(p)})
\]
cannot lower Adams weight, we determine that $\lambda_1$ is detected in bidegree $(2p-1,1)$ in $\pi_{*}\grmot^*\TC(\BP)_p^{\wedge}$ and it maps to a generator $\lambda_1$ of 
$
\pi_{*}\grmot^*\THH(\BP)/p=\bF_p
$
in the same bidegree. 
\end{proof}

\begin{remark}
We will abuse notation and denote by $\lambda_1$ the class in $\pi_*\gr^*_{\mot} \TC(\BP)_p^{\wedge}$ which detects $\lambda_1 \in \pi_*\TC(\BP)_p^{\wedge}$ according to Proposition~\ref{prop:lambda1}.  
\end{remark}

\begin{remark}
Since $\BP_p^{\wedge}$ is a retract of $\MU_{p}^{\wedge}$ as an $\bE_2$ ring by \cite[Corollary~1.3]{CM15} and $\TC(\MU)_p^{\wedge}\simeq \TC(\MU_p^{\wedge})_p^{\wedge}$, there is also a preferred class 
$\lambda_1\in \pi_{2p-1}\TC(\MU)_p^{\wedge}$ detected by a class in bidegree $(2p-1,1)$ of $\pi_{*}\grmot^*\TC(\MU)_p^{\wedge}$ and mapping to a preferred generator of $\pi_{*}\grmot^*\THH(\MU)/p$ in bidegree $(2p-1,1)$. We abuse notation and simply write $\lambda_1$ to denote any of these elements.  This is compatible, under the natural trace map from $\TC$ to $\THH$, with the notations from previous sections of the paper. 
\end{remark}

\begin{defin}\label{def:lambdas}
For $p>2$, we inductively define  
\[
    \lambda_k:=Q^{p^{k-1}}(\lambda_{k-1})\in \pi_{*}\grmot^{*}\TC(\MU)/p
\] 
in bidegree $(2p^{k}-1,1)$ where $\lambda_1$ is as defined in Proposition~\ref{prop:lambda1} and $Q^{p^{k-1}}$ is defined as in Definition~\ref{def:Dyer-Lashof}. At $p=2$, we inductively define $\lambda_k:=Q^{2^{k}}(\lambda_{k-1})$. 
\end{defin}

\begin{lem}
The class $\lambda_k\in \pi_{*}\grmot^{*}\TC(\MU)/p$ from Definition~\ref{def:lambdas} maps to the class $\lambda_k\in \pi_{*}\grmot^{*}\THH(\MU)/p$ from Definition~\ref{lambda-thh-cobar} for each $k\ge 1$.  
\end{lem}

\begin{proof}
The trace map is compatible with the Dyer--Lashof operation $Q^{p^{k}}$ (respectively $Q^{2^{k+1}}$ at $p=2$) from Definition~\ref{def:Dyer-Lashof}, so it suffices to show that inductively 
\[
Q^{p^{k}}(\lambda_k)=\lambda_{k+1}\in \pi_{2p^{k+1}-1}\grmot^{*}\THH(\MU)/p
\] 
(resp. $Q^{2^{k+1}}(\lambda_k)=\lambda_{k+1}$ at $p=2$). This follows from the discussion in Remark~\ref{rem:Dyer--Lashof}. 
\end{proof}

\begin{cor}\label{cor:lambdas-perm}
For each $k\ge 1$, the class
\[
\lambda_k \in \pi_*\gr_{\mot}^*\THH(\BP)/p\subset  \pi_*\gr_{\mot}^*\THH(\BP)/p[t]\subset \pi_*\gr_{\mot}^*\THH(\BP)/p[t,t^{-1}]
\]
is a permanent cycle in the $t$-Bockstein and periodic $t$-Bockstein spectral sequences. 
\end{cor}
\begin{proof}
The classes $\lambda_k$ are in the image of trace map for each $k\ge 1$.
\end{proof}

\subsubsection{Differentials on powers of $t$}\label{sec:diff-on-powers-t}
We now determine differentials on powers of $t$ in the periodic $t$-Bockstein spectral sequence. 

\begin{prop}\label{prop:image-of-alpha-1}
We can identify 
\[ \lambda_1=\sigma^2\alpha_1\in\pi_{2p-1}\THH(\BP)_p^{\wedge} 
\]
where $\alpha_1$ is a generator of the abelian group $\pi_{2p-3}\mathbb{S}_{(p)}=\bF_p$. Consequently, there is a class  $\alpha_1 \in \pi_{2p-3}\gr_{\mot}^{*}\TC^{-}(\MU)_p^{\wedge}$ detected by $t\lambda_1$ in the $t$-Bockstein spectral sequence and detected by the $1$-cocycle 
\[ 
[t_1] \in \pi_*\TC^{-}(\MU/\MU^{\otimes 2})_p^{\wedge}
\]
in the motivic spectral sequence that is in the image of the $1$-cocycle 
\[ [t_{1} ] \in \pi_*\MU  \otimes  \MU \,,
\]
which detects $\alpha_1\in \pi_{2p-3}\grmot^{*}\mathbb{S}$
in bidegree $(2p-3,1)$.
\end{prop}
\begin{proof}
    By \cite{BM94,Rog00}, the class $[t_{1}]$ maps to $t\lambda_1$ in the approximate $t$-Bockstein spectral sequence 
    \[
    \grmot^*\THH(\bZ)_p^{\wedge}[t]/t^2 \implies \pi_*\grmot^*F(S^3_+,\THH(\bZ)_p^{\wedge})^{S^1}
    \]
    computing $\pi_*\grmot^*F(S^3_+,\THH(\bZ)_p^{\wedge})^{S^1}$ in bidegree~$(2p-3,1)$. 
Moreover, we know that $[t_{1}]$ maps to zero in $\pi_*\grmot^*\THH(\BP)_{p}^{\wedge}$ by connectivity of the map 
\[
\grmot^*\THH(\BP)_{p}^{\wedge}\to \grmot^*\THH(\bZ_{(p)})_{p}^{\wedge}
\]
so it must be in positive Nygaard filtration. The map 
\[ \grmot^*\TC^{-}(\BP;\bZ_p)\to  \grmot^*\TC^{-}(\bZ_{(p)};\bZ_p)
\]
cannot strictly raise Nygaard filtration. Since $[t_1]$ is a cocycle that detects $\alpha_1$ and $t\sigma^2x=x$ by~\cite[Lemma A.4.1]{HW22}, we determine that $t\lambda_1=\alpha_1$ and therefore we can choose our class $\lambda_1$ to satisfy $\lambda_1=\sigma^2\alpha_1$. 
\end{proof}

\begin{remark}\label{rem:Dyer--Lashof}
For simplicity, we focus on odd primes $p$. The same argument applies at $p=2$ by replacing $Q^{p^{k}}$ with $Q^{2^{k+1}}$. Note that 
\[ Q^{p^{k}}\sigma^{2}t_{k}=\sigma^{2}Q^{p^{k}}t_{k}=\sigma^{2}t_{k+1} \in \pi_{2p^{k+1}}\THH(\MU/\MU^{\otimes 2})/p
\]
by Lemma~\ref{lem:sigma-squared-commutes-with-power-operations} and~\cite[Theorem~6]{Koc74}. We can therefore alternatively define cocycles representing $\lambda_{k+1}$ by iterated Dyer--Lashof operations 
\[  Q^{p^{k}}\lambda_{k}=\lambda_{k+1}\,.\]
This has an explicit description as the $p$-th power operation in the cosimplicial $\bE_{\infty}$-$\bF_{p}$-algebra 
\[
 \pi_{*}(\THH(\MU/\MU^{\otimes \bullet+1})/p)
\]
as we now explain. Consider the diagram of non-unital $\bE_{\infty}$-$\bF_{p}$-algebras 
\[
\begin{tikzcd}
 \Sigma^{-1}\pi_*(\THH(\MU/\MU^{\otimes 2})/p) \ar[r] \ar[d]&  0 \ar[d]\ar[r,shift right=1ex]\ar[r,shift left=1ex] & \pi_*(\THH(\MU/\MU^{\otimes 2})/p) \ar[d] \\
F \ar[r] & \pi_*(\MU/p)\ar[r,shift right=1ex]\ar[r,shift left=1ex] 
 & \pi_*(\THH(\MU/\MU^{\otimes 2})/p)
\end{tikzcd}
\]
where $F$ denote the limit of the bottom line. Choosing a cobar complex representative
\[ \sigma^{2}t_{1}\in \pi_{*}\THH(\MU/\MU^{\otimes 2})/p\]
for $\lambda_{1}$ the $\mathbb{E}_{2}$-Dyer--Lashof operation $Q^{p^{k}}$ can be computed as the $\mathbb{E}_{1}$-Dyer--Lashof operation on $ \Sigma^{-1}\pi_*(\THH(\MU/\MU^{\otimes 2})/p$ and in this case it is exactly the Steenrod operation given by taking the $p$-th power. We can therefore choose a cocycle representative 
\[ 
(\sigma^{2}t_{1})^{p^{k}}\in \pi_{*}\THH(\MU/\MU^{\otimes 2})/p
\]
for $\lambda_{k+1}$ in the cobar complex computing $\pi_*\grmot^{*}\THH(\MU)/p$. 
\end{remark}

\begin{prop}\label{prop:differential-bp}
In the spectral sequence 
\begin{align}\label{eq:per-t-bockstein}
\pi_{*}\gr_{\mot}^{*}\THH(\mathrm{BP})/(p,v_1,\cdots ,v_{n+1})[t^{\pm 1}]\implies  \pi_{*}\gr_{\mot}^{*}\mathrm{TP}(\mathrm{BP})/(p,v_1,\cdots ,v_{n+1})
\end{align}
we have that 
\begin{enumerate}
    \item the class $t^{p^\ell}$ survives until the $E_{p^{\ell+1}}$-page for every $\ell\ge 0$, and
    \item there are differentials
\[  
d_{p^{j+1}}(t^{p^{j}})=t^{p^{j+1}+p^{j}}\lambda_{j+1}
\]
for each $0\le j\le n$. 
\end{enumerate} 
\end{prop}

\begin{proof}
We know that the differential in the cobar complex 
\[ 
\pi_*\TP(\MU^{\otimes \bullet+1})
\]
satisfies 
\[
d(t^{p^{j-1}})=\eta_{R}(t^{p^{j-1}})-t^{p^{j-1}}=-t^{p^{j+1}}t_{1}^{p^j} \mod (p,v_1, t^{p^j+2p^{j-1}})
\] 
by~\cite[Theorem~3.13]{Wil82} (cf.~\cite[Remark~6.2.5]{HRW22}). Moreover, we know that $t^{p^{j-1}}$ maps to the class with the same name in the cobar complex 
\[
\pi_*\TP(\BP/\MU^{\otimes \bullet+1})/(p,v_1,\cdots ,v_n) \,.
\]
This implies that $d_r(t^{p^{\ell}})=0$ for $1\le r< p^{\ell+1}$ in the (periodic) $t$-Bockstein spectral sequence \eqref{eq:per-t-bockstein} for all $\ell\ge 0$. To see this, note that the differential in the (periodic) $t$-Bockstein spectral sequence is computed by 
\[ 
d_r(t^{p^j})=d(t^{p^j}) \mod (t^{p^j+r+1})
\]
since it is the spectral sequence associated to a filtered chain complex, where $d$ denotes the differential in the cobar complex. The same formula implies the differentials 
\[ 
d_{p^{j+1}}(t^{p^j})=-t^{p^{j+1}+p^j}\lambda_{j+1}  \mod (p,v_1)
\]
for $0\le j\le n$ as we now explain. 
Definition~\ref{def:lambdas},  Remark~\ref{rem:Dyer--Lashof}, and  Proposition~\ref{prop:image-of-alpha-1} imply that $t^{p^j}\lambda_{j+1}$ is represented in the cobar complex by the image $(t\sigma^2 t_1)^{p^j}$ of 
\[ 
[t_1^{p^j}] \in \pi_*\MU^{\otimes 2}\,. 
\] 
Together, Corollary~\ref{cor:lambdas-perm} and Proposition~\ref{prop:differential-bp} imply that $t^{p^{j+1}+p^j}\lambda_{j+1}$ is not a boundary of a $d_r$-differential for $1\le r<p^j$ 
and $0\le j\le n$. 
\end{proof}

\subsection{Hodge--Tate cohomology of $k(n)$}\label{sec:Hodge--Tatekn}

The main theorem of this subsection is the following calculation of mod $(p,v_1,\cdots,v_{n+1})$ Nygaard completed Hodge--Tate cohomology:

\begin{thm}\label{thm:Hodge-Tate}
The commutative diagram 
\[
\begin{tikzcd}
    K_n  \ar[r,"\varphi_p |_{K_n}"]\ar[d,"f_n/(v_{n+1})|_{K_n}",swap] & K_n^t \ar[d,"f_n^{tC_p}/(v_{n+1}|_{K_n^t}"]\\ 
    K_{\infty}  \ar[r,"\varphi_p |_{K_\infty}",swap] & K_{\infty}^t 
\end{tikzcd} 
\]
from Definition~\ref{Kn-Knt}, 
may be identified with the commutative diagram of $\Lambda(\lambda_{n+1})$-modules 
\[
\begin{tikzcd}
   ( L \oplus M)  \otimes \Lambda( \lambda_{n+1})\otimes \bF_p[\mu^{p^{n+1}}] \ar[r] \ar[d] & (
  L \oplus M ) \otimes \Lambda( \lambda_{n+1})\otimes \bF_p[\mu^{\pm p^{n+1}}] \ar[d] \\
L \otimes \Lambda (\bar{\varepsilon}_{n+1})\otimes \bF_p[\mu^{p^n+1}] \ar[r] &  L\otimes \Lambda(\overline{\varepsilon}_{n+1}) \otimes \bF_p[\mu^{\pm p^{n+1}}] 
\end{tikzcd}
\]
where $L:=\Lambda (\overline{\varepsilon}_1,\cdots ,\overline{\varepsilon}_n)$, $M:=\bF_p\{x\mu ,xz_{n+1}: x\in B \}$, and $B$ is defined as in Definition~\ref{def:B}. Here the vertical maps are given by the canonical quotient by $(\lambda_{n+1})$ composed with the canonical inclusion, the bottom horizontal map is the localization map and the top map is the restriction of the localization map.   
\end{thm}

As we will explain, this implies:

\begin{thm}\label{thm:Hodge-Tateresults}
The spectral sequence beginning with
\[\left(\pi_*\gr^*_{\mot} \THH(k(n))^{tC_p} / (p,v_1,\cdots,v_{n+1})\right)[t]\] and converging to $\pi_*\gr^*_{\mot} \TP(k(n)) / (p,v_1,\cdots,v_{n+1})$ collapses at the $\EE_2$-page. All permanent cycles are on the $0$-line. Moreover, the $0$-line is 
\[K_n^t= \left (
   \Lambda (\overline{\varepsilon}_1,\cdots ,\overline{\varepsilon}_n)\oplus \bF_p\{x\mu ,xz_{n+1}: x\in B \} \right ) \otimes \Lambda( \lambda_{n+1})\otimes \bF_p[\mu^{\pm p^{n+1}}] \]
   where $B$ is defined in Definition~\ref{def:B}. 
\end{thm}

\begin{rmk}
The collapse at the $\EE_2$-page determined in
Theorem~\ref{thm:Hodge-Tateresults} admits a comparatively simple proof in the special case that $k(n)$ is an $\bE_1$-$\BP \langle n \rangle$-algebra for an $\bE_3$-$\MU$-algebra form of $\BP \langle n \rangle$.  In such a case, the spectral sequence of Theorem~\ref{thm:Hodge-Tateresults} is a module over the corresponding spectral sequence computing 
\[
\pi_* \grmot^* \TP(\BP\langle n \rangle) / (p,v_1,\cdots,v_{n+1})\,,
\]
where $t$ dies via a $d_1$ differential.
\end{rmk}

We proceed in the remainder of this subsection with proofs of Theorems \ref{thm:Hodge-Tate} and \ref{thm:Hodge-Tateresults} in general. We would like to thank Tristan Yang in particular for several discussions about these two results.

\subsubsection{Integral Hodge--Tate cohomology of $k(n)$}

Here, we prove the following:

\begin{proposition} \label{integral-Hodge-Tate-is-even}
The bigraded homotopy groups \[\pi_*\gr^*_{\mot} \THH(k(n))^{tC_p}\] are concentrated on the $0$-line and $1$-line.  Thus, the motivic spectral sequence computing $\pi_*\THH(k(n))^{tC_p}$ collapses with no differentials or extension problems.
\end{proposition}

The starting point for this result is that $\pi_*\gr^*_{\mot} \THH(k(n))$ is concentrated on the $0$-line and $1$-line, as proved in the previous sections. It follows formally from the spectral sequence
\[
\pi_*\grmot^*\THH(k(n)) [t,t^{-1}] \implies \pi_*\grmot^*\TP(k(n))
\]
that $\pi_*\grmot^*\TP(k(n))$ is concentrated on the $0$-line and $1$-line, and therefore that 
\[
\pi_* \grmot^*\THH(k(n))^{tC_p}
\]
is concentrated on the $0$-line, $1$-line, and ($-1$)-line.

Ruling out classes on the $(-1)$-line would be automatic if every term in the Hopf algebroid
\[\pi_*\THH(k(n)/\MU^{\otimes \bullet+1})^{tC_p}\] were even, but this does not seem clear. To proceed, we will use a slightly different Hopf algebroid that also computes $\pi_*\grmot^* \THH(k(n))^{tC_p}$.

\subsubsection{Constructing an alternative cobar complex}

Choose an $\mathbb{E}_2$-algebra map $\mathbb{S}[y] \to \MU$ sending $y$ to $v_n$, which may then be extended to an $\mathbb{E}_2$-$\MU$-algebra map
\[
\MU[y] \to \MU\,.
\]
The authors would like to thank Piotr Pstr\k{a}gowski for discussions about the result below.  

\begin{proposition}\label{prop:change-of-eff-cover}
The $\mathbb{E}_{\infty}$-$\MU$-algebra unit map $\MU \to \MU[y]$ induces a map
\[
\pi_*\THH(k(n)/\MU^{\otimes \bullet+1}) \to \pi_*\THH(k(n)/\MU[y]^{\otimes \bullet+1})
\]
that is a quasi-isomorphism of cobar complexes.  In particular, the homology of the cosimplicial abelian group
\[\pi_*\THH(k(n)/\MU[y]^{\otimes \bullet+1})\]
computes $\pi_*\gr^*_{\mot} \THH(k(n)).$
\end{proposition}

\begin{proof}
Since the composite 
\[\mathrm{id}_{\MU} : \MU \to \MU[y] \to \MU \] 
is a map of $\bE_2$-$\MU$-algebras, the induced map
\[\fil_{\mot}^*\THH(k(n))\to  \fil_{\ev/\THH(\MU[y])}^*\THH(k(n))\to \fil_{\mot}^*\THH(k(n))
\]
is the identity by functoriality of the even filtration~\cite[Definition 1.1]{Pst23}, which agrees with the even filtration considered in this paper by~\cite[Theorem~7.3]{Pst23}.  It will therefore suffice to prove that the map 
\[\fil_{\ev/\THH(\MU[y])}^*\THH(k(n))\to \fil_{\mot}^*\THH(k(n)) \]
is an equivalence. 

In the language of \cite[Definition~2.15,~Definition~6.16]{Pst23}, we observe that 
\[
\THH(\MU/\MU[y])=\THH(\MU) \otimes_{\THH(\MU[y])} \MU[y] 
\] 
has even homotopy groups so it is a homologically even $\THH(\MU[y])$-module by~\cite[Lemma~2.27]{Pst23}, and therefore $\THH(\MU)$ is a homologically even $\THH(\MU[y])$-module by \cite[Lemma~6.21]{Pst23}. By~\cite[Theorem~6.15]{Pst23} applied with $N=\THH(k(n))$ and $f: R\to S$ given by the map $\THH(\MU[y])\longrightarrow \THH(\MU)$ of $\bE_1$-rings discussed above, the result follows.
\end{proof} 

\subsubsection{Deduction} 
Proposition~\ref{prop:change-of-eff-cover} directly states that the cohomology of the cosimplicial abelian group 
\[
\pi_*\THH(k(n)/\MU[y]^{\otimes \bullet+1})
\] 
computes $\pi_*\gr^*_{\mot} \THH(k(n))$.  However, since the comparison of cosimplicial objects induced by the $\mathbb{E}_{\infty}$-ring map $\MU \to \MU[y]$ is circle equivariant, it follows also that we may compute
\[\pi_* \gr^*_{\mot} \THH(k(n))^{tC_p}\]
using the cohomology of the cosimplicial abelian group
\[\pi_* \THH(k(n)/\MU[y]^{\otimes \bullet+1})^{tC_p}.\]

Proposition~\ref{integral-Hodge-Tate-is-even} will follow from:

\begin{proposition}\label{prop:evenness}
For each integer $k \ge 1$, 
\[\THH(k(n)/\MU[y]^{\otimes k})^{tC_p}\]
is even.
\end{proposition}

\begin{proof}
Write $\MU[y]^{\otimes k}$ as $\MU^{\otimes k}[y_1,y_2,\cdots,y_k]$.  The $\bE_{\infty}$ inclusion of the first tensor factor 
\[\MU[y_1] \to \MU[y]^{\otimes k},\]
together with the $\mathbb{E}_{\infty}$-augmentation $\MU[y_1] \to \MU$ sending $y_1$ to $0$, gives an isomorphism of $S^1$-equivariant spectra

\begin{eqnarray*}    
\THH\left(k(n)/\MU[y]^{\otimes k}\right) \otimes_{\MU[y_1]} \MU &\cong& \THH\left( \left(k(n) \otimes_{\MU[y_1]} \MU\right) / \MU^{\otimes k}[y_2,\cdots,y_k]\right) \\
&\cong& \THH(\mathbb{F}_p / \MU^{\otimes k}[y_2,\cdots,y_k]) \,. 
\end{eqnarray*}

Our goal is then to show that
\begin{enumerate}
    \item \label{it1-HT} $\THH(\mathbb{F}_p / \MU^{\otimes k}[y_2,\cdots,y_k])^{tC_p}$ is even.
    \item \label{it2-HT} $\THH(k(n)/\MU^{\otimes k} [y_1,\cdots,y_k])^{tC_p}$ is $y_1$-complete,
\end{enumerate}
after which the result follows from the fact the $y_1$-Bockstein spectral sequence is concentrated in even degrees.

To prove \eqref{it1-HT}, note that the $\mathbb{E}_1$-$\MU^{\otimes k}[y_2,\cdots,y_k]$-algebra structure on $\mathbb{F}_p$ admits a unique lift to an $\mathbb{E}_{\infty}$-$\MU^{\otimes k}[y_2,\cdots,y_k]$-algebra structure, since $\mathbb{F}_p$ is discrete and $\MU^{\otimes k}[y_2,\cdots,y_k]$ is connective.  In particular, there is a natural $\mathbb{E}_{\infty}$-ring map 
\[\THH(\mathbb{F}_p) \to \THH(\mathbb{F}_p/\MU^{\otimes k}[y_2,\cdots,y_k]),\]
and there is an isomorphism
\[\THH(\mathbb{F}_p/\MU^{\otimes k}[y_2,\cdots,y_k]) = \THH(\mathbb{F}_p) \otimes_{\THH(\MU^{\otimes k}[y_2,\cdots,y_k])} \MU^{\otimes k}[y_2,\cdots,y_k].\]
Recalling that
\[\pi_*\THH(\MU^{\otimes k}[y_2,\cdots,y_k])\]
is an exterior algebra on odd degree classes over $\pi_*\MU^{\otimes k}[y_2,\cdots,y_k]$, the K\"unneth spectral sequence associated to this relative tensor product proves that $\THH(\mathbb{F}_p/\MU^{\otimes k}[y_2,\cdots,y_k])$ is a free $\THH(\mathbb{F}_p)$-module, on even degree generators.  We may thus calculate that \[\THH(\mathbb{F}_p/\MU^{\otimes k}[y_2,\cdots,y_k])^{tC_p}\] is concentrated in even degrees, by considering the Tate spectral sequence as a module over the Tate spectral sequence for $\THH(\mathbb{F}_p)^{tC_p}$.

To prove \eqref{it2-HT}, observe that any bounded below $\MU^{\otimes k}[y_1,\cdots,y_k]$-module is automatically $y_1$ complete, since $y_1$ is in positive degree.  In particular, both \[\THH\left(k(n)/\MU^{\otimes k}[y_1,\cdots,y_k]\right)\] and \[\THH\left(k(n)/\MU^{\otimes k}[y_1,\cdots,y_k]\right)_{hC_p}\] are $y_1$-complete.  It remains to check that \[\THH\left(k(n)/\MU^{\otimes k}[y_1,\cdots,y_k]\right)^{hC_p}\] is $y_1$-complete, which follows from the fact that it is a limit of $y_1$-complete objects.
\end{proof}

\begin{proof}[Proof of Proposition~\ref{integral-Hodge-Tate-is-even}]
Since $\pi_*\grmot^*\THH(k(n);\bF_p)$ is concentrated on the zero and $1$-line and $v_n$ is in Adams weight zero, we know 
\[
\pi_*\grmot^*\THH(k(n))
\] 
is concentrated on the $0$ and $1$ line and therefore $\pi_*\grmot^*\TP(k(n))$ is concentrated on the $0$ and $1$ line and $\pi_*\grmot^*\THH(k(n))^{tC_p}$
 is concentrated on the $-1$, $0$, and $1$ lines. 

By Proposition~\ref{prop:change-of-eff-cover}, we can compute $\gr_{\mot}^*\THH(k(n))^{tC_p}$ using the cobar complex 
\[
\pi_*\THH(k(n)/\MU[y]^{\otimes \bullet+1})^{tC_p}
\]
and by Proposition~\ref{prop:evenness} we know $\THH(k(n)/\MU[y]^{\otimes \bullet+1})^{tC_p}
$
is even and therefore the motivic spectral sequence is concentrated in non-negative Adams weights, proving the claim.
\end{proof}

\subsubsection{Hodge--Tate cohomology of $k(n)$ with coefficients in $\bF_p$}

We next consider the diagram
\[
\begin{tikzcd}
\MU \arrow{r}{\Delta}
 \arrow{d} & (\MU^{\otimes p})^{tC_p} \arrow{d} \\
\THH(\MU) \arrow{r}{\varphi_p} \arrow{d} & \THH(\MU)^{tC_p} \arrow{d} \\
\THH(k(n)) \arrow{r}{\varphi_p} & \THH(k(n))^{tC_p},
\end{tikzcd}
\]
in which the top square is the natural one of $\mathbb{E}_{\infty}$-rings, the lower vertical maps are units of $\mathbb{E}_{0}$-algebra structures, and the lower square constitutes a map of $\mathbb{E}_{0}$-algebras.  This diagram equips 
\[\varphi_p:\THH(k(n)) \to \THH(k(n))^{tC_p}\] with the structure of an $\MU$-module map, which in particular lets us make sense of the expression
\[\THH(k(n)) / v_n \to \THH(k(n))^{tC_p} / 
\varphi_p(v_n) \,.
\]

For notational clarity, it can be convenient to rewrite this as 
\[
\THH(k(n);\mathbb{F}_p) \to \THH(k(n);\mathbb{F}_p^{\otimes p})^{tC_p} \,,
\]
where we recall that, if $R$ is an $\bE_1$-ring and $M$ is an $R$-bimodule, we may form
$\THH(R;M)$ and a Frobenius map
\[
\THH(R;M) \to \THH(R;M^{\otimes p})^{tC_p}\,.
\]
Here, $M^{\otimes p}$ refers to the $p$-th tensor power of $M$ in the category of $R$-bimodules. This construction is functorial in pairs $(R,M)$, and exact in the variable $M$ (cf.~\cite[Theorem~1.7]{NS18} and~\cite[Theorem~6.31]{KMN23}).

\begin{proposition}
On homotopy groups, the commutative diagram
\[
\begin{tikzcd}
 \THH(k(n);\mathbb{F}_p) \arrow{r}{\varphi_p} \arrow{d} &  \THH(k(n);\mathbb{F}_p^{\otimes p})^{tC_p}  \arrow{d} \\ 
\THH(\mathbb{F}_p) \arrow{r} & \THH(\mathbb{F}_p)^{tC_p}
\end{tikzcd}
\]
may be identified as $\Lambda(\lambda_{n+1})$-modules with the square
\[
\begin{tikzcd}
\mathbb{F}_p[\mu]/\mu^{p^{n}} \otimes \mathbb{F}_p[\mu^{p^{n+1}}] \otimes \Lambda(\lambda_{n+1}) \arrow{r} \arrow{d} & \mathbb{F}_p[\mu]/\mu^{p^{n}} \otimes \mathbb{F}_p[\mu^{\pm p^{n+1}}] \otimes \Lambda(\lambda_{n+1}) \arrow{d} \\
\mathbb{F}_p[\mu] \arrow{r} & \mathbb{F}_p[\mu^{\pm 1}].
\end{tikzcd}
\]
\end{proposition}

\begin{proof}
Filtering $\BP$ and $k(n)$ by the (d\'ecalage of) their $\mathbb{F}_p$-Adams towers following~\cite[Appendix~C.2]{HW22}, we arrive at a diagram of Hochschild--May spectral sequences, in the sense of \cite{AKS18,Kee20,JLL23}, with $\EE_1$-terms
\[
\begin{tikzcd}
\THH(\mathbb{F}_p[v_0,v_1,\cdots];\mathbb{F}_p) \arrow{d} \arrow{r} & \THH(\mathbb{F}_p[v_1,\cdots,v_n];\mathbb{F}_p^{\otimes p})^{tC_p} \arrow{d} \\
\pi_*\THH(\mathbb{F}_p[v_n];\mathbb{F}_p) \arrow{r} \arrow{d}& \pi_*\THH(\mathbb{F}_p[v_n];\mathbb{F}_p^{\otimes p})^{tC_p} \arrow{d} \\
\pi_*\THH(\mathbb{F}_p) \arrow{r} & \pi_*\THH(\mathbb{F}_p)^{tC_p}
\end{tikzcd}
\]
which is isomorphic to
\[
\begin{tikzcd}
\bF_p[\mu] \otimes \Lambda(dv_0,dv_1,\cdots) \arrow{r} \arrow{d} & \bF_p[\varphi_p(\mu)^{\pm 1}] \otimes  \Lambda(\varphi_p(dv_0),\varphi_p(dv_1),\cdots) \arrow{d} \\
\bF_p[\mu] \otimes \Lambda(dv_n) \arrow{d} \arrow{r} & \bF_p[\varphi_p(\mu)^{\pm 1}] \otimes \Lambda(\varphi_p(dv_n)) \arrow{d} \\
\bF_p[\mu] \arrow{r} & \bF_p[\varphi_p(\mu)^{\pm 1}]\,.
\end{tikzcd}
\]
The spectral sequences beginning with the middle horizontal arrow do not have multiplicative structure, but are modules over the spectral sequences beginning with the upper horizontal arrow.  In particular, the result follows from the description of the upper left spectral sequence in Proposition~\ref{prop:Hocschild-MaySSBP} as in the proof of Proposition~\ref{prop:Hochschild-MaySSkn}, where we abuse notation and write $\mu$ for $\varphi_p(\mu)$ and $\lambda_{n+1}$ for $\varphi_p(\lambda_{n+1})$. 
\end{proof}

We can identify
\[
\gr^*_{\mot} \THH(k(n);\mathbb{F}_p)
\]
with the quotient of $\gr^*_{\mot} \THH(k(n))$ by the element $v_n \in \pi_*\gr^*_{\mot} \THH(\MU)$.  We can similarly identify 
\[
\gr^*_{\mot} \THH(k(n);\mathbb{F}_p^{\otimes p})^{tC_p}
\]
with the quotient of $\gr^*_{\mot} \THH(k(n))^{tC_p}$ by $\varphi_p(v_n)$.  This allows us to make sense of the following:

\begin{corollary}
On bigraded homotopy groups, the commutative diagram
\[
\begin{tikzcd}
\gr^*_{\mot} \THH(k(n);\mathbb{F}_p) \arrow{r} \arrow{d} & \gr^*_{\mot} \THH(k(n);\mathbb{F}_p^{\otimes p})^{tC_p}  \arrow{d} \\ 
\gr^*_{\mot}\THH(\mathbb{F}_p) \arrow{r} & \gr^*_{\mot} \THH(\mathbb{F}_p)^{tC_p}
\end{tikzcd}
\]
may be identified as $\Lambda(\lambda_{n+1})$-modules with the square
\[
\begin{tikzcd}
\mathbb{F}_p[\mu]/\mu^{p^{n}} \otimes \mathbb{F}_p[\mu^{p^{n+1}}] \otimes \Lambda(\lambda_{n+1}) \arrow{r} \arrow{d} & \mathbb{F}_p[\mu]/\mu^{p^{n}} \otimes \mathbb{F}_p[\mu^{\pm p^{n+1}}] \otimes \Lambda(\lambda_{n+1}) \arrow{d} \\
\mathbb{F}_p[\mu] \arrow{r} & \mathbb{F}_p[\mu^{\pm 1}].
\end{tikzcd}
\]
\end{corollary}

\begin{proof}
Note that the motivic spectral sequence for $\gr^*_{\mot} \THH(k(n);\mathbb{F}_p^{\otimes p})^{tC_p}$ is a priori concentrated on the $0$, $1$, and $-1$ lines.  Thus, no motivic differentials are possible, and the motivic spectral sequence collapses on the $\EE_2$-page. To see this, note that the $1$-line and the $-1$-line are only nontrivial in degrees of the same parity and therefore the shortest possible differential in the motivic spectral sequence is a $d_3$. This is a feature of the motivic filtration, i.e. we are working in the even synthetic category, cf.~\cite[\S~5.2]{Pst23a} There are also a priori no filtration jumps possible in any of the maps.  Thus, it suffices to prove the result nonmotivically, where it follows from the previous proposition.
\end{proof}

\subsubsection{Hodge--Tate cohomology of $k(n)$ modulo $(p,v_1,\cdots,v_{n+1})$}

In this section, we finally deduce the main theorems about the mod $(p,v_1,\cdots,v_{n+1})$ Nygaard completed Hodge--Tate cohomology. 

\begin{proof}[Proof of Theorem~\ref{thm:Hodge-Tate}]
Consider the diagram
\[
\begin{tikzcd}
\MU \arrow{r}{\varphi_p} \arrow{d} & (\MU^{\otimes p})^{tC_p} \arrow{d} \\
\THH(\MU;\mathbb{F}_p) \arrow{r}{\varphi_p} \arrow{d} & \THH(\MU;\mathbb{F}_p^{\otimes p})^{tC_p} \arrow{d} \\
\THH(k(n);\mathbb{F}_p) \arrow{r}{\varphi_p} & \THH(k(n);\mathbb{F}_p^{\otimes p})^{tC_p},
\end{tikzcd}
\]
where the upper square is one of $\mathbb{E}_{\infty}$-algebras, the lower vertical maps are units of $\mathbb{E}_0$-algebra structures, and the lower square constitutes a map of $\mathbb{E}_0$-algebras.
The diagram equips 
\[\varphi_p:\THH(k(n);\mathbb{F}_p)\cong \THH(k(n))/v_n \to \THH(k(n))^{tC_p}/\varphi_p(v_n) \cong \THH(k(n);\mathbb{F}_p^{\otimes p})^{tC_p}\]
with the structure of a $\THH(\MU;\mathbb{F}_p)$-module map.  In particular, since $\THH(\MU;\mathbb{F}_p)$ is an $\mathbb{F}_p$-algebra this gives nullhomotopies of the actions of $p,v_1,\cdots,v_{n-1}$ and $v_{n+1}$. Similarly, applying $\gr^*_{\mot}$ and using the fact that $p,v_1,\cdots,v_{n+1}$ are all zero in $\pi_*\gr^*_{\mot} \THH(\MU;\mathbb{F}_p)$, we obtain for each $i$ combatible nullhomotopies of $v_i$ acting on $\gr^*_{\mot} \THH(k(n);\mathbb{F}_p)$ and $\varphi_p(v_i)$ acting on $\gr^*_{\mot} \THH(k(n)) / \varphi(v_n)$.

Using our preferred basis for $\pi_*\grmot^*\THH(k(n))/(p,v_1,\cdots,v_{n+1})$, we conclude from this and the previous subsection that the commutative diagram 
\[
\begin{tikzcd}
    K_n  \ar[r]\ar[d] & K_n^t \ar[d]\\ 
    K_{\infty}  \ar[r] & K_{\infty}^t 
\end{tikzcd} 
\]
from Definition~\ref{Kn-Knt}, 
may be identified with the commutative diagram of $\Lambda(\lambda_{n+1})$-modules 
\[
\begin{tikzcd}
   ( L \oplus M)  \otimes \Lambda( \lambda_{n+1})\otimes \bF_p[\mu^{p^{n+1}}] \ar[r] \ar[d] & (
  L \oplus M ) \otimes \Lambda( \lambda_{n+1})\otimes \bF_p[\mu^{\pm p^{n+1}}] \ar[d] \\
L \otimes \Lambda (\bar{\varepsilon}_{n+1})\otimes \bF_p[\mu^{p^n+1}] \ar[r] &  L\otimes \Lambda(\overline{\varepsilon}_{n+1}) \otimes \bF_p[\mu^{\pm p^{n+1}}] 
\end{tikzcd}
\]
where $L:=\Lambda (\overline{\varepsilon}_1,\cdots ,\overline{\varepsilon}_n)$, $M:=\bF_p\{x\mu ,xz_{n+1}: x\in B \}$, and $B$ is defined as in Definition~\ref{def:B}. Here the vertical maps are given by the canonical quotient by $(\lambda_{n+1})$ composed with the canonical inclusion, the bottom horizontal map is the localization map and the top map is the restriction of the localization map.  
\end{proof}

\begin{proof}[Proof of Theorem \ref{thm:Hodge-Tateresults}]
We consider the spectral sequence beginning with
\[\pi_*\grmot^*\THH(k(n))^{tC_p}/(p,\cdots ,v_{n+1})[t] \,.\]
We wish to check that the $\mathrm{E}_2$-page is concentrated on the $0$-line (which in particular implies that the spectral sequence collapses after the $d_1$ differential).

The entire spectral sequence is a module over $\Lambda(\lambda_{n+1})$, with $\lambda_{n+1}$ a permanent cycle.  In particular, the $d_1$ differential is compatible with the (2-stage) $\lambda_{n+1}$-adic filtration on the $\mathrm{E}_1$-page.  At the level of associated graded of the $\lambda_{n+1}$-adic filtration, the $d_1$ differential can be calculated using the injection 
\[\pi_*\grmot^*\THH(k(n))^{tC_p}/(p,\cdots ,v_{n+1})[t]\subset \bF_p[\mu^{\pm 1}]\otimes \Lambda (\varepsilon_0,\overline{\varepsilon}_1,\cdots,\overline{\varepsilon}_{n+1}) \,.\]
In the target spectral sequence the differentials are generated by the differential  
$d_1(\epsilon_{0}) = t\mu$ using the Leibniz rule and the fact that all other classes are $d_1$-cycles. Since the source is closed under the $\sigma$ operator, the map remains an injection on $\mathrm{E}_2$-pages. This implies that the spectral sequence collapses at the $\EE_2$-term and is concentrated on the $0$-line. To determine the $0$-line, we use our preferred basis and consider the inclusion of an $\bF_p[\mu^{\pm p^{n+1}},t]\otimes \Lambda (\varepsilon_{n+1})$-submodule 
\[
(\Lambda(\overline{\varepsilon}_1,\cdots ,\overline{\varepsilon}_n)\oplus \bF_2\{ x\mu,x\varepsilon_0 : x\in B\})\otimes \Lambda( \epsilon_{n+1}) \otimes \Lambda(\lambda_{n+1})\otimes \bF_2[\mu^{p^{n+1}}] \otimes [t] 
\]
in 
\[
\mathbb{F}_p[\mu] \otimes \Lambda(\epsilon_0,\overline{\varepsilon}_1,\cdots,\overline{\epsilon}_{n+1})[t]
\]
determined by Theorem~\ref{thm:height-n}. We can alternatively describe the $0$-line as $K_n^t$ from Definition~\ref{Kn-Knt}, which was computed in Theorem~\ref{thm:Hodge-Tate}. 
\end{proof}

\subsection{Prismatic cohomology of $k(n)$}\label{sec:prismatickn}

In the previous subsection, we used the identity
\[(\THH(k(n))^{tC_p})^{hS^1} \cong \TP(k(n))\]
to name every class in
\[\pi_*\gr^*_{\mot} \TP(k(n))/(p,v_1,\cdots,v_{n+1}).\]
In this subsection, we shall first determine exactly where those classes are detected in the periodic $t$-Bockstein spectral sequence 
\[\mathrm{E}_1=\left(\pi_*\gr^*_{\mot} \THH(k(n)) / (p,v_1,\cdots,v_{n+1})\right)[t^{\pm 1}] \implies \pi_* \gr^*_{\mot} \TP(k(n)) / (p,v_1,\cdots,v_{n+1}).\]

Once we know which classes are permanent in the periodic $t$-Bockstein spectral sequence, we will be able to compute its exact pattern of differentials.  Knowing this pattern of differentials will tell us not only $\pi_*\gr^*_{\mot} \TP(k(n)) / (p,v_1,\cdots,v_{n+1})$, which was already computed in the previous subsection, but also $\pi_* \gr^*_{\mot} \TC^{-}(k(n)) / (p,v_1,\cdots,v_{n+1})$.

To begin, we note a general fact about differentials in the (periodic) $t$-Bockstein spectral sequence. 

\begin{lem}\label{lem:lambda-leibniz}
The differentials in the (periodic) $t$-Bockstein spectral sequence satisfy 
\[ d_r(\lambda_{n+1} x)=\lambda_{n+1}d_r(x) 
\]
for all $r\ge 1$. 
\end{lem}
\begin{proof}
This follows from the module structure over the (periodic) $t$-Bockstein spectral sequence 
\[
\pi_*\grmot^*\THH(\BP)/(p,v_1,\cdots,v_{n+1})[t^{\pm 1}]\implies  \pi_*\grmot^*\TP(\BP)/(p,v_1,\cdots,v_{n+1})
\]
in which $\lambda_{n+1}$ is known to be a permanent cycle by Corollary~\ref{cor:lambdas-perm}. 
\end{proof}
We now determine the $\EE_{2}$-page of the periodic $t$-Bockstein spectral sequence:  

\begin{prop}\label{prop: E4 page}
The $\EE_2$-page of the periodic $t$-Bockstein spectral sequence \begin{equation}\label{eq: periodic t-bss}
\pi_*\grmot^*\THH(k(n))/(p,\cdots ,v_{n+1})[t^{\pm 1}]\implies \pi_*\grmot^*\TP(k(n))/(p,\cdots ,v_{n+1})\,.
\end{equation}
can be identified with
\[ 
\Lambda(\bar{\varepsilon}_1,	\dots, \bar{\varepsilon}_{n},\lambda_{n+1})[t^{\pm 1}]\,.
\] 
\end{prop}

\begin{proof}
The map of periodic $t$-Bockstein spectral sequences induced by the map $k(n)\to \bF_p$ is an injection on $\EE_1$-pages modulo $\lambda_{n+1}$. 
The periodic $t$-Bockstein spectral sequence for $\bF_p$ has differentials $d_1(\varepsilon_i)=t\mu^{p^i}$ by Proposition~\ref{prop:sigmaFp}. The result then follows from the description of $\pi_{*}\gr_{\mot}^*\THH(k(n))/(p,v_1,\cdots ,v_{n+1})$ in Theorem~\ref{thm:height-n} along with the Leibniz rule for multiplication by $\lambda_{n+1}$ from Lemma~\ref{lem:lambda-leibniz}. 
\end{proof}
 
 

\begin{rem2}
By Theorem~\ref{thm:Hodge-Tateresults}, we have an isomorphism
\[ \pi_*\grmot^*\TP(k(n))/(p,v_1,\cdots,v_{n+1}) \cong \pi_*\grmot^*\THH(k(n))^{tC_p}/(p,v_1,\cdots,v_n) \,.
\]
where both sides of this isomorphism compute the mod $(p,v_1,\cdots,v_{n+1})$ Nygaard completed prismatic cohomology of $k(n)$.  Theorem~\ref{thm:Hodge-Tateresults} gives an explicit basis for this $\mathbb{F}_p$ vector space, but in the proposition below we will give a different basis that identifies Nygaard filtrations on elements.
\end{rem2}

\begin{rem2} \label{rem:FpGinfinity}
The spectral sequence beginning with 
\[ \pi_* \gr^*_{\mot} \THH(\mathbb{F}_p)^{tC_p} / (p,v_1,\cdots,v_{n+1}) [t]\]
and converging to $\pi_*\gr^*_{\mot} \mathrm{TP}(\mathbb{F}_p) / (p,v_1,\cdots,v_{n+1})$ collapses at the $\mathrm{E}_2$=$\mathrm{E}_{\infty}$-page, with all permanent cycles on the $0$-line.  This $0$-line is isomorphic to
\[\Lambda(\overline{\varepsilon}_1,\cdots,\overline{\varepsilon}_n) \otimes \mathbb{F}_p[\mu^{\pm 1}].\]

Since the $\mathrm{E}_2=\mathrm{E}_{\infty}$-page is entirely concentrated on the $0$-line, any class $\overline{\varepsilon}_{S} \mu^j$ on the $0$-line refers to a uniquely specified element of $\mathrm{TP}(\mathbb{F}_p) / (p,v_1,\cdots,v_{n+1})$.  This class is detected in the periodic $t$-Bockstein spectral sequence 
\[\mathrm{E}_1=\Lambda(\epsilon_0,\overline{\varepsilon}_1,\cdots,\overline{\varepsilon}_n) [t^{\pm 1}] \implies \pi_*\gr^* \TP(\mathbb{F}_p) / (p,v_1,\cdots,v_{n+1})\]
by $\overline{\varepsilon}_{S} t^{-j}$.
\end{rem2}


\begin{prop}\label{prop:perm-cycles-tp}
The mod $(p,v_1,\cdots ,v_{n+1})$ Nygaard completed prismatic cohomology of $k(n)$ can be identified, as a $\Lambda(\lambda_{n+1})$-module, with 
\[
\left ( \Lambda (\overline{\varepsilon}_1,\cdots ,\overline{\varepsilon}_n) \oplus \mathbb{F}_p\{xt^{-1},y: x\in N \,,y\in P\} \right ) \otimes \Lambda (\lambda_{n+1})\otimes \bF_p[t^{ \pm p^{n+1}}].
\] 
where 
\begin{align*}
N:= \coprod_{\overset{S\subset \{1,\cdots ,n\}}{ |S|<n}} \{ \overline{\varepsilon}_St^{-j}  &:0\le j \le p^{n}-1+f(S)\}  
\end{align*}
and 
\begin{align*}
P=\coprod_{S\subset \{1,\cdots ,n\}}\{ \overline{\varepsilon}_St^j &: 0<j<-f(S)\} 
\end{align*}
Here, elements are named according to where they are detected on the $\mathrm{E}_{\infty}$-page of the periodic $t$-Bockstein spectral sequence, so that $\overline{\varepsilon}_i\in \{\overline{\varepsilon}_i\}$, $\overline{\varepsilon}_{S}t^{j} \in \{\overline{\varepsilon}_St^{j}\}$, $t^{p^{n+1}}\in \{t^{p^{n+1}}\}$ and $\lambda_{n+1}$ is the unique class in $\{\lambda_{n+1}\}$ in the periodic $t$-Bockstein spectral sequence. The quantity $f(S)$ is specified according to Definition~\ref{e(S)-f(S)}. 
\end{prop}

\begin{proof}
As we will spell out in more detail below, the entire situation is determined by Theorem~\ref{thm:Hodge-Tateresults} and the reduction map between periodic $t$-Bockstein spectral sequences converging to $\pi_*\grmot^*\TP(k(n)) / (p,v_1,\cdots,v_{n+1})$ and $\pi_*\grmot^* \TP(\mathbb{F}_p) / (p,v_1,\cdots,v_{n+1})$, respectively.  As in Remark \ref{rem:FpGinfinity}, we understand very well the codomain of this map, and this turns out to be enough to prove the proposition.

In more detail, according to Theorem~\ref{thm:Hodge-Tateresults}, every class in $\pi_*\grmot^* \TP(k(n)) / (p,v_1,\cdots,v_{n+1})$ is uniquely specified by a class in
\[K_n^t= \left (
   \Lambda (\overline{\varepsilon}_1,\cdots ,\overline{\varepsilon}_n)\oplus \bF_p\{x\mu ,xz_{n+1}: x\in B \} \right ) \otimes \Lambda( \lambda_{n+1})\otimes \bF_p[\mu^{\pm p^{n+1}}],\]
where
\begin{align*}
B:= \coprod_{\overset{S\subset \{1,\cdots ,n\}}{|S|<n}} \{ \bar{\varepsilon}_S\mu^{j}  &:0\le j < p^{n}-1+f(S)\} 
\end{align*}
This naming scheme is in particular designed to specify where classes map to in the homotopy fixed point spectral sequence for $(\THH(\mathbb{F}_p)^{tC_p})^{hS^1}$.

Classes in $K_n^{t}$ of the form $x \mu$, where $x=\overline{\epsilon}_{S} \mu^j$, map to elements $\overline{\epsilon}_{S}\mu^{j+1}$ that are detected by $\overline{\epsilon}_{S} t^{-j-1}$.   It follows that the Nygaard filtration (i.e., periodic $t$-Bockstein filtration) of $x\mu$ in $\pi_* \gr^*_{\mot} \TP(k(n)) / (p,v_1,\cdots,v_{n+1})$ is \emph{at most} $-j-1$.  However, we might worry that $x\mu$ is detected in the periodic $t$-Bockstein spectral sequence by some class that maps to zero.  In fact, this is not possible for degree reasons.  The only classes mapping to zero are multiples of $\lambda_{n+1}$, which are in higher Nygaard filtration.

The situation is similar for classes of the form $xz_{n+1}$ where $x=\overline{\varepsilon}_S\mu^j$, but it is slightly more difficult to understand where these are detected after reduction.

Specifically, one must take care to note that each summand in the expression
\[ z_{n+1} :=\overline{\varepsilon}_{n+1}\mu+\overline{\varepsilon}_{n}\mu^{p^{n+1}-p^{n}+1}+\cdots +\overline{\varepsilon}_1\mu^{p^{n+1}-p+1}  
\]
is detected in different Nygaard filtration in 
\[
\pi_*\grmot^*\TP(\bF_p)/(p,v_1,\cdots ,v_n)\,.
\]  
In $\pi_*\grmot^*\TP(\bF_p)/(p,v_1,\cdots ,v_n)$, we know that
\[ [[xz_{n+1}]]\in \{x\overline{\varepsilon}_st^{-p^{n+1}+p^s-1}\} 
\]
where $s\in \{1,2,\cdots ,n+1\}$ is the smallest integer such that $x\overline{\varepsilon}_s\mu^{p^{n+1}-p^s+1}\ne 0\in K_{\infty}^t$ for $x\in B$. Note that $x\overline{\varepsilon}_{n+1}\ne 0$ for all $x\in B$. This is sufficient when $xz_{n+1}$ is a monomial in $K_n^t$. 
When $xz_{n+1}$ is not a monomial in $K_n^t$, then we follow the same algorithm for each term in the sum. 

For example, when $n=2$ the class $\overline{\varepsilon}_1+t^{p^2-p}\overline{\varepsilon}_2+\overline{\varepsilon}_3t^{p^3-p}=t^{p^{3}-p+1}z_{3}$ is the sum of $\overline{\varepsilon}_1$ and $t^{p^2-p}\overline{\varepsilon}_2+\overline{\varepsilon}_3t^{p^3-p}$ in $\grmot^*\TP(k(n))/(p,v_1,\cdots ,v_{n+1})$. So we first choose the class $\overline{\varepsilon}_1\in \{\varepsilon_1\}$ and then we choose a class 
\[
t^{p^2-p}\overline{\varepsilon}_2+t^{p^3-p}\overline{\varepsilon}_3\in \{t^{p^2-p}\overline{\varepsilon}_2\}\,.
\]
Following this algorithm produces the description in the statement of the proposition for classes in the image of the reduction map. 

Note that we can always choose  $[[\lambda_{n+1}xt^{-j}]]\in \{\lambda_{n+1}xt^{-j}\}$ as the product $[[\lambda_{n+1}]]\cdot [[xt^{-j}]]$ where $[[\lambda_{n+1}]]\in \{\lambda_{n+1}\}$ and $[[xt^{-j}]]\in \{xt^{-j}\}$. To see this, note that by 
inspection there is a unique class $\lambda_{n+1}$ in $\{\lambda_{n+1}\}$, for example see Figure~\ref{fig:prismatic}. Moreover, the indeterminacy of $[[\lambda_{n+1}xt^{-j}]]$ is exactly $\lambda_{n+1}$ multiplied by the indeterminacy of $[[xt^{-j}]]$. 
\end{proof}

The following result will be important for identifying $\pi_*\grmot^*\TC^{-}(k(n))/(p,v_1,\cdots,v_{n+1})$
\begin{thm}\label{thm:prismatic}
There is a unique possible pattern of differentials in the periodic $t$-Bockstein spectral sequence computing the mod $(p,v_1,\cdots,v_{n+1})$ Nygaard completed prismatic cohomology of $k(n)$, starting with $\mathrm{E}_2$-page 
\[\Lambda(\overline{\varepsilon}_1,\cdots,\overline{\epsilon}_n,\lambda_{n+1})[t^{\pm 1}]
\]
and ending with $\EE_{\infty}$-page 
\[
\left ( \bigoplus_{S\subset \{1,\cdots ,n \}}\bF_p\{\overline{\varepsilon}_St^{j} : -p^n-f(S)<j \le -f(S) \} \right )\otimes \Lambda (\lambda_{n+1})\otimes \bF_p[t^{\pm p^{n+1}}] \,.
\] 
The quantity $f(S)$ is specified according to Definition~\ref{e(S)-f(S)}. The differentials are all of length $p^{n+1}$, and given by the formula
\begin{align*}
    d_{p^{n+1}}(t^{sp^{n+1}+j}\bar{\varepsilon}_{S})&=t^{(s+1)p^{n+1}+j}\lambda_{n+1}\bar{\varepsilon}_{S}  
\end{align*}
where $s$ is an integer and $-f(S)<j \le p^{n+1}-p^n-f(S)$. 
\end{thm}

\begin{proof}
We claim that there is only one pattern of differentials that is consistent with each of the permanent cycles described in Proposition~\ref{prop:perm-cycles-tp}
being permanent cycles. 

By $t^{p^{n+1}}$-linearity, it suffices to determine the differentials in the range of stems $0\le s \le 2p^{n+1}$. Here $t^{p^{n+1}}$-linearity of the differentials follows from the module structure of the periodic $t$-Bockstein spectral sequence for $k(n)$ over the periodic $t$-Bockstein spectral sequence for $\BP$ and the fact that all differentials land in zero groups by the $\EE_{p^{n+2}}$-page. 

To start, we give an argument for the differentials on powers of $t$.
We know that 
\[
d_r(t^{-p^n})=0
\] 
for $1\le r<p^{n+1}$ and 
\[
d_{p^{n+1}}(t^{-p^n})\dot{=}t^{p^{n+1}-p^n}\lambda_{n+1}
\] 
using the module structure of the periodic $t$-Bockstein spectral sequence for $k(n)$ over the periodic $t$-Bockstein spectral sequence for $\BP$. In particular, there is a Leibniz rule for multiplication by $t^{p^n}$ up until (and including) the $\EE_{p^{n+1}}$-page. This implies that $t^{-\ell}$ is a permanent cycle for all $\ell\equiv 0,\cdots,p^{n}-1\mod p^{n+1}$ and consequently there are differentials 
\[
d_r(t^{-\ell})\dot{=}t^{p^{n+1}-\ell}\lambda_{n+1}
\]
whenever $\ell\not\equiv 0,\cdots,p^{n}-1\mod p^{n+1}$. Or in other words, there are differentials 
\[
d_r(t^{sp^{n+1}+j})\dot{=}t^{(s+1)p^{n+1}+j}\lambda_{n+1}
\]
whenever $0<j\le p^{n+1}-p^n$ and $s$ is an integer.  

We now provide an argument for differentials whose source is in stem $1$, for exposition. Let 
\[
    w(S):=\sum_{s\in S}p^s\,.
\]
We note that each other class in stem $1$, that is not a multiple of $\lambda_{n+1}$, is the source of a differential besides $t^{p-1}\bar{\varepsilon}_1$. Note that $\lambda_{n+1}x$ is a permanent cycle if and only if $x$ is a permanent cycle by Lemma~\ref{lem:lambda-leibniz}, so a similar argument for each stem $0\le s <2p^{n+1}$ will suffice by $t^{p^{n+1}}$-linearity. 

In the case $n=1$, this is clear, so it suffices to consider the case $n>1$. 

Note that the classes besides $t^{p-1}\bar{\varepsilon}_1$ and multiples of $\lambda_{n+1}$ are all of the form 
\[
    t^{w(S)-(1+|S|)/2}\bar{\varepsilon}_S
\] 
where $|S|\equiv 1\mod 2$ and $S\ne \{1\}$ since 
\[
|\overline{\varepsilon}_S|-2w(S)+1+|S|=1\,.
\]

We know that for each such $S$ the classes 
\[
t^{-p^n+w(S)-(1+|S|)/2}\bar{\varepsilon}_S
\] 
are permanent cycles since
\[
 -\sum_{s\in S}p^{s-1}\le -p <-1\le -(1+|S|)/2
\]
when $S\ne \{1\}$ and consequently
\begin{align*}
-p^n-f(S)&=-p^n+\sum_{s\in S}p^s-p^{s-1} \\
&< -p^n+\sum_{s\in S}p^s-(1+|S|)/2 \\
&=-p^n+w(S)-(1+|S|)/2  \\
&\le -f(S) 
\end{align*}
where the last inequality holds because
\[
-p^n-(1+|S|)/2 \le -\sum_{s\in S}p^{s-1}
\] 
for $|S|$ odd. Suppose that 
\begin{equation}\label{eq:non-existent-diff}
d_r(
t^{w(S)-(1+|S|)/2-j+1}\bar{\varepsilon}_S)=y\ne 0
\end{equation}
for each $|S|\equiv 1\mod 2$ with $S\ne \{1\}$ and some $r<p^{n+1}$. Then for bidegree reasons, 
\[
 y\in \{ t^{p^{n+1}+w(S)-(1+|S|)/2}\bar{\varepsilon}_S\lambda_{n+1} : |S|\equiv 1\mod 2\} 
\]
up to multiplication by a unit. 

Since 
\[
d_r(t^{-p^n})=0
\]
for $r<p^{n+1}$ and $t^{-p^n}$ acts on the $\EE_r$-page of spectral sequence for $r<p^{n+1}$, the existence of a differential \eqref{eq:non-existent-diff} implies a differential of the form 
\begin{equation}\label{eq:non-existent-diff-2}
d_r(t^{-p^n+w(S)-(1+|S|)/2}\bar{\varepsilon}_S)=t^{-p^n}y\ne 0 
\end{equation}
leading to a contradiction. 

Note that the finite set 
$\{\bar{\varepsilon}_S: |S|\equiv 1\mod 2\}$
is totally ordered by the degree of $\bar{\varepsilon}_S$. Let $\bar{\varepsilon}_T$ be the class of maximal degree in the finite set $\{\bar{\varepsilon}_S: |S|\equiv 1\mod 2\}$ equipped with this total ordering. Since 
\[
t^{w(T)-(1+|T|)/2}\bar{\varepsilon}_T
\] 
is not a permanent cycle and we know it cannot be a boundary, by considering the map of spectral sequences induced by $k(n)\to \bF_p$, the only remaining possibility is a differential 
\[
d_{p^{n+1}}(t^{w(S)-1}\bar{\varepsilon}_S)=t^{p^{n+1}+w(S)-j}\bar{\varepsilon}_S\lambda_{n+1}\] 
up to multiplication by a unit. 
By a downward induction on the total order of the set $\{\bar{\varepsilon}_S : |S|\equiv 1\mod 2\,, S\ne \{1\}\}$, we determine that  
\[
d_{p^{n+1}}(t^{p^n+w(S)-1}\bar{\varepsilon}_S)\dot{=}t^{p^{n+1}+w(S)-j}\bar{\varepsilon}_S\lambda_{n+1}
\] 
for each $|S|\equiv 1\mod 2$ and $S\ne \{1\}$. This proves the claim.

We now provide a more general argument for even stems and odd stems in two separate cases. The idea of the argument is essentially the same as the one above. 

In stems $0\le 2j<2p^{n+1}$, the permanent cycles that are not $\lambda_{n+1}$-multiples are 
\[ 
\{t^{-j}\} 
\]
for $0\le j\le p-1$
\[ 
\{
t^{-j} ,t^{w(S)-|S|/2-j}\varepsilon_S  : |S|\equiv 0 \mod 2, w(S)-|S|/2-j\le-f(S)
\} 
\]
for $p\le j <p^n$
and 
\[ 
\{
t^{w(S)-|S|/2-j}\varepsilon_S  : |S|\equiv 0 \mod 2, w(S)-|S|/2-j\le-f(S)
\} 
\]
for $p^n<j<w(S)-\lfloor n/2 \rfloor$ by Theorem~\ref{prop:perm-cycles-tp}. 

The remaining classes in stem $0\le 2j<2p^{n+1}$ that are not $\lambda_{n+1}$-multiples are therefore 
\[
\{t^{w(S)-|S|/2-j}\varepsilon_S  : |S|\equiv 0 \mod 2, w(S)-|S|/2-j>-f(S) \} 
\]
for $0\le j <p^n$
and 
\[
\{t^j , t^{w(S)-|S|/2-j}\varepsilon_S  : |S|\equiv 0 \mod 2, w(S)-|S|/2-j>-f(S) \} 
\]
for $p^n \le j <p^{n+1}$.

Multiplying by $t^{-p^n}$ produces permanent cycles
\[
\{t^{-p^n+w(S)-|S|/2-j}\varepsilon_S  : |S|\equiv 0 \mod 2, w(S)-|S|/2-j>-f(S) \} 
\]
in degree $2j+p^n$ since 
\[ -p^n-f(S) < -p^n+w(S)-|S|/2-j \le -f(S) \]
whenever 
\[ -f(S) < w(S)-|S|/2-j \le -f(S)+p^n 
\,.
\]
Since $d_{p^{n+1}}(t^{p^n})=t^{p^{n+1}+p^n}\lambda_{n+1}$ the Leibniz rule implies a differential 
\[ d_{p^{n+1}}(t^{w(S)-|S|/2-j}\varepsilon_S)=t^{p^{n+1}+w(S)-|S|/2-j}\varepsilon_S
\]
provided that the class $t^{w(S)-|S|/2-j}\varepsilon_S$ survives until the $\EE_{p^{n+1}}$-page. We then observe that for a fixed integer $j$, the classes 
\[ \{t^{w(S)-|S|/2-j}\varepsilon_S : |S|\equiv 0 \mod 2
\}\]
are totally ordered by the degree of $\varepsilon_S$. Starting with largest class using this total ordering, we can use $t^{p^n}$-linearity, which holds up until the $\EE_{p^{n+1}}$-page, to determine that a differential on this largest class would contradict the fact that the $t^{-p^n}$-multiple of this class is a permanent cycle. Consequently, we must have the specified length $p^{n+1}$-differential. A finite downward induction using this total order then implies the stated differentials.

In stems $1< 2j-1<2p^{n+1}$, the permanent cycles that are not $\lambda_{n+1}$-multiples are 
\[ 
\{
t^{w(S)-(1+|S|)/2-j}\varepsilon_S  : |S|\equiv 1 \mod 2, w(S)-(1+|S|)/2-j\le-f(S)
\} 
\]
for $1\le j<w(S)-\lfloor n/2 \rfloor$ by Theorem~\ref{prop:perm-cycles-tp}. 

The remaining classes in stems $1<2j-1<2p^{n+1}$ that are not $\lambda_{n+1}$-multiples are therefore 
\[
\{t^{w(S)-(1+|S|)/2-j}\varepsilon_S  : |S|\equiv 1 \mod 2, w(S)-(1+|S|)/2-j>-f(S) \} \,.
\]

Multiplying by $t^{-p^n}$ produces permanent cycles
\[
\{t^{-p^n+w(S)-(1+|S|)/2-j}\varepsilon_S  : |S|\equiv 1 \mod 2, w(S)-|S|/2-j>-f(S) \} 
\]
in degree $2j-1+p^n$ since 
\[ -p^n-f(S) < -p^n+w(S)-(1+|S|)/2-j \le -f(S) \]
whenever 
\[ -f(S) < w(S)-(1+|S|)/2-j \le -f(S)+p^n 
,.
\]
Since $d_{p^{n+1}}(t^{p^n})=t^{p^{n+1}+p^n}\lambda_{n+1}$ the Leibniz rule implies a differential 
\[ 
d_{p^{n+1}}(t^{w(S)-(1+|S|)/2-j}\varepsilon_S)=t^{p^{n+1}+w(S)-(1+|S|)/2-j}\varepsilon_S
\]
provided that the class $t^{w(S)-(1+|S|)/2-j}\varepsilon_S$ survives until the $\EE_{p^{n+1}}$-page. 

We then observe that for a fixed integer $j$, the classes 
\[ \{t^{w(S)-(1+|S|)/2-j}\varepsilon_S : |S|\equiv 1 \mod 2\}\]
are totally ordered by $|\varepsilon_S|$. Starting with largest class using this total ordering, we determine that a shorter differential on this class would lead to a contradiction by $t^{-p^n}$-linearity at pages up until $\EE_{p^{n+1}}$ and the fact that the $t^{-p^{n}}$-multiple of this class is a permanent cycle. There are therefore no other possible differentials and therefore we must have the specified length $p^{n+1}$-differential. A finite downward induction using this total order then implies the stated differentials. 

\end{proof}

\begin{figure}[ht!]
\resizebox{\textwidth}{!}{ 
\begin{tikzpicture}[radius=1,yscale=1.5]
\foreach \n in {-2,-1,...,27} \node [below] at (\n,-.8-12) {$\n$};
\foreach \s in {-12,-11,...,5} \node [left] at (-.3-2,\s) {$\s$};
\draw [thin,color=lightgray] (-2,-12) grid (27,5);
\node [below] at (0,0) {$1$};
\node [below] at (2,-1) {$t^{-1}$};
\node [below] at (4,-2) {$t^{-2}$};
\node [below] at (6,-3) {$t^{-3}$};
\node [below] at (1,1) {$t\bar{\varepsilon}_1$};
\node [below] at (3,0) {$\bar{\varepsilon}_1$};
\node [below] at (5,-1) {$t^{-1}\bar{\varepsilon}_1$};
\node [below] at (7,-2) {$t^{-2}\bar{\varepsilon}_1$};
\node [below] at (3,2) {$t^{2}\bar{\varepsilon}_2$};
\node [below] at (5,1) {$t\bar{\varepsilon}_2$};
\node [below] at (7,0) {$\bar{\varepsilon}_2$};
\node [below] at (9,-1) {$t^{-1}\bar{\varepsilon}_2$};
\node [below] at (4,3) {$t^{3}\bar{\varepsilon}_1\bar{\varepsilon}_2$};
\node [below] at (6,2) {$t^{2}\bar{\varepsilon}_1\bar{\varepsilon}_2$};
\node [below] at (8,1) {$t\bar{\varepsilon}_1\bar{\varepsilon}_2$};
\node [below] at (10,0) {$\bar{\varepsilon}_1\bar{\varepsilon}_2$};

\node [below] at (15,0) {$\lambda_3$};
\node [below] at (17,-1) {$t^{-1}\lambda_{3}$};
\node [below] at (19,-2) {$t^{-2}\lambda_{3}$};
\node [below] at (21,-3) {$t^{-3}\lambda_{3}$};
\node [below] at (16,1) {$t\bar{\varepsilon}_1\lambda_{3}$};
\node [below] at (18,0) {$\bar{\varepsilon}_1\lambda_{3}$};
\node [below] at (20,-1) {$t^{-1}\bar{\varepsilon}_1\lambda_{3}$};
\node [below] at (22,-2) {$t^{-2}\bar{\varepsilon}_1\lambda_{3}$};
\node [below] at (18,2) {$t^{2}\bar{\varepsilon}_2\lambda_{3}$};
\node [below] at (20,1) {$t\bar{\varepsilon}_2\lambda_{3}$};
\node [below] at (22,0) {$\bar{\varepsilon}_2\lambda_{3}$};
\node [below] at (24,-1) {$t^{-1}\bar{\varepsilon}_2\lambda_{3}$};
\node [below] at (19,3) {$t^{3}\bar{\varepsilon}_1\bar{\varepsilon}_2\lambda_{3}$};
\node [below] at (21,2) {$t^{2}\bar{\varepsilon}_1\bar{\varepsilon}_2\lambda_{3}$};
\node [below] at (23,1) {$t\bar{\varepsilon}_1\bar{\varepsilon}_2\lambda_{3}$};
\node [below] at (25,0) {$\bar{\varepsilon}_1\bar{\varepsilon}_2\lambda_{3}$};

\node [below] at (16,-8) {$t^{-8}$};
\node [below] at (18,-9) {$t^{-9}$};
\node [below] at (20,-10) {$t^{-10}$};
\node [below] at (22,-11) {$t^{-11}$};
\node [below] at (17,-7) {$t^{-7}\bar{\varepsilon}_1$};
\node [below] at (19,-8) {$t^{-8}\bar{\varepsilon}_1$};
\node [below] at (21,-9) {$t^{-9}\bar{\varepsilon}_1$};
\node [above] at (23,-10) {$t^{-10}\bar{\varepsilon}_1$};
\node [below] at (19,-6) {$t^{-6}\bar{\varepsilon}_2$};
\node [below] at (21,-7) {$t^{-7}\bar{\varepsilon}_2$};
\node [below] at (23,-8) {$t^{-8}\bar{\varepsilon}_2$};
\node [below] at (25,-9) {$t^{-9}\bar{\varepsilon}_2$};
\node [below] at (20,-5) {$t^{-5}\bar{\varepsilon}_1\bar{\varepsilon}_2$};
\node [below] at (22,-6) {$t^{-6}\bar{\varepsilon}_1\bar{\varepsilon}_2$};
\node [below] at (24,-7) {$t^{-7}\bar{\varepsilon}_1\bar{\varepsilon}_2$};
\node [below] at (26,-8) {$t^{-8}\bar{\varepsilon}_1\bar{\varepsilon}_2$};
\end{tikzpicture}
}
\caption{The $\EE_\infty$-page of the periodic t-Bockstein spectral sequence computing the 
mod $(2,v_1,v_2,v_3)$-prismatic cohomology of $k(2)$. Each named class is a generator of a copy of $\mathbb{F}_2$.}
\label{fig:prismatic}
\end{figure}


\subsection{Syntomic cohomology of $k(n)$}\label{sec:syntomic-kn}

Now we have all the ingredients to compute the syntomic cohomology of Morava K-theory. 

\begin{definition}\label{MS}
Given $S\subset \{1,\cdots ,n\}$, we write 
\[ 
M_{S}:=
	\bF_p\{t^d\bar{\varepsilon}_S\lambda_{n+1} : -f(S) < d\le p^{n+1}-p^n-f(S)\}\,.
\]
using Definition~\ref{e(S)-f(S)}.
\end{definition}
With this notation in hand, we prove the following proposition, which will be important for understanding $\pi_*\grmot^*\TC^{-}(k(n))/(p,v_1,\cdots,v_{n+1})$.
\begin{prop}\label{prop:syntomic-A11}
The nontrivial classes in the vector spaces 
\begin{align}
\bigoplus_{S\subset \{1,\cdots , n\}} M_S 
\end{align}
using Definition~\ref{MS} are the only nontrivial classes that are hit by differentials that cross from Nygaard filtration $<0$ to Nygaard filtration $>0$ in the periodic $t$-Bockstein spectral sequence \eqref{eq: periodic t-bss}.
\end{prop}

\begin{proof}
By Theorem~\ref{thm:prismatic}, the periodic $t$-Bockstein spectral sequence collapses at the $\EE_{p^{n+1}+1}$-page and $d_r(x)\ne 0$ for $r\ge 2$ if and only if $d_r(x)=\lambda_{n+1}y$ for some $y\ne 0$. Therefore, the boundaries of differentials crossing from Nygaard filtration $<0$ to Nygaard filtration $>0$ are contained in the bigraded $\mathbb{F}_p$-vector space
\[
\bF_p\{t^{d+jp^{n+1}}\lambda_{n+1}\bar{\varepsilon}_S : 1\le d\le p^{n+1}-1\,, j\in \mathbb{Z} \,, S\subset \{1,2,\cdots ,n\}
\} \\
\]
using Definition~\ref{e(S)-f(S)}. 
We will now explain how the result follows from the proof of Theorem~\ref{thm:prismatic}. 

By $t^{p^{n+1}}$-periodicity, it suffices to consider the cases $t^{d}\lambda_{n+1}\bar{\varepsilon}_S$  for $0\le d\le p^{n+1}-1$. 

We know that the classes $t^d\lambda_{n+1}$ for $p^{n+1}-p^{n} <d\le p^{n+1}-1$ are not hit by differentials because the only possible sources of $d_r$ differentials  for $2\le r\le p^{n+1}$ that could hit these classes for bidegree reasons are the classes $t^{-d}$ for $1\le d<p^n$ and these classes are known to be permanent cycles by Proposition~\ref{prop:perm-cycles-tp}. Moreover, there are differentials 
\[
d_{p^{n+1}}(t^{j})=t^{p^{n+1}+j}\lambda_{n+1}
\] 
for each $-p^{n+1}< j \le  -p^{n}$. This proves the case $S=\emptyset$. 

If $S\ne \emptyset$ the class $t^{d}\bar{\varepsilon}_S\lambda_{n+1}$ for $-f(S)< d \le p^{n+1}-p^n-f(S)$ is hit by a differential
\[
d_{p^{n+1}}(t^{-p^{n+1}+d}
)=t^{d}\overline{\varepsilon}_i\lambda_{n+1}
\]
that crosses from strictly negative Nygaard filtration to strictly positive Nygaard filtration by Theorem~\ref{thm:prismatic}.  

For $0\le m\le p^{n+1}$, there are no further classes that can hit the classes $t^{p^{n+1}-m}\bar{\varepsilon}_{S}\lambda_{n+1}$ by consideration of bidegrees and permanent cycles. 
\end{proof}

\begin{defin}\label{Nyg}
Let 
\[ 
\textup{Nyg}_{\ge 1} \subset \pi_*\grmot^*\TC^{-}(k(n))/(p,v_1,\cdots,v_{n+1}))
\]
denote the subset of elements in Nygaard filtration $\ge 1$ and write 
\[
\textup{Nyg}_{=0} :=\pi_*\grmot^*\TC^{-}(k(n))/(p,v_1,\cdots,v_{n+1}))/\textup{Nyg}_{\ge 1} 
\]
for the canonical quotient. 
\end{defin}

\begin{proposition}\label{prop:Frobkn}
The Frobenius map 
\[ 
\varphi : \pi_*\grmot^*\TC^{-}(k(n))/(p,v_1,\cdots,v_{n+1})\to \pi_*\grmot^*\TP(k(n))/(p,v_1,\cdots,v_{n+1})
\]
satisfies $\varphi |_{\mathrm{Nyg}_{\ge 1}}=0$ and the induced map 
\[ 
\overline{\varphi} : \textup{Nyg}_{=0}\to \pi_*\grmot^*\TP(k(n))/(p,v_1,\cdots,v_{n+1}))
\]
is a monomorphism. 
\end{proposition}
\begin{proof}
This follows from Theorem~\ref{thm:Hodge-Tate}, Theorem~\ref{thm:Hodge-Tateresults}, and Definition~\ref{Nyg}. 
\end{proof}

\begin{defin}\label{def:ABCD}
There is a short exact sequence 
\[ 
0 \to T \to \pi_*\grmot^*\TC^{-}(k(n))/(p,v_1,\cdots,v_{n+1}) \to F \to 0 
\]
of $\Lambda(\lambda_{n+1})$-modules where $F$ is defined to be the image 
\[ 
F= \can \left (\pi_*\grmot^*\TC^{-}(k(n))/(p,v_1,\cdots,v_{n+1}) \right ) 
\]
inside of 
\[
\pi_*\grmot^*\TP(k(n))/(p,v_1,\cdots,v_{n+1})
\]
and $T$ is the kernel. This produces a diagram of short exact sequences 
\[
\begin{tikzcd}
& 0 \ar[d] & 0 \ar[d] & 0 \ar[d] & \\
0 \ar[r] & A_{11}\ar[d]  \ar[r] & \mathrm{Nyg}_{\ge 1} \ar[d]  \ar[r] & A_{10} \ar[d]   \ar[r] & 0 \\
0 \ar[r] &T \ar[r]\ar[d]  & \pi_*\grmot^*\TC^{-}(k(n))/(p,v_1,\cdots,v_{n+1} )\ar[d] \ar[r] &  F  \ar[d] \ar[r] & 0 \\
0 \ar[r] & A_{01} \ar[d] \ar[r] & \mathrm{Nyg}_{=0} \ar[d] \ar[r] &  A_{00} \ar[d] \ar[r] & 0 \\
& 0 & 0  & 0  & \\
\end{tikzcd}
\]
of $\Lambda(\lambda_{n+1})$-modules where 
\begin{align*}
A_{11}:=&\mathrm{Nyg}_{\ge 1}\cap T \\  
A_{10}:=& \mathrm{Nyg}_{\ge 1}/A_{11}\\
A_{01}:=&  T/A_{11}  \\
A_{00}:=& \mathrm{Nyg}_{=0}/A_{01}
\end{align*}
using Definitions~\ref{MS} and~\ref{Nyg}. 
\end{defin}

\begin{prop}
We can identify 
\begin{align*}
A_{00}& = \Lambda (\bar{\varepsilon}_1,\cdots ,\bar{\varepsilon}_n)\otimes \Lambda (\lambda_{n+1}) \\ 
A_{10} & = \bF_p\{\bar{\varepsilon}_{S}t^{j}  : S\subset \{1,\cdots ,n\}, -p^n-f(S)\le j <-f(S) \}
\otimes \bF_p[t^{p^{n+1}}]\{t^{p^{n+1}}\}\otimes \Lambda (\lambda_{n+1}) \\
A_{01} & =  \bF_p\{x\mu ,xz_{n+1}: x\in B \}  \otimes \Lambda( \lambda_{n+1})\otimes \bF_p[\mu^{p^{n+1}}] \\
A_{11}& =\bigoplus_{S\subset \{1,\cdots,n\}}M_S 
\end{align*}
and consequently there is an isomorphism of $\Lambda (\lambda_{n+1})$-modules
\[ 
\pi_{\ast}\grmot^{\ast}\TC^{-}(k(n))/(p,v_1,\cdots ,v_{n+1}) \cong A_{00} \oplus A_{01} \oplus A_{10} \oplus A_{11} \,.
\]
Here, elements are named according to where they are detected on the $\mathrm{E}_{\infty}$-page of the $t$-Bockstein spectral sequence, so that $\bar{\varepsilon}_i\in \{\bar{\varepsilon}_i\}$, $\bar{\varepsilon}_{S}t^{j} \in \{\bar{\varepsilon}_St^{j}\}$, $t^{p^{n+1}}\in \{t^{p^{n+1}}\}$, $t^d\bar{\varepsilon}_S\lambda_{n+1}\in \{t^d\bar{\varepsilon}_S\lambda_{n+1}\}$ and $\lambda_{n+1}$ is the unique class in $\{\lambda_{n+1}\}$ in the $t$-Bockstein spectral sequence.  The quantity $f(S)$ is specified according to Definition~\ref{e(S)-f(S)} and the notation $M_S$ is specified according to Definition~\ref{MS}. The notation $z_{n+1}$ is defined in Theorem~\ref{thm:Kn}. 
\end{prop} 

\begin{proof}
We first observe that $\lambda_{n+1}\in \mathrm{Nyg}_{=0}$ and $\lambda_{n+1}\in T$ and there is a unique class $\lambda_{n+1}$ in $\{\lambda_{n+1}\}$. By direct inspection of Theorem~\ref{thm:prismatic}, we can identify 
\[ 
A_{00}:=\Lambda( \bar{\varepsilon}_1,\cdots,\bar{\varepsilon}_n )\otimes \Lambda (\lambda_{n+1}) 
\]
as a free $\Lambda(\lambda_{n+1})$-module. Consequently, the exact sequences
\[ 
0 \to A_{01}\to \mathrm{Nyg}_{=0} \to A_{00} \to 0 
\]
and 
\[
0\to A_{10}\to F\to A_{00} \to 0
\]
split as $\Lambda (\lambda_{n+1})$-modules. Since $\mathrm{Nyg}_{=0}=K_n$ as defined in Definition~\ref{Kn-Knt} and described in Theorem~\ref{thm:Kn}, it is also a free $\Lambda( \lambda_{n+1})$-module. Consequently, the short exact sequence 
\[ 0 \to \mathrm{Nyg}_{\ge 1}\to \pi_{\ast}\grmot^{\ast}\TC^{-}(\BP)/(p,v_1,\cdots ,v_{n+1}) \to \mathrm{Nyg}_{=0}\to 0 
\]
is also split. From our prismatic cohomology computation in Theorem~\ref{thm:prismatic}, we can also identify $A_{10}$ and observe that it is free as a $\Lambda(\lambda_{n+1})$-module. Consequently, the short exact sequence 
\[ 0 \to A_{11}\to \mathrm{Nyg}_{\ge 1} \to A_{10} \to 0 
\]
also splits. We identified $A_{11}$ in Proposition~\ref{prop:syntomic-A11}. 
\end{proof}

\begin{figure}[ht!]
\resizebox{\textwidth}{!}{ 
\begin{tikzpicture}[radius=1,yscale=2]
\foreach \n in {-2,-1,...,27} \node [below] at (\n,-.8-1) {$\n$};
\foreach \s in {-1,0,...,12} \node [left] at (-.3-2,\s) {$\s$};
\draw [thin,color=lightgray] (-2,-1) grid (27,12);
\node [below] at (0,0) {$1$};
\node [below] at (1,1) {$t\bar{\varepsilon}_1$};
\node [below] at (3,0) {$\bar{\varepsilon}_1$};
\node [below] at (3,2) {$t^{2}\bar{\varepsilon}_2$};
\node [below] at (5,1) {$t\bar{\varepsilon}_2$};
\node [below] at (7,0) {$\bar{\varepsilon}_2$};
\node [below] at (4,3) {$t^{3}\bar{\varepsilon}_1\bar{\varepsilon}_2$};
\node [below] at (6,2) {$t^{2}\bar{\varepsilon}_1\bar{\varepsilon}_2$};
\node [below] at (8,1) {$t\bar{\varepsilon}_1\bar{\varepsilon}_2$};
\node [below] at (10,0) {$\bar{\varepsilon}_1\bar{\varepsilon}_2$};

\node [above] at (15,0) {$\lambda_3$};
\node [below] at (16,1) {$t\bar{\varepsilon}_1\lambda_{3}$};
\node [below] at (18,0) {$\bar{\varepsilon}_1\lambda_{3}$};
\node [below] at (18,2) {$t^{2}\bar{\varepsilon}_2\lambda_{3}$};
\node [below] at (20,1) {$t\bar{\varepsilon}_2\lambda_{3}$};
\node [below] at (22,0) {$\bar{\varepsilon}_2\lambda_{3}$};
\node [below] at (19,3) {$t^{3}\bar{\varepsilon}_1\bar{\varepsilon}_2\lambda_{3}$};
\node [below] at (21,2) {$t^{2}\bar{\varepsilon}_1\bar{\varepsilon}_2\lambda_{3}$};
\node [below] at (23,1) {$t\bar{\varepsilon}_1\bar{\varepsilon}_2\lambda_{3}$};
\node [below] at (25,0) {$\bar{\varepsilon}_1\bar{\varepsilon}_2\lambda_{3}$};

\node [above] at (-1,8) {$t^{8}\lambda_3$};
\node [above] at (1,7) {$t^{7}\lambda_3$};
\node [above] at (3,6) {$t^{6}\lambda_3$};
\node [above] at (5,5) {$t^{5}\lambda_3$};
\node [below] at (0,9) {$t^{9}\bar{\varepsilon}_1\lambda_{3}$};
\node [below] at (2,8) {$t^{8}\bar{\varepsilon}_1\lambda_{3}$};
\node [below] at (4,7) {$t^{7}\bar{\varepsilon}_1\lambda_{3}$};
\node [below] at (6,6) {$t^{6}\bar{\varepsilon}_1\lambda_{3}$};
\node [below] at (2,10) {$t^{10}\bar{\varepsilon}_2\lambda_{3}$};
\node [below] at (4,9) {$t^{9}\bar{\varepsilon}_2\lambda_{3}$};
\node [below] at (6,8) {$t^{8}\bar{\varepsilon}_2\lambda_{3}$};
\node [below] at (8,7) {$t^{7}\bar{\varepsilon}_2\lambda_{3}$};
\node [below] at (3,11) {$t^{11}\bar{\varepsilon}_1\bar{\varepsilon}_2\lambda_{3}$};
\node [below] at (5,10) {$t^{10}\bar{\varepsilon}_1\bar{\varepsilon}_2\lambda_{3}$};
\node [below] at (7,9) {$t^{9}\bar{\varepsilon}_1\bar{\varepsilon}_2\lambda_{3}$};
\node [below] at (9,8) {$t^{8}\bar{\varepsilon}_1\bar{\varepsilon}_2\lambda_{3}$};

\begin{scope}[color=blue]
\node [above] at (7,4) {$t^{4}\lambda_3$};
\node [above] at (9,3) {$t^{3}\lambda_3$};
\node [above] at (11,2) {$t^{2}\lambda_3$};
\node [above] at (13,1) {$t\lambda_3$};
\end{scope}
\begin{scope}[color=green]
\node [above] at (8,5) {$t^{5}\bar{\varepsilon}_{1}\lambda_3$};
\node [above] at (10,4) {$t^{4}\bar{\varepsilon}_{1}\lambda_3$};
\node [above] at (12,3) {$t^{3}\bar{\varepsilon}_{1}\lambda_3$};
\node [above] at (14,2) {$t^{2}\bar{\varepsilon}_{1}\lambda_3$};
\end{scope}
\begin{scope}[color=orange]
\node [above] at (10,6) {$t^{6}\bar{\varepsilon}_{2}\lambda_3$};
\node [above] at (12,5) {$t^{5}\bar{\varepsilon}_{2}\lambda_3$};
\node [above] at (14,4) {$t^{4}\bar{\varepsilon}_{2}\lambda_3$};
\node [above] at (16,3) {$t^{3}\bar{\varepsilon}_{2}\lambda_3$};
\end{scope}
\begin{scope}[color=purple]
\node [above] at (11,7) {$t^{7}\bar{\varepsilon}_{1}\bar{\varepsilon}_{2}\lambda_3$};
\node [above] at (13,6) {$t^{6}\bar{\varepsilon}_{1}\bar{\varepsilon}_{2}\lambda_3$};
\node [above] at (15,5) {$t^{5}\bar{\varepsilon}_{1}\bar{\varepsilon}_{2}\lambda_3$};
\node [above] at (17,4) {$t^{4}\bar{\varepsilon}_{1}\bar{\varepsilon}_{2}\lambda_3$};
\end{scope}
\end{tikzpicture}
}
\caption{The $\EE_\infty$-page of the $t$-Bockstein spectral sequence converging to $\pi_{*}\grmot^{*}\TC^{-}(k(2))/(2,v_{1},v_{2},v_{3})$. Each named class represents a copy of $\mathbb{F}_2$. We draw $M_{\emptyset}$ in \textcolor{blue}{blue}, $M_{1}$ in \textcolor{green}{green}, $M_{2}$ in \textcolor{orange}{orange}, and $M_{\{1,2\}}$ in \textcolor{red}{red}.}
\label{fig:TC-}
\end{figure}


\begin{thm}\label{thm:syntomic-cohomology-kn}
The syntomic cohomology of $k(n)$ modulo $(p,v_{1},\cdots ,v_{n+1})$ is isomorphic to
\begin{align}\label{eq: syntomic}
\Lambda(\delta ,\bar{\varepsilon}_{1},\cdots,  \bar{\varepsilon}_{n})\otimes \Lambda (\lambda_{n+1})\oplus \bigoplus_{S\subset \{1,\cdots,n\}}M_S
\end{align}
as a $ \Lambda (\lambda_{n+1})$-module, where $M_{S}$ is defined in Definition~\ref{MS}. 
The bidegrees are 
\begin{align*}
\|\lambda_{n+1}\|& =(2p^n-1,1)\\
\| \partial \| & = (-1,1)  \\
\| \bar{\varepsilon}_i\| & = (2p^i-1,1) \,, \text{ for } 1\le i \le n \\
\|t^d\bar{\varepsilon}_S\lambda_{n+1}\| & =(\sum_{s\in S}2p^s-|S|-2d+2p^{n+1}-1,1-|S|)
\end{align*}
for $S\subset \{1,\cdots ,n\}$. In particular, in Adams weights $\le -n$ the syntomic cohomology is
\[
\mathbb{F}_{p}\{\bar{\varepsilon}_{1}\cdot \ldots \cdot\bar{\varepsilon}_{n}\}
\]
concentrated in bidegree $(\sum_{s=1}^n2p^s-n,-n)$. 
\end{thm}

\begin{proof}
In this proof, we will write 
\[
\pi_*\grmot^*\TC^{-} := \pi_*\grmot^*\TC^{-}(k(n))/(p,v_1,\cdots,v_{n+1})
\]
and 
\[
\pi_*\grmot^*\TP := \pi_*\grmot^*\TP(k(n))/(p,v_1,\cdots,v_{n+1})
\]
for conciseness. We observe that $\mathrm{can}(\mathrm{Nyg}_{\ge 1})$ consists exactly of  elements of $\pi_*\grmot^*\TP$ in Nygaard filtration $\ge 1$, so we write $\mathrm{Nyg}_{\ge 1}^{t}$ for this subgroup. We consider the diagram 
\[
\begin{tikzcd}
\mathrm{ker}_{\ge 1}\ar[r, rightarrowtail] \ar[d,rightarrowtail]   & \mathrm{Nyg}_{\ge 1} \ar[r,"\textup{can}|-0", twoheadrightarrow]\ar[d,rightarrowtail] & \mathrm{Nyg}_{\ge 1}^t \ar[d,rightarrowtail] \ar[r] &  0 \ar[d] \\
\mathrm{ker}(\textup{can}-\varphi)  \ar[r] \ar[d,twoheadrightarrow]  &  \pi_*\grmot^*\TC^{-} \ar[r,"\textup{can}-\varphi"] \ar[d,twoheadrightarrow]  & 
\pi_*\grmot^*\TP 
 \ar[r,twoheadrightarrow] \ar[d,twoheadrightarrow] & \mathrm{coker}(\textup{can}-\varphi) \ar[d, equal] \\ 
\mathrm{ker}_{=0}\ar[r, rightarrowtail]  & \mathrm{Nyg}_{=0} \ar[r,"\overline{\textup{can}}-\varphi"] & \pi_*\grmot^*\TP/\mathrm{Nyg}_{\ge 1}^t \ar[r,twoheadrightarrow] & \mathrm{coker} (\textup{can}-\varphi)
\end{tikzcd}
\]
where each column and row is exact. We claim that we can identify 
\[ 
\pi_*\grmot^*\TP \cong A_{00}\oplus A_{01}\oplus A_{10}
\]
as a $\Lambda(\lambda_{n+1})$-module. To see this, consider the short exact sequences 
\[ 
0 \to \mathrm{Nyg}_{\ge 1}^t \to \pi_*\grmot^*\TP  \to  \mathrm{Nyg}_{\le 0}^t  \to 0 \,,
\]
where $\mathrm{Nyg}_{\le 0}^t:=\pi_*\grmot^*\TP/\mathrm{Nyg}_{\ge 1}^t$ is defined as the cokernel of the canonical inclusion, and  
\[
0 \to \mathrm{Nyg}_{=0}^t \to \mathrm{Nyg}_{\le 0}^t \to   \mathrm{Nyg}_{<0}^t \to 0
\]
where $\mathrm{Nyg}_{=0}^t$ consists of elements in $\pi_*\grmot^*\TP$ in Nygaard filtration exactly zero. Observe that $\mathrm{Nyg}_{<0}^t\cong A_{01}$, $\mathrm{Nyg}_{=0}^t\cong A_{00}$, and $\mathrm{Nyg}_{\ge 1 }^t\cong A_{10}$ as $\Lambda (\lambda_{n+1})$-modules. The splitting then follows because  $A_{00}$ and $A_{01}$ are free $\Lambda( \lambda_{n+1})$-modules. We can therefore identify 
$\mathrm{ker}_{\ge 1} = A_{11}$ and $\mathrm{Nyg}_{\le 0}^t= A_{00}\oplus A_{01}$. 

Using our identifications, we also determine that 
\[
\mathrm{can}(A_{00}) =A_{00}
=\varphi(A_{00})
\]
so $A_{00}\subset \mathrm{ker}_{=0}$ and $A_{00}\subset \mathrm{coker}(\textup{can}-\varphi)$. We also observe that 
$\varphi(A_{01})=A_{01}$ and $\mathrm{can}(A_{01})=0$ so we determine that $\mathrm{ker}_{=0}=A_{00}$. 

We observed that $\varphi|_{A_{01}}$ is a bijection onto its image and $\mathrm{can}(A_{01})=0$,  $\varphi(A_{10})=0$ and $\can |_{A_{10}}$ is a bijection onto its image, and $\can(A_{00})=A_{00}=\varphi(A_{00})$
so $\overline{\can}-\varphi$ factors as 
\[ 
\begin{tikzcd}
\mathrm{Nyg}_{=0} \ar[r,twoheadrightarrow] & A_{10} \subset \mathrm{Nyg}_{\le 0}^t
\end{tikzcd}
\]
where 
$\begin{tikzcd}
\mathrm{Nyg}_{=0} \ar[r,twoheadrightarrow] &A_{10} \end{tikzcd}$ 
is a surjection. Consequently, we determine that 
\[
A_{00}=\mathrm{coker}(\overline{\textup{can}}-\varphi)= \mathrm{coker}(\textup{can}-\varphi)\,.
\]

Since $\ker_{=0}=A_{00}$ is a free $\Lambda (\lambda_{n+1})$-module, we determine that 
\[ 
\mathrm{ker}(\textup{can}-\varphi)=A_{00}\oplus A_{11}
\]
as $\Lambda (\lambda_{n+1})$-modules and the extension of $\Lambda(\lambda_{n+1})$-modules 
\[ \Sigma^{-1,1}A_{00}=\mathrm{coker} (\textup{can}-\varphi) \to \pi_*\gr_{\mot}^*\TC(k(n))/(p,v_1,\cdots ,v_n) \to \mathrm{ker} (\textup{can}-\varphi)=A_{00}\oplus A_{11} \]
 is split. This proves the claim, where $\partial$ denotes the image of $\Sigma^{-1,1}1\in \Sigma^{-1,1}A_{00}$. 
\end{proof}


\begin{figure}[ht!]
\resizebox{\textwidth}{!}{ 
\begin{tikzpicture}[radius=1,yscale=2]
\foreach \n in {-2,-1,...,26} \node [below] at (\n,-.8-2) {$\n$};
\foreach \s in {-2,-1,...,3} \node [left] at (-.3-2,\s) {$\s$};
\draw [thin,color=lightgray] (-2,-2) grid (26,3);
\node [below] at (5,-1) {$\bar{\varepsilon}_1$};

\node [below] at (0,0) {$1$};
\node [below] at (4,0) {$\partial \bar{\varepsilon}_1$};
\node [below] at (6,0) {$t^7\bar{\varepsilon}_1\lambda_2$};
\node [below] at (8,0) {$t^7\bar{\varepsilon}_1\lambda_2$};
\node [below] at (10,0) {$t^6\bar{\varepsilon}_1\lambda_2$};
\node [below] at (12,0) {$t^5\bar{\varepsilon}_1\lambda_2$};
\node [below] at (14,0) {$t^4\bar{\varepsilon}_1\lambda_2$};
\node [below] at (16,0) {$t^3\bar{\varepsilon}_1\lambda_2$};

\node [above] at (22,0) {$\bar{\varepsilon}_1\lambda_2$};

\node [below] at (-1,1) {$\partial$};

\node [above] at (5,1) {$t^6\lambda_3$};
\node [above] at (7,1) {$t^5\lambda_3$};
\node [above] at (9,1) {$t^4\lambda_3$};
\node [above] at (11,1) {$t^3\lambda_3$};
\node [above] at (13,1) {$t^2\lambda_3$};
\node [above] at (15,1) {$t\lambda_2$};
\node [above] at (17,1) {$\lambda_2$};
\node [above] at (21,1) {$\partial\bar{\varepsilon}_1\lambda_3$};

\node [above] at (16,2) {$\partial\lambda_2$};
\end{tikzpicture}
}
\caption{The mod $(3,v_1,v_2)$-syntomic cohomology of $k(1)$}
\label{fig:syntomic-k1}
\end{figure}

\section{Algebraic K-theory and topological cyclic homology}\label{sec:k-theory}
In Section~\ref{sec:low-heights}, we compute the mod $(p,v_1)$ algebraic K-theory of $k(1)$ at all primes $p\ge 3$, extending work of~\cite{AR12}. We then prove the Lichtenbaum---Quillen property for the algebraic K-theory of $k(n)$ in Section~\ref{sec:LQ}. We prove that the telescope conjecture is true for the algebraic K-theory of $k(n)$ in Section~\ref{telescope}. We prove redshift for $k(n)$ in Section~\ref{sec:redshift}. We conclude by discussing the case of Morava K-theory over more general perfect fields in Section~\ref{sec:perfect-fields}.

\subsection{Algebraic K-theory at low heights}\label{sec:low-heights}
In this section, we extend a previous computation of Ausoni--Rognes~\cite{AR12} regarding $V(1)_*\TC(k(1))$. Note that, at primes $p \ge 3$, there exists a self map 
\[\Sigma^{2p-2}\mathbb{S}/p \to \mathbb{S}/p
\]
inducing multiplication by $v_1$ on $K(1)_*$ whose cofiber is denoted $V(1)$ ~\cite{Ada66,Tod71}. When $p\ge 5$ there exists a self map 
\[ v_2 : \Sigma^{2p^2-2}V(1)\to V(1)
\]
inducing multiplication by $v_2$ on $K(2)_*$, and at $p=3$ there exists a self-map 
\[ v_2^9 : \Sigma^{36}V(1)\to V(1)
\]
inducing $v_2^9$ on $K(2)_*$ \cite{Tod71,BP04}. 

If $p\ge 5$, then $V(1)$ is homotopy commutative and associative by~\cite{Oka84}, but a key feature of our work is that such structure is not necessary for our computation. In particular, we determine that 
\[
\gr_{\ev}^*\bS/(p,v_1) =\gr_{\ev/\bS}^*V(1)
\]
as in~\cite[Remark~6.0.5.]{HRW22} and $\gr_{\ev}^*\bS/(p,v_1)$ is an $\bE_\infty$-$\gr_{\ev}^*\bS$-algebra for $p\ge 3$. 

As in Section~\ref{sec:setup}, there is a motivic spectral sequence 
\[
\EE_2^{s,2t-s}=\pi_s\grmot^t\TC(k(1))/(p,v_1)\implies V(1)_s\TC(k(1))   
\]
with differentials 
\[ 
d_r : \EE_2^{s,2t-s} \longrightarrow  \EE_2^{s-1,2t-s+r}
\]
for $r\ge 2$ associated to the filtration
\[ \fil_{\ev/\mathbb{S}}^*V(1)\otimes_{\fil_{\ev}^*\mathbb{S}}\fil_{\mot}^*\TC(k(1)) 
\]
on $V(1)\otimes \TC(k(1))$. We also note that at all primes $p$ there is a $v_2$-Bockstein spectral sequence 
\[ E_2^{s,2t-s,f}=\pi_s\grmot^t\TC(k(1)/(p,v_1,v_2)[v_2] \implies \pi_s\grmot^t\TC(k(1)/(p,v_1)
\]
with differentials satisfying
\[ d_r :E_2^{s,2t-s,f} \to E_2^{s-1,2t-s+1,f+r}, 
\]
associated to the tower 
\[
\begin{tikzcd}
\dots \ar[r,"v_2"] & \grmot^*\TC(k(1))/(p,v_1) \ar[r,"v_2"]\ar[d]  & \grmot^*\TC(k(1))/(p,v_1)  \ar[d]  \\
& \grmot^*\TC(k(1))/(p,v_1,v_2) & \grmot^*\TC(k(1))/(p,v_1,v_2)\,.  
\end{tikzcd}
\]
\begin{thm}\label{thm:tck1}
At $p\ge 5$, there is a preferred isomorphism
\begin{align*} 
\pi_*\TC(k(1))/(p,v_1)\cong& \bF_{p}[v_{2}]\otimes \Lambda(\delta ,\bar{\varepsilon}_1,\lambda_{2})\oplus  \\
&\bF_{p}[v_{2}]\otimes    \bF_{p}\{ 
t^{k}\bar{\varepsilon}_1\lambda_2
 :  p-1 < k \le p^{2}-1\}\oplus \\
& \bF_{p}[v_{2}]\otimes \bF_p\{t^{k}\lambda_2 : 0 <k\le p^{2}-p \} 
\end{align*}
of $\bF_p[v_2]$-modules. At $p=3$, there is a preferred isomorphism  
\begin{align*} 
\pi_*\TC(k(1))/(3,v_1)\cong& \bF_{3}[v_{2}^9]\otimes \bF_{3}[v_{2}]/v_2^9\otimes \Lambda(\delta ,\bar{\varepsilon}_1,\lambda_{2})\oplus  \\
&\bF_{3}[v_{2}^9]\otimes \bF_{3}[v_{2}]/v_2^9\otimes   \bF_{3}\{ t^{k}\bar{\varepsilon}_1\lambda_2  :  2< k\le 8\}\oplus \\
& \bF_{3}[v_{2}^9]\otimes \bF_{3}[v_{2}]/v_2^9\otimes \bF_p\{t^{k}\lambda_2 : 0 <k\le 6 \} 
\end{align*}
$\bF_3[v_2^{9}]$-modules. Here we write $t^{k}\overline{\varepsilon}_1\lambda_2$ and $t^k\lambda$ for the classes detected by $t^{k}\overline{\varepsilon}_1\lambda_2$ and $t^k\lambda_2$ respectively in the $\EE_\infty$-page of the $t$-Bockstein spectral sequence.
\end{thm}
\begin{proof}
There are no possible $v_{2}$-Bockstein spectral sequence differentials for bidegree reasons. The reader can observe this directly at the prime $p=3$ in Figure~\ref{fig:syntomic-k1}. 

The mod $(p,v_1)$-syntomic cohomology of $k(1)$ is therefore
\begin{align}
\pi_*\grmot^*\TC(k(1))/(p,v_1) \cong \mathbb{F}_{p}[v_{2}] \otimes \left ( \Lambda(\partial ,\bar{\varepsilon}_{1} ,\lambda_{2})\oplus M_{\emptyset}\oplus M_{1} \right )
\end{align}
for all primes $p$ where
\begin{align*}
M_{\{1\}}&= \bF_{p}\{ t^{k}\bar{\varepsilon}_1\lambda_2 :  p-1 < d \le p^{2}-1 \} \text{ and } \\
M_{\emptyset}&=\bF_p\{t^k\lambda_2 : 0< k\le p^{2}-p \} \,,
\end{align*}
using Definition~\ref{MS}.

The motivic spectral sequence is concentrated in Adams weights $[-1,2]$. There are no possible $d_2$-differentials because the $-1$-line and the $1$-line are concentrated in odd degrees and the $0$-line and the $1$-line are concentrated in even degrees. It therefore suffices to check that there are no possible $d_3$-differentials. To see this, note that the $-1$-line is supported in degrees $2p-1+(2p^{2}-2)j$ for $j\ge 0$ and the $2$-line is supported in degrees $(2p^2-2)\ell$ for $\ell\ge 1$ so there are no possible $d_3$-differentials for bidegree reasons. Since $\TC(k(1))/(p,v_1)$ is a $V(1)\otimes D(V(1))$-module and $V(1)$ has a $v_2^1$-self map when $p\ge 5$ and a $v_2^9$-self map when $p=3$, the result follows. 
\end{proof}

\begin{thm}\label{thm:kk1}
At $p\ge 5$, there is preferred isomorphism of $\bF_p[v_2]$-modules 
\begin{align*}
V(1)_*\K(k(1))\cong & \bF_p[v_2]\otimes (\Lambda (\bar{\varepsilon}_1)\otimes \mathbb{F}_p\{1,\partial \lambda_2,\lambda_2,\partial v_2
\} \oplus M_{\emptyset}\oplus M_{\{1\}})\,.
\end{align*}
At $p=3$, there is an isomorphism of $\bF_3[v_2^9]$-modules
\begin{align*}
V(1)_*\K(k(1))\cong & \mathbb{F}_3[v_2^9] \otimes \left ( \mathbb{F}_3[v_2]/v_2^9\otimes (\Lambda (\bar{\varepsilon}_1)\otimes \mathbb{F}_p\{1,\partial \lambda_2,\lambda_2,\partial v_2
\} \oplus M_{\emptyset}\oplus M_{\{1\}}) \right )\,.
\end{align*}
Here $M_{S}$ for $S\subset \{1\}$ is defined as in Definition~\ref{MS}. 
\end{thm}

\begin{proof}
The proof of \cite[Theorem~7.7]{AR12} applies using Theorem~\ref{thm:tck1} in place of~\cite[Theorem~7.6]{AR12}, but we repeat the proof for completeness. By~\cite[Theorem 7.3.1.8]{DGM13}, there is a map of fiber sequences 
\[ 
\begin{tikzcd}
V(1)\otimes \mathrm{K}(k(1))\ar[r] \ar[d] & V(1)\otimes \mathrm{TC}(k(1))\ar[r] \ar[d]& V(1)\otimes \Sigma^{-1}H\mathbb{Z}_p \ar[d,"="] \\ 
V(1)\otimes \mathrm{K}(\mathbb{F}_p )\ar[r] & 
V(1)\otimes \mathrm{TC}(\mathbb{F}_p)\ar[r] & V(1)\otimes \Sigma^{-1}H\mathbb{Z}_p  
\end{tikzcd}
\]
and the bottom fiber sequence induces the split exact sequence  
\[
\Lambda(\bar{\varepsilon}_1)\to \Lambda(\partial, \bar{\varepsilon}_1) \to \Sigma^{-1}\Lambda(\bar{\varepsilon}_1)
\]
on homotopy. The reduction map $k(1)\to \mathbb{F}_p$ is $2p-3$-connective so the induced map
\[
V(1)\otimes F(k(1))\to V(1)\otimes F(\mathbb{F}_p)
\] 
is $2p-2$-connective for $F\in \{\K,\TC\}$ by~\cite[Proposition~10.9]{BM94}. Since $\Sigma^{-1}\Lambda (\bar{\varepsilon}_1)$
is concentrated in degrees $\le 2p-2$, 
we determine that the map 
\[ 
V(1)_*\TC(k(1))\to \Sigma^{-1}\Lambda (\bar{\varepsilon}_1)
\]
is surjective with kernel given by the stated formula.
\end{proof}

\begin{remark}
This extends~\cite[Theorem~7.6]{AR12} to arbitrary $\bE_1$-$\MU$-algebra forms of $k(1)$ in the sense of Definition~\ref{def:form-kn} and to the prime $p=3$. 
To compare notation, note that 
\[\bF_p\{t^{k}\bar{\varepsilon}_{1}\lambda_{2} :  p-1< k \le p^{2}-1\}=\bF_p\{t^{dp}\bar{\varepsilon}_{1}\lambda_{2} :  0<d<p \} \oplus \bF_p\{t^{d}v_{2}: 0<d<p^{2}-p : p \nmid d \} \]
where we write $t^{p+d}\bar{\varepsilon}_{1}\lambda_{2}$ for the class $t^{d}v_{2}$ (up to multiplication by a unit) and 
\[
\bF_p\{t^{k}\lambda_{2} : 0< k\le p^{2}-p \} = \bF_p\{t^{dp}\lambda_{2} : 0<d<p \}\oplus   \bF_p\{\textup{dlog}v_{1}\cdot t^{d}v_{2} : 0<d<p^{2}-p : p \nmid d \}
\]
in light of the relation $\textup{dlog} v_{1} \cdot t^{d}v_{2}=t^d\lambda_2$ from \cite[Theorem~7.6]{AR12}.
\end{remark}

\subsection{Lichtenbaum--Quillen}\label{sec:LQ}
In this section, we note that our computations imply that the Lichtenbaum--Quillen property holds for $\K(k(n))$. 
We can fix a choice of integers 
\[
(i_0,i_1,\cdots,i_n,i_{n+1})
\]
such that the finite complex 
\[
\bS/(p^{i_0},v_1^{i_1},\cdots,v_n^{i_n},v_{n+1}^{i_{n+1}}) 
\]
exists by~\cite[Theorem~9]{HS98}. 
\begin{theorem}\label{thm:finite-TC}
Let $n\ge 1$. The graded abelian group 
\[
\pi_*\TC(k(n))/(p^{i_0},v_1^{i_1},\cdots,v_n^{i_n},v_{n+1}^{i_{n+1}})
\] 
is finite. Consequently, the spectrum $\TC(k(n))$ has fp-type $n+1$.
\end{theorem}
\begin{proof}
We determined that 
\[ \pi_*\grmot^*\TC(k(n))/(p,v_1,\cdots,v_{n+1})\]
is finite. This immediately implies that 
\[ \pi_*\grmot^*\TC(R)/(p^{i_0},v_1^{i_1},\cdots,v_{n+1}^{i_{n+1}})
\]
is finite for any integers $i_0,i_1,\cdots, i_{n+1}$ by using the finite length filtrations that produce $p,v_1,\cdots ,v_{n+1}$-Bockstein spectral sequences. This immediately implies that 
\[
\pi_*\TC(k(n))/(p^{i_0},v_1^{i_1},\cdots,v_{n+1}^{i_{n+1}})
\]
is finite since it is a sub-quotient of a finite graded abelian group. Consequently $\TC(k(n))$ has fp-type $n+1$ by~\cite[Proposition~3.2]{MR99}. 
\end{proof}
This answers Question~\ref{quest:fp-type} in the affirmative for all $\mathbb{E}_1$-$\MU$-algebra forms of $k(n)$ at all primes $p$ and heights $n\ge 1$. This also implies the generalization of the version of the Lichtenbaum--Qullen conjecture suggested in~\cite[Conjecture~4.3]{AR08}. 

\begin{cor}\label{cor:LQC-K-theory}
The fiber of the map
\[ \K(k(n))_{(p)}\longrightarrow  L_{n+1}^{f}\K(k(n))_{(p)} \]
is bounded above for integers $n\ge 1$ and primes $p$.
\end{cor}

\begin{proof}
This follows from Theorem~\ref{thm:finite-TC} by~\cite[Theorem~3.1.3]{HW22}. 
\end{proof} 

\begin{cor}[Lichtenbaum--Quillen for algebraic K-theory of Morava K-theory]\label{cor:LQC-K-theory-periodic}
The fiber of the canonical map 
\[ \K(K(n))_{(p)}\longrightarrow  L_{n+1}^{f}\K(K(n))_{(p)} \]
is bounded above for integers $n\ge 1$ all primes $p$. 
\end{cor}
\begin{proof}
This follows from Corollary~\ref{cor:LQC-K-theory} together with the fiber sequence 
\[ \K(\bF_p)\to  \K(k(n))\to \K(K(n))
\]
of~\cite{BM08} (cf.~\cite[8.8.4]{BarwickHeart}) and the fact that the localization map
\[ H\mathbb{Z}_{(p)}\simeq \K(\bF_p)_{(p)}\to L_{n+1}^f\K(\bF_p)_{(p)}\simeq L_0\K(\bF_p)_{(p)}\simeq H\mathbb{Q}
\]
has bounded above fiber by~\cite{Qui72}. 
\end{proof}

\subsection{Telescope}\label{telescope}
We now note that our results also imply the telescope conjecture for the topological cyclic homology of $k(n)$. In fact, we prove a more general statement. 

\begin{defin}
Let $E$ be an $\bE_{\infty}$-ring. Given an $\bE_1$-$E$-algebra $R$, we say that the $E$-based motivic cohomology of $R$ has a horizontal vanishing line if there exists an integer $N$ such that $\pi_*\gr_{\mot/E}^*\THH(R)$ vanishes in Adams weights $\ge N$. 
\end{defin}

\begin{prop}\label{prop:telescope}
Let $E$ be an $\bE_{\infty}$-ring such that there exists an eff map of cyclotomic $\bE_{\infty}$-rings $\THH(E)\to B$ where $B$ is even. Let $R$ be a connective $\bE_1$-$E$-algebra such that the $E$-based motivic cohomology of $R$ has a horizontal vanishing line.
Then the canonical map 
\[ 
L_{n+1}^f\K(R)_{(p)}\to L_{n+1}\K(R)_{(p)}
\]
is an equivalence. 
\end{prop}

\begin{proof}
By the arithmetic fracture square 
\[
\begin{tikzcd}
\K(R)_{(p)} \ar[r] \ar[d] & \K(R)_{p}^{\wedge}  \ar[d] \\ 
L_0\K(R) \ar[r]  & L_0\K(R)_{p}^{\wedge}
\end{tikzcd}
\]
and the fact that $L_m^fX\simeq L_nX$ for all $\bQ$-modules $X$, it suffices to prove the result for $\K(R)_p^\wedge$. 

We know the height-$1$ telescope conjecture holds by~\cite{Mil81,Mah82}. 
Together with~\cite{Mit90}, this implies that 
\[ 
L_m^f\K(\pi_0R)_p^{\wedge}\to  L_m\K(\pi_0R)_p^{\wedge} \text{ and }
L_m^f\TC(\pi_0R)_p^{\wedge}\to  L_m\TC(\pi_0R)_p^{\wedge}
\]
are equivalences for all $m\ge 0$.
It therefore suffices to prove that 
\[ 
L_{n+1}^f\TC(R)_p^{\wedge}\to  L_{n+1}\TC(R)_p^{\wedge}
\]
is an equivalence. 

The horizontal vanishing line implies that $\pi_*\gr_{\mot/E}^*\TC^{-}(R)_p^{\wedge}$ and $\pi_*\gr_{\mot/E}^*\TP(R)_p^{\wedge}$ also vanish in Adams weights $\ge N$. 
Since there is an eff map $\THH(E)\to B$ of $S^1$-equivariant $\bE_{\infty}$-rings where $B$ is even, and we assumed that $\THH(R)$ is bounded below, this implies that the cobar complexes $\pi_{2*}X^{\bullet}$ and $\pi_{2*}Y^{\bullet}$ have bounded cohomology where 
\[
X^{\bullet}:=\TC^{-}(R \otimes_{E}B^{\otimes_{E}\bullet+1})
\]
and 
\[
Y^{\bullet}:=\TP(R \otimes_{E}B^{\otimes_{E}\bullet+1}).
\]
By~\cite{CM15}, we determine inductively that for each $2\le m\le n+1$ and integers $i_0,i_1,\cdots,i_m$ such that $\bS/(p^{i_0},\cdots,v_m^{i_m})$ exists we have 
\begin{align*}
v_{m+1}^{-1}\mathbb{S}/(p^{i_0},\cdots ,v_m^{i_m})\otimes \mathrm{Tot}\thinspace \bigl( W^{\bullet} \bigr) & \simeq \mathrm{Tot}\thinspace v_{m+1}^{-1} \bigl( W^{\bullet} /(p^{i_0},\cdots ,v_m^{i_m})\bigr) \\
\simeq \mathrm{Tot}\thinspace L_n \bigl( W^{\bullet} /(p^{i_0},\cdots ,v_m^{i_m})\bigr) 
\end{align*}
for $W^{\bullet}\in \{X^{\bullet},Y^{\bullet}\}$. The last equivalence holds since $B$ is an even $\bE_{\infty}$-ring and consequently any $B$ module is an $\MU$-module and we know the telescope conjecture holds for $\MU$-modules. 

Since limits of $L_m$-local spectra are $L_m$-local, we determine the telescope conjecture holds for $\TC^{-}(R)_p^{\wedge}$ and $\TP(R)_p^{\wedge}$ for all $2\le m \le n+1$. By the five lemma and~\cite[Corollary~1.5]{NS18}, this implies the height-$n+1$ telescope conjecture for $\TC(R)_p^{\wedge}$ as desired. 
\end{proof}

\begin{thm}[Telescope for algebraic K-theory of Morava K-theory]\label{thm:telescope}
The canonical map 
\[ 
L_{n+1}^f\K(k(n))_{(p)}\to L_{n+1}\K(k(n))_{(p)}
\]
is an equivalence for all integers $n\ge 1$ and primes $p$. 
\end{thm}

\begin{proof}
There is a cyclotomic map~$\THH(\MU)\to \THH(\MU/\MW)$ where $\THH(\MU/\MW)$ is even by~\cite[Theorem~3.2.10]{HRW22}. This map is also eff by basechange along the map $\THH(\MW)\to \MW$, which is eff by the same argument as~\cite[Example~4.2.3]{HRW22}. Since $k(n)$ is connective, it is clear that it has bounded below topological Hochschild homology. 
We determined that the $\MU$-based motivic cohomology of $k(n)$ has a horizontal vanishing line in Corollary~\ref{cor:vanishing-line}. The result then follows from Proposition~\ref{prop:telescope}.
\end{proof}

\begin{cor}\label{thm:telescope-periodic}
The canonical map 
\[ 
L_{n+1}^f\K(K(n))_{(p)}\to L_{n+1}\K(K(n))_{(p)}
\]
is an equivalence for all integers $n\ge 1$ and primes $p$. 
\end{cor}
\begin{proof}
This follows from the localization sequence 
\[\K(\bF_p)\to \K(k(n))\to \K(K(n))
\]
of~\cite{BM08} (cf.~\cite[8.8.4]{BarwickHeart}) and the fact that 
\[
L_0\K(\bF_p)_{(p)}\simeq L_{n+1}^f\K(\bF_p)_{(p)}\simeq L_{n+1}\K(\bF_p)_{(p)} 
\]
by~\cite{Qui72}. 
\end{proof}

\subsection{Redshift}\label{sec:redshift}
By Theorem~\ref{thm:telescope}, the canonical map  
 \[
 L_{T(n+1)} \K(k(n))\overset{\simeq}{\longrightarrow} L_{K(n+1)}\K(k(n))
 \] 
 is an equivalence. In this section, we prove that $L_{T(n+1)} \K(k(n))$ is non-zero.  In particular, we show that algebraic K-theory of Morava $K$-theory redshifts.  Our argument is substantially more elaborate than known arguments proving that Morava $E$-theories and truncated Brown--Peterson spectra redshift.\footnote{There are essentially two proofs of redshift for Morava $E$-theories. The first, due to Allen Yuan \cite{Yua21}, crucially uses the universal property of $\THH$ of $\mathbb{E}_{\infty}$-rings. The second is to first prove redshift for truncated Brown--Peterson spectra, which is accomplished in~\cite{HW22} by comparison with $\TC^{-}_*(\BPn/\MU)$. The root adjunction formalism of \cite{ABM22} then allows one to deduce redshift for Morava $E$-theories, as explained in \cite[\S 8]{ABM22}.} 
Unfortunately, any direct attempt to imitate those arguments seems fraught with subtle difficulties.\footnote{For example, we do not fully understand the $\MU_*$-module structure on $\TC^{-}_*(k(n)/\MU)$.  The $\MU_*$-module structure on $\TC^{-}_*(\BP\langle n \rangle/\MU)$ is substantially simpler, because there $(p,v_1,\cdots,v_n,v_{n+1})$ is a regular sequence.} 

So far we have studied syntomic cohomology modulo $(p,v_1,\cdots,v_n,v_{n+1})$.  To prove redshift, it will be important also to have some understanding of syntomic cohomology without coefficients, and with coefficients modulo $(p^{a_0},v_1^{a_1},\cdots,v_n^{a_n})$.

\begin{proposition}\label{prop:integral-collapse}
The motivic spectral sequence converging to $\pi_*\TC(k(n))_p^{\wedge}$ has $\mathrm{E}_2$-page supported on the $0$-line, $1$-line, and $2$-line.  In particular, it collapses at the $\mathrm{E}_2$-page.
\end{proposition}

\begin{proof}
By Corollary~\ref{cor:vanishing-line}, the motivic spectral sequence converging to $\pi_*\THH(k(n))$ admits a horizontal vanishing line at the $\mathrm{E}_2$-page, such that all classes are supported on the $0$-line and the $1$-line. Thus, the motivic spectral sequences converging to $\pi_*\TC^{-}(k(n))$ and to $\pi_*\TP(k(n))$ are also supported on the $0$-line and the $1$-line. It then formally follows that all classes on the $\mathrm{E}_2$-page of the motivic spectral sequence for $\pi_*\TC(k(n))$ are supported on lines $0$, $1$, and $2$.  By parity considerations, no differentials are possible.
\end{proof}

\begin{cor}\label{cor:connectivity}
In degrees $\le 2p^n-2$, the $\mathrm{E}_2$-page of the motivic spectral sequence for $\pi_*\TC(k(n))_p^{\wedge}$ contains only a copy of $\mathbb{Z}_p$ in bidegree $(0,0)$ and $\mathbb{Z}_{p}$ in bidegree $(-1,1)$. 
\end{cor}

\begin{proof}
The reduction map $k(n)\to \bF_{p}$ is $(2p^n-3)$-connective, so the induced map $\TC(k(n))\to \TC(\bF_{p})$ is $(2p^n-2)$-connective by~\cite[Proposition~10.9]{BM94} and~\cite[Theorem~7.2.2.1]{DGM13}. 

We know $\pi_*\grmot^*\TC(\mathbb{F}_p)_p^{\wedge}=\Lambda(\partial)$ where $\| \partial \| =(-1,1)$. By Proposition~\ref{prop:integral-collapse} and the fact that maps of motivic spectral sequences cannot lower motivic filtration, we see that the map on $\EE_{2}$-pages of motivic spectral sequences induced by the reduction map $k(n)\to \bF_{p}$ is an isomorphism in in bidegrees $(a,b)$ where $a\le 2p^{n}-2$. 
\end{proof}

In previous sections, we computed
\[\pi_*\gr^*_{\mot} \mathrm{TC}(k(n)) / (p,v_1,\cdots,v_n,v_{n+1}).\]

In order to study redshift, it is necessary to understand something about the $\mathbb{F}_p[v_{n+1}]$-module structure on
\[\pi_* \gr^*_{\mot} \TC(k(n) / (p,v_1,...,v_n),\]
which we can access via a $v_{n+1}$-Bockstein spectral sequence.

\begin{theorem}[Unital redshift in syntomic cohomology] \label{thm:unital-syntomic-redshift}
In 
\[\pi_*\gr^*_{\mot} \TC(k(n)) / (p,v_1,\cdots,v_n),\]
the unit $1$ is $v_{n+1}$-torsion free.
\end{theorem}

\begin{proof}
Consider the $v_{n+1}$-Bockstein spectral sequence of signature
\[\left(\pi_*\gr^*_{\mot} \TC(k(n)) / (p,v_1,\cdots,v_{n+1})\right)[v_{n+1}] \implies \pi_* \gr^*_{\mot} \TC(k(n)) / (p,v_1,\cdots,v_n).\]

We must prove that there is no differential in this spectral sequence killing a power of $v_{n+1}$ times the unit.  For bidegree reasons, the only possible differentials to worry about are $d_1$ differentials by Theorem~\ref{thm:syntomic-cohomology-kn}.  Therefore, if we can show that $v_{n+1}$ times the unit is non-zero in $\pi_* \gr^*_{\mot} \TC(k(n)) / (p,v_1,\
\cdots,v_n)$, it will automatically follow that every power of $v_{n+1}$ times the unit is non-zero.

To check that $v_{n+1}$ times the unit is non-zero, we can check that it has non-zero image under the composite map
\[
\gr^*_{\mot}\TC(k(n))\longrightarrow \gr^*_{\mot}\TC^{-}(k(n))  \longrightarrow \gr^*_{\mot}\mathrm{fiber}(\sigma)
\]
mod $(p,v_1,\cdots ,v_n)$. Here we break with convention and write $\sigma$ for the nonmotivic map 
\[ \sigma : \THH(k(n))\to \Sigma^{-1}\THH(k(n))\,.\]

Since $\mu=\sigma^2v_{n+1}$, the image of $v_{n+1}$ in $\mathrm{fiber}(\sigma)$ agrees with the image of $\mu$ under the map
\[
\Sigma^{-2}\grmot^*\THH(k(n)) / (p,v_1,\cdots,v_{n}) \to \grmot^*\mathrm{fiber}(\sigma)/(p,v_1,\cdots ,v_n)\,.
\]

We note that this cannot map to zero for bidegree reasons, by Theorem~\ref{thm:height-n}.
\end{proof}

\begin{rmk}
Theorem~\ref{thm:unital-syntomic-redshift} implies strong redshift results at the level of syntomic cohomology.  For example, suppose that $A \to k(n)$ is any map of $\mathbb{E}_1$-$\MU$-algebras. Then there is an induced unital map 
\[\pi_* \gr^*_{\mot} \TC(A) / (p,v_1,\cdots,v_n) \to \pi_* \gr^*_{\mot} \TC(k(n)) / (p,v_1,\cdots,v_n),\]
from which it follows that 
\[v_{n+1}^{-1}\gr^*_{\mot}\TC(A) / (p,v_1,\cdots,v_n) \ne 0.\]
\end{rmk}

It is unfortunately not immediate that $L_{T(n+1)} \TC(k(n)) \ne 0$, because of the distinction between syntomic cohomology and topological cyclic homology.  We deal with this distinction using the following proposition, which depends crucially on the purity theorem \cite{LMMT24}.

\begin{proposition}\label{prop:kappa-classes}
There exists a sequence of positive integers $(a_0,\cdots,a_n,a_{n+1})$ and a sequence of classes $(\kappa_0,\kappa_1,\cdots,\kappa_{n-1})$ such that:
\begin{enumerate}
\item \label{it1:redshift} There are maps of homotopy commutative and associative generalized Moore algebras
\[
\mathbb{S} \to \mathbb{S}/(p^{a_0}) \to \mathbb{S}/(p^{a_0},v_1^{a_1}) \to \cdots \to \mathbb{S}/(p^{a_0},v_1^{a_1},\cdots,v_{n+1}^{a_{n+1}})\,,\]
corresponding to a sequence of well-defined even $\MU_*\MU$-comodule algebras
\[
\MU_* \to \MU_* / (p^{a_0}) \to \MU_* / (p^{a_0},v_1^{a_1}) \to \cdots \to \MU_* / (p^{a_0},v_1^{a_1},\cdots,v_{n+1}^{a_{n+1}})\,.
\]
\item \label{it2:redshift} For each $0 \le j \le n-1$, 
\[
\kappa_j \in \pi_{j+(2p-2)a_1+(2p^2-2)a_2+\cdots+(2p^j-2)a_{j}} \TC(k(n))/(p^{a_0},v_1^{a_1},\cdots,v_j^{a_j})\,.
\]
The class $\kappa_0 \in \pi_0 \TC(k(n))/p^{a_0}$ is equal to the unit
\[
\mathbb{S} \to \mathrm{TC}(\mathbb{S}) \to \TC(k(n)) \to \mathrm{TC}(k(n))/p^{a_0} \,.
\]
For each $1 \le j \le n-1$, $\kappa_j$ is carried to $\kappa_{j-1}$ under $\pi_*$ of the Bockstein map
\[
\TC(k(n))/(p^{a_0},v_1^{a_1},\cdots,v_j^{a_j}) \to \Sigma^{(2p^j-2)a_j+1} \TC(k(n)) / (p^{a_0},v_1^{a_1},\cdots,v_{j-1}^{a_{j-1}})\,.
\]
\item \label{it3:redshift} For each $0 \le j \le n-1$, in the motivic spectral sequence 
\[\mathrm{E}_2=\pi_*\gr^*_{\mot} \TC(k(n)) / (p^{a_0},v_1^{a_1},\cdots,v_j^{a_j}) \implies \pi_* \TC(k(n)) / (p^{a_0},v_1^{a_1},\cdots,v_j^{a_j})\,,\]
$\kappa_j$ is detected on the $-j$ line (i.e., in Adams weight $-j$).  The entire $\mathrm{E}_2$-page of this spectral sequence is concentrated on the $-j-1$ line and above. We will denote the class on the $\mathrm{E}_2$-page detecting $\kappa_j$ by $\overline{\kappa_j}$.
\end{enumerate}
\end{proposition}

\begin{proof}
For simplicity, thanks to \cite{burklundMoore}, we will choose a sequence of integers $(a_0,a_1,\cdots,a_{n+1})$ so that each generalized Moore spectrum
\[\mathbb{S}/(p^{a_0},v_1^{a_1},\cdots,v_j^{a_j})\]
is an $\mathbb{E}_{n+3-j}$-algebra for $0\le j\le n+1$.

In particular, we will construct the integers $a_j$ and classes $\kappa_j$ by induction on $j$.  For each $j \ge 0$, we will pick $a_j$ to be a sufficiently large integer, which guarantees the claim in the above paragraph by \cite[Theorem 1.4 and Remark 5.5]{burklundMoore}.  The main point will be to check that, when $0 \le j \le n-1$ and we have already selected $(a_0,a_1,\cdots,a_{j-1})$ and $(\kappa_0,\kappa_1,\cdots,\kappa_{j-1})$, any sufficiently large choice of $a_j$ will allow us to pick $\kappa_j$ satisfying \eqref{it2:redshift} and \eqref{it3:redshift}.

When $j=0$ we can take $a_0$ to be any integer such that $\mathbb{S}/p^{a_0}$ is highly structured, and  $\kappa_0$ to be the composite
\[\mathbb{S} \to \mathrm{TC}(\mathbb{S}) \to \mathrm{TC}(k(n)) \to \TC(k(n))/p^{a_0}\]
It remains to check \eqref{it3:redshift} when $j=0$, which follows from Proposition \ref{prop:integral-collapse} and Corollary \ref{cor:connectivity}.

For $j>0$, we proceed by induction.  First, we observe that the previously constructed $\mathbb{E}_{n+4-j}$ ring $\mathbb{S}/(p^{a_0},v_1^{a_1},\cdots,v_{j-1}^{a_{j-1}})$ admits a self-map acting as a power of $v_j$. By taking sufficiently large powers of that self-map we may find a self-map $v_j^{b_j}$ such that, for every integer $k>0$,
\[\mathbb{S}/(p^{a_0},v_1^{a_1},\cdots,v_{j-1}^{a_{j-1}},v_j^{kb_j})\]
exists as an $\mathbb{E}_{n+3-j}$-algebra, and $\MU_*/(p^{a_0},v_1^{a_1},\cdots,v_{j-1}^{a_{j-1}},v_j^{kb_j})$ is a well-defined $\MU_*\MU$ comodule algebra.  In a moment we will set $a_j=kb_j$ for a specific integer $k$.

Now, by induction we have previously defined a class 
\[\kappa_{j-1} \in \pi_*\TC(k(n)) / (p^{a_0},v_1^{a_1},\cdots,v_{j-1}^{a_{j-1}}).\]
According to the results of \cite{LMMT24}, since $0< j \le n-1$, $L_{T(j)} \TC(k(n)) = 0$.  This means that there exists some integer $k>0$ such that \[v_j^{kb_j} \kappa_{j-1} = 0 \in \pi_*\TC(k(n)) / (p^{a_0},v_1^{a_1},\cdots,v_{j-1}^{a_{j-1}}).\]
We define $a_j$ to be this specific multiple of $b_j$, so that 
\begin{align}\label{eq:relation}
v_j^{a_j} \kappa_{j-1}=0\,.
\end{align}

The relation \eqref{eq:relation} of the previous sentence ensures that we may choose a class
\[\kappa_j \in \pi_*\TC(k(n)) / (p^{a_0},\cdots,v_j^{a_j})\]
which maps to $\kappa_{j-1}$ under the Bockstein, satisfying condition \eqref{it2:redshift}.

It remains to check \eqref{it3:redshift}. First note that the claim that the entire $\mathrm{E}_2$-page is concentrated on the $-j-1$ line and above is immediate, since modding out an element can only introduce one additional line in the spectral sequence.  The remaining claim in \eqref{it3:redshift} follows from checking that $\overline{\kappa_{j-1}}$ is annihilated by $v_j^{a_j}$ already at the $\mathrm{E}_2$-page of the motivic spectral sequence
\[\pi_* \gr^*_{\mot} \TC(k(n)) / (p^{a_0},v_1^{a_1},\cdots,v_{j-1}^{a_{j-1}} ) \implies \pi_* \TC(k(n)) / (p^{a_0},v_1^{a_1},\cdots,v_{j-1}^{a_{j-1}}).\]
This is immediate from the fact that $\overline{\kappa_{j-1}}$ is a permanent cycle exactly one line above the lowest line of the $\mathrm{E}_2$-page and the $d_r$ differentials in the motivic spectral sequence raise Adams weight by exactly $r$. Consequently, no $v_{j}^{b_j}$-power multiple of $\overline{\kappa_{j-1}}$ can be the target of a motivic differential and so the relation $v_j^{a_j} \overline{\kappa}_{j-1} =0 $ that must be present at the $\mathrm{E}_{\infty}$-page is also present at the $\mathrm{E}_2$-page. 
\end{proof}

\begin{theorem}[Redshift for Morava $K$-theory]
The localization $L_{T(n+1)} \K(k(n))$ is not equal to zero for all integers $n\ge 1$ and primes $p$. 
\end{theorem}

\begin{proof}
Choose $(a_0,\cdots,a_n,a_{n+1})$ according to Proposition~\ref{prop:kappa-classes}, and consider the map of motivic spectral sequences

\[
\begin{tikzcd}
    \text{E}_2=\pi_* \gr^*_{\mot} \TC(k(n)) / (p^{a_0},v_1^{a_1},\cdots,v_{n-1}^{a_{n-1}}) \arrow[r,Rightarrow] \arrow{d} & \pi_*\TC(k(n)) / (p^{a_0},v_1^{a_1},\cdots,v_{n-1}^{a_{n-1}}) \arrow{d} \\
    \text{E}_2=\pi_* \gr^*_{\mot} \TC(k(n)) / (p^{a_0},v_1^{a_1},\cdots,v_n^{a_n}) \arrow[r,Rightarrow] & \pi_*\TC(k(n)) / (p^{a_0},v_1^{a_1},\cdots,v_n^{a_n})
\end{tikzcd}
\]

The previous proposition identifies a class $\kappa_{n-1}$ in the upper right of the above square, which is detected by a class $\overline{\kappa_{n-1}}$ in the upper left of Adams weight $-n+1$.  Let $\kappa$ and $\overline{\kappa}$ denote the images of $\kappa_{n-1}$ and $\overline{\kappa_{n-1}}$, respectively, under the vertical maps.

In order to prove the theorem, it suffices to prove that $\kappa$ is $v_{n+1}^{a_{n+1}}$-torsionfree, by which we mean that for no integer $j \ge 0$ is $v_{n+1}^{j a_{n+1}}\kappa$ equal to zero.  For this, it suffices to prove that $\overline{\kappa}$ is $v_{n+1}^{a_{n+1}}$-torsion free on the $\mathrm{E}_{\infty}$-page of the motivic spectral sequence. Applying Theorem~\ref{thm:syntomic-cohomology-kn}, note that
\[
\pi_* \gr^*_{\mot} \TC(k(n)) / (p^{a_0},v_1^{a_1},\cdots,v_n^{a_n})
\]
is concentrated in Adams weight $-n$ and above, because of its finite filtration with associated graded copies of $\pi_* \gr^{*}_{\mot} \TC(k(n)) / (p,v_1,\cdots,v_n)$ and the $v_{n+1}$-Bockstein spectral sequence.  Since motivic $d_r$ differentials raise Adams weight by $r$, it follows that $\overline{\kappa}$ is $v_{n+1}^{a_{n+1}}$-torsion free on the $\mathrm{E}_{\infty}$-page of the motivic spectral sequence if and only if it is $v_{n+1}^{a_{n+1}}$-torsion free on the $\mathrm{E}_2$-page of the motivic spectral sequence.

It remains only to check that $\kappa$ is $v_{n+1}^{a_{n+1}}$-torsion free on the $\mathrm{E}_2$-page of the motivic spectral sequence.  For this, consider the composite map (with some bigraded suspensions ommitted)
\[\gr^*_{\mot} \TC(k(n)) / (p^{a_0},v_1^{a_1},\cdots,v_{n-1}^{a_{n-1}}) \to \gr^*_{\mot} \TC(k(n))/(p^{a_0}) \to \gr^*_{\mot} \TC(k(n)) / (p,\cdots,v_n)\]
of $\gr_{\ev}^*\bS$-modules, 
where the first map in the composite is the iterated Bockstein and the second is projection on to the quotient. 
By Lemma~\ref{lem:redshift-final} below, it then suffices to check that the unit is $v_{n+1}$-torsion free in $\pi_*\gr^*_{\mot} \TC(k(n)) / (p,v_1,\cdots,v_{n})$, which is exactly Theorem~\ref{thm:unital-syntomic-redshift}.

\end{proof}

\begin{lem} \label{lem:redshift-final}
Let $(a_0,a_1,\cdots ,a_{n+1})$  satisfy \eqref{it1:redshift} of Proposition~\ref{prop:kappa-classes}. 
Suppose that $X$ and $Y$ are two $\gr^*_{\ev} \mathbb{S}$-modules equipped with a $\gr^*_{\ev} \mathbb{S}$-module map
\[f:X / (p^{a_0},v_1^{a_1},\cdots,v_{n-1}^{a_{n-1}}) \to Y / (p,v_1,\cdots,v_n).\]
Furthermore suppose that $\kappa$ is a class in $\pi_* \left(X  /(p^{a_0},v_1^{a_1},\cdots,v_{n-1}^{a_{n-1}})\right)$ such that $(\pi_*f)(\kappa)$ is $v_{n+1}$-torsion free.  Then the image of $\kappa$ under the quotient map
\[\pi_* \left(X  /(p^{a_0},v_1^{a_1},\cdots,v_{n-1}^{a_{n-1}}) \right)\to \pi_* \left(X  /(p^{a_0},v_1^{a_1},\cdots,v_{n-1}^{a_{n-1}},v_n^{a_n})\right)\]
is $v_{n+1}^{a_{n+1}}$-torsion free.

\end{lem}

\begin{proof}
The target $Y/(p,v_1,\cdots,v_n)$ is by definition equal to $Y \otimes_{\gr^*_{\ev} \mathbb{S}} (\gr^*_{\ev} \mathbb{S} / (p,v_1,\cdots,v_n))$, as a module over the $\mathbb{E}_{\infty}$-$\gr^*_{\ev} \mathbb{S}$-algebra $\gr^*_{\ev} \mathbb{S} / (p,v_1,\cdots,v_n)$.  Thus, the map $f$ factors through the projection 
\[
\begin{tikzcd}
X \otimes_{\gr^*_{\ev} \mathbb{S} } (\gr^*_{\ev} \mathbb{S}/ (p^{a_0},v_1^{a_1}, \cdots, v_{n-1}^{a_{n-1}})) \arrow{d}{\mathrm{id}\otimes \eta_L} \\
X \otimes_{\gr^*_{\ev} \mathbb{S} } (\gr^*_{\ev} \mathbb{S}/ (p^{a_0},v_1^{a_1}, \cdots, v_{n-1}^{a_{n-1}})) \otimes_{\gr^*_{\ev} \mathbb{S} } (\gr^*_{\ev} \mathbb{S}/ (p,v_1, \cdots, v_n)) \,.
\end{tikzcd}
\]
Let us thus consider 
\[(\gr^*_{\ev} \mathbb{S}/ (p^{a_0},v_1^{a_1}, \cdots, v_{n-1}^{a_{n-1}})) \otimes_{\gr^*_{\ev} \mathbb{S} } (\gr^*_{\ev} \mathbb{S}/ (p,v_1, \cdots,v_n))\]
together with its two unit maps
\[\eta_L:\gr^*_{\ev} \mathbb{S}/ (p^{a_0},v_1^{a_1}, \cdots, v_{n-1}^{a_{n-1}}) \to (\gr^*_{\ev} \mathbb{S}/ (p^{a_0},v_1^{a_1}, \cdots, v_{n-1}^{a_{n-1}})) \otimes_{\gr^*_{\ev} \mathbb{S} } (\gr^*_{\ev} \mathbb{S}/ (p,v_1, \cdots,v_n))\]
and
\[\eta_R:\gr^*_{\ev} \mathbb{S}/ (p,v_1, \cdots, v_n) \to (\gr^*_{\ev} \mathbb{S}/ (p^{a_0},v_1^{a_1}, \cdots, v_{n-1}^{a_{n-1}})) \otimes_{\gr^*_{\ev} \mathbb{S} } (\gr^*_{\ev} \mathbb{S}/ (p,v_1, \cdots, v_n))\,.
\]

First, we claim that $(\pi_*\eta_L) (v_n^{a_n})$ is nilpotent, so that there exists some $k$ such that $\eta_L$ factors over $\gr^*_{\ev} \mathbb{S} / (p^{a_0},v_1^{a_1},\cdots,v_{n-1}^{a_{n-1}},v_n^{ka_n})$.  Second, we claim that some power of $v_{n+1}^{a_{n+1}}$ in $\pi_*\gr^*_{\ev} \mathbb{S} / (p^{a_0},v_1^{a_1},\cdots,v_n^{a_n})$ lifts to $\pi_*\gr^*_{\ev} \mathbb{S} / (p^{a_0},v_1^{a_1},\cdots,v_{n-1}^{a_{n-1}},v_n^{ka_n})$, such that its image in the homotopy of 
\[ (\gr^*_{\ev} \mathbb{S}/ (p^{a_0},v_1^{a_1}, \cdots, v_{n-1}^{a_{n-1}})) \otimes_{\gr^*_{\ev} \mathbb{S} } (\gr^*_{\ev} \mathbb{S}/ (p,v_1, \cdots, v_n))\]
agrees with a power of $(\pi_* \eta_R)(v_{n+1})$.  The lemma statement immediately follows.

To see the claims, note that given any sequence of integers $(b_0,\cdots,b_j)$ such that \[\MU_* / (p^{b_0},v_1^{b_1},\cdots,v_j^{b_j})\] is a well-defined $\MU_*\MU$-comodule algebra, classes in the kernel of the multiplication map
\[\gr^*_{\ev} \mathbb{S} / (p^{b_0},v_1^{b_1},\cdots,v_j^{b_j}) \otimes_{\gr^*_{\ev} \mathbb{S}} \gr^*_{\ev} \mathbb{S} / (p^{b_0},v_1^{b_1},\cdots,v_j^{b_j}) \to \gr^*_{\ev} \mathbb{S} / (p^{b_0},v_1^{b_1},\cdots,v_j^{b_j})\]
are nilpotent, because the associated Adams spectral sequence has a horizontal vanishing line (the domain of the multiplication map is in the thick subcategory generated by the codomain).
\end{proof}

\begin{cor}
The localization $L_{T(n+1)} \K(K(n))$ is not equal to zero for all integers $n\ge 1$ and primes $p$.
\end{cor}

\begin{proof}
For $n\ge 1$, we know $L_{T(n+1)}\K(\bF_p)\simeq 0$ by~\cite{Qui72} so this follows from the localization sequence 
\[ \K(\bF_p)\to\K(k(n))\to \K(K(n))
\]
of~\cite{BM08} (cf.~\cite[8.8.4]{BarwickHeart}). 
\end{proof}

\subsection{Morava K-theories over a perfect field}\label{sec:perfect-fields}
Let $k$ be a perfect field of characteristic $p$, with ring of Witt vectors $\W(k)$. 
Then we might consider the base change $k(n) \otimes \bS_{\W(k)}$, where $\bS_{\W(k)}$ denotes the spherical Witt vectors from~\cite[Example~5.2.7]{Lur21}.

\begin{remark}~\label{alg-closed-kn}
Choose a perfect field $k$ and height $n$ formal group law $\bG$ over $k$, so that there is an associated Lubin--Tate $\mathbb{E}_\infty$-ring $E(k,\mathbb{G})$.  Then, a very natural project is to study the algebraic $K$-theories of $\mathbb{E}_1$-$E(k,\mathbb{G})$-algebras $K_n=E(k,\mathbb{G})/(p,u_1,\cdots,u_{n-1})$, where $\pi_*(K_n) \cong k[u^{\pm}]$, $|u|=2$.  

In the case where $k$ is algebraically closed, and $k(n)$ is an $\mathbb{E}_1$-$\mathrm{BP}\langle n \rangle$-algebra, then there exists a ring structure on $K_n$ such that $\mathrm{K}(K_n)$ has $\mathrm{K}(k(n) \otimes \mathbb{S}_{W(k)})$ as a retract.  For this particular ring structure, redshift for $\mathrm{K}(k(n) \otimes \mathbb{S}_{W(k)})$ implies redshift for $\mathrm{K}(K_n)$ and thus for $\mathrm{K}(E(k,\mathbb{G}))$ \cite[Corollary 9.11]{ABM22}. 
\end{remark}

In order to understand the algebraic $K$-theory of $k(n) \otimes \mathbb{S}_{W(k)}$, we can replace $\MU$ with $\MU_{\W(k)}:=\MU\otimes \bS_{\W(k)}$.


In this section, we write $\gr^*_{\mot/\MU_{W(k)}}$ 
to denote the functor defined analogously to $\grmot^*$, but with each instance of $\MU$ replaced by $\MU_{\W(k)}$. Note that one can also replace $\mathrm{MW}$ with $\mathrm{MW}_{\W(k)}:=\mathrm{MW}\otimes \bS_{\W(k)}$ in order to produce an eff map of cyclotomic $\mathbb{E}_{\infty}$-rings $\THH(\MU_{\W(k)})\to \THH(\MU_{\W(k)}/\MW_{\W(k)})$ whose target is even. All of the results for $\THH$, $\THH^{tC_p}$, $\TC^{-}$, and $\TP$ can be computed mutatis mutandis, providing the following theorem. 
\begin{thm}
There is an isomorphism of bigraded $k$-vector spaces   
\[ 
\pi_*\gr^*_{\mot/\MU_{W(k)}}F(k(n)\otimes \mathbb{S}_{W(k)})/(p,v_1,\cdots ,v_n)\cong \pi_*(\grmot^*F(\mathrm{k}(n)))/(p,v_1,\cdots ,v_n)\otimes_{\mathbb{F}_p} k
\]
for $F\in \{\THH, \THH^{tC_p},\TC^{-},\TP\}$. 
\end{thm}
Some care must be taken when considering topological cyclic homology. Instead of simply applying base-change, we have the following result. 
\begin{thm}\label{thm:syntomic-for-k}
Let $k$ be a perfect field of characteristic $p$. Then there is an isomorphism of bigraded $\bF_p$-vector spaces 
\begin{align*}\label{eq: syntomic over k} \pi_*\gr_{\mot/\MU_{\W(k)}}^*\TC(k(n)\otimes \mathbb{S}_{W(k)})/(p,v_1,\cdots ,v_n)\cong& \Lambda(\bar{\varepsilon}_{1},\cdots  ,\bar{\varepsilon}_{n})\otimes \Lambda (\lambda_{n+1})\oplus\\
  &k_{\textup{Fr}_p} \otimes_{\bF_p}\Lambda(\bar{\varepsilon}_{1},\cdots  \bar{\varepsilon}_{n})\otimes \Lambda (\lambda_{n+1})\{\partial\}\oplus  \\
 &\bigoplus_{S\subset \{1,\cdots ,n\}} k\otimes_{\mathbb{F}_{p}}M_S 
\end{align*}
where $M_{S}$ is defined in Definition~\ref{MS}. 
In particular, the only non-trivial group in Adams weights $\le -n$ is $\bF_p\{\bar{\varepsilon}_1\cdots\bar{\varepsilon}_n\}$ in bidegree 
\[ 
(\sum_{i=1}^n2p^i-n,-n) \,.
\] 
\end{thm}
 Here, we write $k_{\textup{Fr}_p}$ for the coinvariants of the action of Frobenius on $k$. 
 
\begin{proof}
The proof is the same except that, on the summand $\Lambda(\varepsilon_1,\cdots ,\varepsilon_n)\otimes k$, the map $\can-\varphi$ is given by $1-\Fr_p$. This has kernel $k^{\textup{Fr}_p}\otimes_{\mathbb{F}_p}\Lambda(\varepsilon_1,\cdots ,\varepsilon_n)=\Lambda(\varepsilon_1,\cdots ,\varepsilon_n)$  and cokernel $k_{\textup{Fr}_p}\otimes_{\mathbb{F}_p}\Lambda(\varepsilon_1,\cdots ,\varepsilon_n)$
\end{proof}

\begin{cor}
The mod $(p,v_{1},\cdots ,v_{n+1})$ $\MU_{W(\overline{\mathbb{F}}_p)}$-based syntomic cohomology of the $\bE_1$-$\MU_{W(\overline{\mathbb{F}}_p)}$-algebra $k(n)\otimes \bS_{W(\bar{\bF}_p)}$ is 
\begin{align*}
\pi_*\gr^*_{\mathrm{mot}/\MU_{W(\overline{\mathbb{F}}_p)}}\TC(\mathrm{k}_{\bar{\mathbb{F}}_{p}}(n))/(p,v_1,\cdots ,v_{n+1})\cong & \bar{\bF}_{p}\otimes_{\bF_p}\Lambda(\bar{\varepsilon}_{1},\cdots  \bar{\varepsilon}_{n})\otimes \Lambda (\lambda_{n+1})\oplus \\
 & \bigoplus_{S\subset \{1,\cdots ,n\}} \bar{\mathbb{F}}_{p}\otimes_{\bF_{p}}M_S 
\end{align*}
where $M_{S}$ is defined in Definition~\ref{MS}. 
\end{cor}

\begin{cor}
The mod $(p,v_{1},v_{2})$ $\MU_{W(\overline{\mathbb{F}}_p)}$-based  syntomic cohomology of the $\bE_1$-$\MU_{W(\overline{\mathbb{F}}_p)}$-algebra $k(2)\otimes \mathbb{S}_{W(\overline{\mathbb{F}}_p)}$ is 
\begin{equation*}\label{eq:syntomic}
\bF_{p}[v_{3}] \otimes \left ( \Lambda(\bar{\varepsilon}_{1}, \bar{\varepsilon}_{2})\otimes \Lambda (\lambda_{3})\oplus \bigoplus_{S\subset \{1,2\}} \bar{\bF}_{p}\otimes_{\bF_{p}} M_S \right )
\end{equation*}
where $M_{S}$ is defined in Definition~\ref{MS}. 
\end{cor}

\begin{proof}
From Figure~\ref{fig:alg-closure}, we observe that the only possible $v_{3}$-Bockstein spectral sequence differentials are 
\begin{align}
\label{first d1} d_{1}(t^{p^{2}+p-1}\bar{\varepsilon}_{1}\bar{\varepsilon}_{2}\lambda_{3}) & \in \bF_{p}\{v_{3}\} \\
\label{second d1}d_{1}(\bar{\varepsilon}_{1} \bar{\varepsilon}_{2} \lambda_{3}) & \in \bF_{p}\{v_{3}t^{p^{3}-p^{2}}\bar{\varepsilon}_{1}\lambda_{3},v_{3}t^{p^{3}-p}\bar{\varepsilon}_{2}\lambda_{3}\} \\
\label{third d1} d_{1}(\bar{\varepsilon}_{2}\lambda_{3}) & \in \bF_{p}\{v_{3}t^{p^{3}-p^{2}}\lambda_{3}\} 
\end{align}
and since in each case both source and target are in $\ker (\can - \varphi)$, we can consider the same differentials in the motivic spectral sequence 
\begin{equation}\label{eq:TC-Bockstein}
\pi_{*}\grmot^{*}\TC^{-}(k(n)\otimes \mathbb{S}_{W(\bar{\bF}_p)})/(p,v_{1},v_{2},v_{3})[v_{3}]\implies  \pi_{*}\grmot^{*}\TC^{-}(k(n)\otimes \mathbb{S}_{W(\bar{\mathbb{F}}_p)})/(p,v_{1},v_{2})\,.
\end{equation}
The differential \eqref{first d1} can be ruled out by considering the map of $v_{3}$-Bockstein spectral sequences induced by the canonical map since $\text{can}(t^{p^{2}+p-1}\bar{\varepsilon}_{1}\bar{\varepsilon}_{2}\lambda_{3})=0$ and $\text{can}(v_{3})=v_{3}$. The differentials in the Bockstein spectral sequence \eqref{eq:TC-Bockstein} are $\lambda_{3}$-linear because $\lambda_{3}$ is a permanent cycle in the corresponding spectral sequence for $\BP$ and the spectral sequence \eqref{eq:TC-Bockstein} is a module over the corresponding spectral sequence for $\BP$. Therefore, if any of the differentials \eqref{second d1} or \eqref{third d1} occur then this would imply that one of the differentials 
\begin{align}
\label{second d1 prime}d_{1}(\bar{\varepsilon}_{1}\bar{\varepsilon}_{2}) & \in \bF_{p}\{v_{3}t^{p^{3}-p^{2}}\bar{\varepsilon}_{1},v_{3}t^{p^{3}-p}\bar{\varepsilon}_{2}\lambda_{3}\} \\
\label{third d1 prime}d_{1}(\bar{\varepsilon}_{2}) & \in \bF_{p}\{v_{3}t^{p^{3}-p^{2}}\} 
\end{align}
and now these latter non-trivial differentials can be ruled out by mapping to $\mathbb{F}_{p}$, where we know that these $d_{1}(\bar{\varepsilon}_{1}\bar{\varepsilon}_{2})=d_{1}(\bar{\varepsilon}_{2})=0$ in the $v_{3}$-Bockstein spectral sequence. 
\end{proof}

By applying the same proof strategy as the previous section and using Theorem~\ref{thm:syntomic-for-k}, we produce the following redshift result. 

\begin{thm}\label{thm:Morava K-theory}
Let $k$ be a perfect field of characteristic $p$. Then $\K(k(n)\otimes \mathbb{S}_{W(k)})$ has height exactly $n+1$ and $\K(K(n)\otimes \mathbb{S}_{W(k)})$ has height exactly $n+1$. 
\end{thm}

\begin{remark}
By Remark \ref{alg-closed-kn}, the above result provides a proof that $E(k,\mathbb{G})$ satisfies redshift that is independent of~\cite[Theorem~A]{Yua21}. It is also independent of (though closely related to) \cite{HW22}, so differs from the proof given in~\cite[Theorem~8.9]{ABM22}.

This is of interest because redshift for $E(k,\mathbb{G})$ is one of the lynchpins in the proof of redshift for all $\bE_{\infty}$-rings by~\cite{BSY22}. We wonder whether Morava $K$-theories, or similar $\mathbb{E}_1$-rings, are algebras over a large class of $\mathbb{E}_2$-rings. 
\end{remark}

\begin{remark}
The fact that the mod $(p,v_1,\cdots,v_{n+1})$ $\MU_{W(\overline{\mathbb{F}}_p)}$ syntomic cohomology of $k(n) \otimes \mathbb{S}_{W(\overline{\mathbb{F}}_p)}$ is concentrated in Adams weights $1$ and below formally implies that the same is true of the integral $\MU_{W(\overline{\mathbb{F}}_p)}$ syntomic cohomology, through Bockstein spectral sequences.  Thus, the motivic spectral sequence computing $\pi_*\TC(k(n) \otimes \mathbb{S}_{W(\overline{\mathbb{F}}_p)})_p^{\wedge}$ is concentrated on only the $0$-line and the $1$-line.
\end{remark}

\begin{figure}[ht!]
\resizebox{\textwidth}{!}{ 
\begin{tikzpicture}[radius=1,yscale=2]
\foreach \n in {-2,-1,...,26} \node [below] at (\n,-.8-3) {$\n$};
\foreach \s in {-3,-2,...,3} \node [left] at (-.3-2,\s) {$\s$};
\draw [thin,color=lightgray] (-2,-3) grid (26,3);
\node [draw,minimum size=.1cm,below] at (0,0) {$1$};
\node [draw,minimum size=.1cm,below] at (3,-1) {$\bar{\varepsilon}_1$};
\node [draw,minimum size=.1cm,below] at (7,-1) {$\bar{\varepsilon}_2$};
\node [draw,minimum size=.1cm,below] at (10,-2) {$\bar{\varepsilon}_1\bar{\varepsilon}_2$};
\node [draw,minimum size=.1cm,above] at (15,1) {$\lambda_3$};
\node [above] at (13,1) {$t\lambda_3$};
\node [above] at (11,1) {$t^2\lambda_3$};
\node [above] at (9,1) {$t^3\lambda_3$};
\node [above] at (7,1) {$t^4\lambda_3$};
\node [below] at (14,0) {$t^2\bar{\varepsilon}_1\lambda_3$};
\node [below] at (12,0) {$t^3\bar{\varepsilon}_1\lambda_3$};
\node [below] at (10,0) {$t^4\bar{\varepsilon}_1\lambda_3$};
\node [below] at (8,0) {$t^5\bar{\varepsilon}_1\lambda_3$};
\node [above] at (16,0) {$t^3\bar{\varepsilon}_2\lambda_3$};
\node [above] at (14,0) {$t^4\bar{\varepsilon}_2\lambda_3$};
\node [above] at (12,0) {$t^5\bar{\varepsilon}_2\lambda_3$};
\node [above] at (10,0) {$t^6\bar{\varepsilon}_2\lambda_3$};
\node [below] at (17,-1) {$t^4\bar{\varepsilon}_1\bar{\varepsilon}_2\lambda_3$};
\node [below] at (15,-1) {$t^5\bar{\varepsilon}_1\bar{\varepsilon}_2\lambda_3$};
\node [below] at (13,-1) {$t^6\bar{\varepsilon}_1\bar{\varepsilon}_2\lambda_3$};
\node [below] at (11,-1) {$t^7\bar{\varepsilon}_1\bar{\varepsilon}_2\lambda_3$};
\node [draw,minimum size=.1cm,above] at (18,0) {$\bar{\varepsilon}_1\lambda_3$};
\node [draw,minimum size=.1cm,above] at (22,0) {$\bar{\varepsilon}_2\lambda_3$};
\node [draw,minimum size=.1cm,above] at (25,-1) {$\bar{\varepsilon}_1\bar{\varepsilon}_2\lambda_3$};
\end{tikzpicture}
}
\caption{The mod $(2,v_1,v_2,v_3)$-syntomic cohomology of $k_{\bar{\bF}_{2}}(2)$. Here boxed classes are generators of $\bF_2$ and unboxed classes are generators of $\bar{\bF}_2$}
\label{fig:alg-closure}
\end{figure}


\appendix
\section{Forms of Morava K-theory}\label{forms}
\setcounter{subsection}{1}
In this appendix, we make a few remarks about the space of choices of the $\bE_1$-$\MU$-algebra $K(n)$ that is fixed at the beginning of the paper.  We begin by discussing the less structured notion of an $\bE_1$-$\bS$-algebra form of $K(n)$.

\begin{defin}
Suppose $A$ is an $\bE_1$-ring with 
\[ 
\pi_*A\cong \mathbb{F}_p[v_n^{\pm 1}],
\]
where $|v_n|=2p^n-2$.  If the height $n$ telescopic localization $L_{T(n)}A \ne 0$, then we refer to $A$ as an $\bE_1$-$\bS$-algebra form of $K(n)$. In this case, we also call $\tau_{\ge 0}A$ an $\bE_1$-$\bS$-algebra form of $k(n)$. 
\end{defin}

\begin{remark}
If $A$ is any $\bE_1$-$\bS$-algebra form of $K(n)$, then every $A$-module is isomorphic to a direct sum of suspensions of $A$.  This follows from the fact that $\pi_*(A)$ is a graded field, i.e. that every graded module over it is free.
\end{remark}

\begin{remark} \label{rmk:K(n)-underlying-spectrum}
The underlying spectra of any two $\bE_1$-$\bS$-algebra forms of $K(n)$ are equivalent.  Indeed, if $A_1$ and $A_2$ are any two $\bE_1$-$\bS$-algebra forms of $K(n)$, then the fact that $A_1$ and $A_2$ are both non-trivial after $T(n)$-localization implies that $A_1 \otimes A_2 \ne 0$.  Since $A_1 \otimes A_2$ is a non-zero free $A_1$-module, we see that $A_1$ is a unital retract of $A_1 \otimes A_2$, and thus that $A_1$ admits the structure of a homotopy $A_2$ module.  In particular, $A_1$ is a direct sum of shifts of $A_2$, and comparing homotopy groups we learn that $A_1 \cong A_2$.  
\end{remark}

\begin{remark}
Given $A$ an $\bE_1$-$\bS$-algebra form of $K(n)$, the associated form of $k(n)$ is $\tau_{\ge 0}A$.  One can also recover $A$ as an $\bE_1$-ring from $\tau_{\ge 0} A$, by $T(n)$-localization.
\end{remark}

By~\cite[Corollary~1.3]{CM15}, the Quillen idempotent defines a map
\[
\BP \longrightarrow \MU_{(p)}
\]
of $\bE_2$-rings. 
While $\pi_*\BP$ is only noncanonically isomorphic to a polynomial ring 
\[
\bZ_{(p)}[v_1,v_2,\cdots]\,,
\] 
for each $i \ge 1$ the subring 
\[
\bZ_{(p)}[v_1,\cdots,v_i] \subset \pi_*\BP
\] 
is well-defined. Indeed, it is the subring generated by all elements of degree at most $2p^i-2$. 

\begin{defin}\label{def:form-kn}
An \emph{$\bE_1$-$\MU$-algebra form of $k(n)$} is a $p$-local $\bE_1$-$\MU$-algebra~$A$ such that, along the composite ring map
\[
\bZ_{(p)}[v_1,\cdots,v_n] \subset \pi_*\BP \subset \pi_*\MU_{(p)} \to \pi_*A,
\]
\begin{itemize}
\item The element~$p$, as well as all elements of degree strictly between $0$ and $2p^n-2$, map to zero.
\item The induced map~$\bF_p[v_n] \to \pi_*A$ is an isomorphism.
\end{itemize}
Here, the map~$\pi_*\MU_{(p)} \to \pi_*A$ arises as the $p$-localization of the $\MU$-algebra unit map.
\end{defin}

In this paper, we are interested in the algebraic K-theories of $\bE_1$-$\MU$-algebra forms of $k(n)$.  Of course, algebraic K-theory depends only on underlying $\bE_1$-$\bS$-algebra structure, but we exploit $\bE_1$-$\MU$-algebra structure when making our computations. 

The following proposition essentially follows from work of Angeltveit~\cite{AngeltveitUniqueness}.  It implies that we are not giving anything up by assuming that an $\bE_1$-$\MU$-algebra structure exists:

\begin{prop}\label{prop:MU-alg-vs-S-alg-on-kn}
Suppose $A$ is an $\bE_1$-$\bS$-algebra form of $k(n)$. Then $A$ is the underlying $\bE_1$-ring of an $\bE_1$-$\MU$-algebra form of $k(n)$.
\end{prop}

\begin{proof}
Fix an $\bE_1$-$\bS$-algebra form $A$ of $k(n)$. Let $P_m A$ denote $\tau_{\le m(2p^n-2)} A$, which canonically inherits an $\bE_1$-ring structure.  Note that $P_0 R \simeq \bF_p$ is an Eilenberg--MacLane spectrum, which admits unique $\bE_1$-ring and $\bE_1$-$\BP$-algebra structures.

Assuming now that we have chosen an $\bE_1$-$\MU$-algebra structure on $P_{m} A$, compatible with the fixed $\bE_1$-ring structure, we will produce a compatible $\bE_1$-$\MU$-algebra structure on $P_{m+1}A$.  Taking the limit in $m$, we will recover that $A$ admits an $\bE_1$-$\MU$-algebra structure.

The unit map
$\BP \to \MU_{(p)} \to P_mA$ is a map of spectra, and the spectrum $P_mA$ is determined by Remark \ref{rmk:K(n)-underlying-spectrum}. Comparing spectrum $k$-invariants of $P_m A$ and $\BP$, we see that the composite
\[
\pi_* \BP \to \pi_*\MU_{(p)} \to \pi_*P_mA
\]
is surjective, and the map $\pi_* \MU_{(p)} \to \pi_*P_mA$ is given by modding out a regular sequence of elements in $\pi_* \MU_{(p)}$.  Thus, exactly as in~\cite[Proposition 5.7]{AngeltveitUniqueness}, we may calculate the $\MU_{(p)}$-based topological Hochschild cohomology of $P_mA$ (with $\bF_p$ coefficients) using the universal coefficient spectral sequence~\cite[Corollary~2.5]{AHL10}. As in \cite[Theorem 5.9]{AngeltveitUniqueness} we may similarly calculate the $\bS$-based topological Hochschild cohomology of $P_mA$ using the map of the universal coefficient spectral sequences induced by $P_mA\to \bF_p$.

 In particular, we observe that the map 
\[ 
\THC^{(m+1)(2p-2)+2}(P_{m}A/\MU;\bF_p) \longrightarrow \THC^{(m+1)(2p-2)+2}(P_{m}A;\bF_p) 
\]
induced by the unit map $\bS\to \MU$ is surjective. Thus, the class representing the $\bE_1$-$\bS$-algebra $k$-invariant lifts to the class representing the  $\bE_1$-$\MU$-algebra $k$-invariant in this degree.

Passing to the limit over $m$ produces an $\bE_{1}$-$\MU$-algebra structure on $\tau_{\ge 0}A$. It remains to check that $\tau_{\ge 0} A$ is an $\bE_1$-$\MU$-algebra form of $k(n)$.  For this, we must consider the composite map
$\pi_*\BP \to \pi_*\MU_{(p)} \to \pi_*A$.  For degree reasons this sends elements of degree strictly between $0$ and $2p^n-2$ in $\bZ_{(p)}[v_1,\cdots,v_{n-1}]$ to zero, and since it is a ring map it sends $p$ to zero. Comparing spectrum $k$-invariants, the class $v_n$ is sent to a non-zero class in $\pi_*R$.
\end{proof}

\begin{rmk}
In~\cite{AngeltveitUniqueness}, Angeltveit proves that any two $\bE_1$-$\MU$-algebra forms of $k(n)$ are equivalent as $\bE_1$-$\bS$-algebras, \emph{so long as the two forms of $k(n)$ have formal group laws isomorphic to the Honda formal group law}. We expect that, in general, Angeltveit's arguments prove that two $\bE_1$-$\MU$-algebra forms of $k(n)$ are equivalent as $\bE_1$-$\bS$-algebras if and only if they have isomorphic formal group laws.  It would be nice to revisit this result.  Finally, if $\bG$ is a height~$n$ formal group over a perfect field $k$, it would be wonderful to classify $\bE_1$-$\bS_{\W(k)}$-algebra structures on $2$-periodic Morava K-theories under the associated Lubin--Tate theory. Here $\bS_{\W(k)}$ denotes the spherical Witt vectors from~\cite[Example~5.2.7]{Lur21}. 
\end{rmk}

In the main body of this paper, we fix a particular (but arbitrary) $\bE_1$-$\MU$-algebra form of $k(n)$, according to the convention below.  Since the choices are arbitrary, the theorems we prove about $k(n)$ hold for all $\bE_1$-$\MU$-algebra forms. In light of Proposition~\ref{prop:MU-alg-vs-S-alg-on-kn}, the main results in this paper therefore hold for all $\bE_1$-$\bS$-algebra forms.

\begin{convention}\label{conv:vi}
Throughout this paper, we use the symbol $k(n)$ to denote a fixed $\bE_1$-$\MU$-algebra form of $k(n)$. We also fix for each $i \ge 1$ an indecomposable polynomial generator 
\[
v_i \in \pi_{2p^i-2} \BP \,,
\]
and denote also by $v_i$ the image of this class under the map $\pi_{2p^i-2} \BP \to \pi_{2p^i-2} \MU_{(p)}$. We make these choices such that, for each $i>n$, $v_i$ maps to zero in $\pi_{2p^i-2} k(n)$.
\end{convention}

\bibliographystyle{alpha}
\bibliography{kkn}
\end{document}